\documentclass[]{article}

\usepackage{xr}
\usepackage[doi=false,isbn=false,url=false,
	hyperref=auto,
	sorting=nyt,
	maxnames=10,
	maxcitenames=3,
	backend=biber,
	texencoding=auto,
	giveninits=true,
	block=space,
	style=numeric,
]{biblatex}

\DefineBibliographyStrings{english}{%
	andothers = {et al}
} 
 
\DeclareFieldFormat{eprint:arxiv}{%
	\href{https://arxiv.org/abs/#1}{\emph{arXiv:\allowbreak #1}}}

\DeclareFieldFormat{eprint:mrnumber}{%
	\href{http://www.ams.org/mathscinet-getitem?mr=MR#1}{MR\allowbreak #1}}

\DeclareFieldFormat{eprint:jstor}{%
	\href{http://www.jstor.org/stable/#1}{JSTOR\allowbreak #1}}

\DeclareFieldFormat{eprint:customeprint}{%
	\em Available at \upshape \ttfamily \href{https://#1}{#1}}

\DeclareFieldFormat{eprint:online}{%
	Available \href{https://#1}{online}}

\DeclareFieldFormat{eprint:inprep}{%
	\emph{In preparation}}

\DeclareFieldFormat{eprint:manual}{%
	\emph{#1}}

\DeclareFieldFormat{eprint:onarxiv}{%
	\emph{#1}}

\DeclareFieldFormat{eprint:toappear}{%
	\mbox{\emph{To appear}}}

\DeclareFieldFormat{eprint:accepted}{%
	\emph{Accepted at {#1}; to appear}}

\DeclareSourcemap{
	\maps[datatype=bibtex]{
		\map{
			\step[fieldsource=mrnumber,		fieldtarget=eprint, final]
			\step[fieldset=eprinttype,		fieldvalue=mrnumber]
		}
		\map{
			\step[fieldsource=arxiv,		fieldtarget=eprint, final]
			\step[fieldset=eprinttype,		fieldvalue=arxiv]
		}
		\map{
			\step[fieldsource=jstor,		fieldtarget=eprint, final]
			\step[fieldset=eprinttype,		fieldvalue=jstor]
		}
		\map{
			\step[fieldsource=customeprint,	fieldtarget=eprint, final]
			\step[fieldset=eprinttype,		fieldvalue=customeprint]
		}
		\map{
			\step[fieldsource=online,		fieldtarget=eprint, final]
			\step[fieldset=eprinttype,		fieldvalue=online]
		}
		\map{
			\step[fieldsource=inprep,		fieldtarget=eprint, final]
			\step[fieldset=eprinttype,		fieldvalue=inprep]
		}
		\map{
			\step[fieldsource=manual,		fieldtarget=eprint, final]
			\step[fieldset=eprinttype,		fieldvalue=manual]
		}
		\map{
			\step[fieldsource=onarxiv,		fieldtarget=eprint, final]
			\step[fieldset=eprinttype,		fieldvalue=onarxiv]
		}
		\map{
			\step[fieldsource=toappear,		fieldtarget=eprint, final]
			\step[fieldset=eprinttype,		fieldvalue=toappear]
		}
		\map{
			\step[fieldsource=accepted,		fieldtarget=eprint, final]
			\step[fieldset=eprinttype,		fieldvalue=accepted]
		}
	}
}

\DeclareFieldFormat{doi}{%
	\href{https://doi.org/#1}{DOI}%
}

\DeclareFieldFormat{url}{%
	\href{https://#1}{\texttt{#1}}%
}

\DeclareFieldFormat{comments}{%
	\emph{#1}%
} 
 
\DeclareFieldFormat
	[article,inbook,incollection,inproceedings,patent,thesis,unpublished]
	{title}{{#1\isdot}}

\DeclareFieldFormat
	[article,book,inbook,incollection,inproceedings,patent,thesis,unpublished, online]
	{date}{\mkbibparens{#1}}
\DeclareFieldFormat
	[article,book,inbook,incollection,inproceedings,patent,thesis,unpublished]
	{volume}{\mkbibbold{#1}}
\DeclareFieldFormat
	{pages}{\mkbibparens{#1}}

\DeclareFieldFormat{edition}{%
	\ifinteger{#1}
	{\mkbibordedition{#1}~\bibstring{edition}\addcomma}
	{#1~\bibstring{edition}\addcomma}}

\renewbibmacro*{volume+number+eid}{%
	\printfield{volume}%
\setunit*{\adddot}%
	\printfield{number}%
\setunit{\space}%
	\printfield{eid}%
}

\DeclareBibliographyDriver{article}{%
	\usebibmacro{bibindex}%
	\usebibmacro{begentry}%
	\usebibmacro{author}%
\setunit{\addspace}%
	\usebibmacro{date}%
\newunit\newblock
	\usebibmacro{title}%
\newunit\newblock
	\printfield{journaltitle}\addperiod%
\setunit*{\addspace}%
	\usebibmacro{volume+number+eid}%
\setunit*{\addspace}\newblock
	\printfield{pages}
\setunit*{\addspace}
	\printfield{comments}%
	\usebibmacro{eprint}%
\setunit*{\addspace}%
	\printfield{doi}%
	\usebibmacro{finentry}%
}

\DeclareBibliographyDriver{online}{%
	\usebibmacro{bibindex}%
	\usebibmacro{begentry}%
	\usebibmacro{author}%
\setunit{\addspace}%
	\usebibmacro{date}%
\newunit\newblock
	\usebibmacro{title}%
\newunit\newblock
	\usebibmacro{eprint}%
	\usebibmacro{finentry}%
}

\DeclareBibliographyDriver{book}{%
	\usebibmacro{bibindex}%
	\usebibmacro{begentry}%
	\usebibmacro{author}%
\setunit{\addspace}\newblock
	\usebibmacro{date}%
\newunit\newblock
	\usebibmacro{title}%
\setunit*{\addspace}\newblock
	\printfield{edition}%
\newunit\newblock
	\printlist{publisher}
\setunit*{\addspace}%
	\printfield{volume}%
\setunit*{\addspace}\newblock%
	\printfield{comments}%
	\usebibmacro{eprint}%
\setunit*{\addspace}
	\printfield{doi}
	\usebibmacro{finentry}%
}

\DeclareBibliographyDriver{thesis}{%
	\usebibmacro{bibindex}%
	\usebibmacro{begentry}%
	\usebibmacro{author}%
\setunit{\addspace}\newblock
	\usebibmacro{date}%
\newunit\newblock
	\usebibmacro{title}.%
\setunit*{\addspace}\newblock
	\space Thesis, \printlist{institution}
	\usebibmacro{eprint}%
	\printfield{comments}%
\setunit*{\addspace}
	\printfield{doi}%
	\usebibmacro{finentry}%
}

\DeclareBibliographyDriver{inproceedings}{%
	\usebibmacro{bibindex}%
	\usebibmacro{begentry}%
	\usebibmacro{author}%
\setunit{\addspace}%
	\usebibmacro{date}%
\newunit\newblock
	\usebibmacro{title}%
\newunit\newblock
	\printfield{booktitle}\addcomma%
\setunit*{\addspace}%
	\printfield{series}\addcomma%
\setunit*{\addspace}\newblock
	\usebibmacro{editor}\addcomma
\setunit*{\addspace}\newblock
	\printlist{publisher}
\setunit*{\addspace}%
	\printfield{volume}%
\setunit*{\addspace}%
	\printfield{pages}
\setunit*{\addspace}
	\usebibmacro{eprint}%
	\printfield{comments}%
\setunit*{\addnnspace}
	\printfield{doi}
	\usebibmacro{finentry}%
}

\DeclareBibliographyDriver{inbook}{%
	\usebibmacro{bibindex}%
	\usebibmacro{begentry}%
	\usebibmacro{author}%
\setunit{\addspace}%
	\usebibmacro{date}%
\newunit\newblock
	\usebibmacro{title}%
\newunit\newblock
	\printfield{booktitle}\addcomma%
\setunit*{\addspace}%
	\printfield{series}\addcomma%
\setunit*{\addspace}\newblock
	\usebibmacro{editor}\addcomma
\setunit*{\addspace}\newblock
	\printlist{publisher}
\setunit*{\addspace}%
	\printfield{volume}%
\setunit*{\addspace}%
	\printfield{pages}
\setunit*{\addspace}
	\usebibmacro{eprint}%
	\printfield{comments}%
\setunit*{\addspace}
	\printfield{doi}
	\usebibmacro{finentry}%
}

\DeclareBibliographyDriver{incollection}{%
	\usebibmacro{bibindex}%
	\usebibmacro{begentry}%
	\usebibmacro{author}%
\setunit{\addspace}%
	\usebibmacro{date}%
\newunit\newblock
	\usebibmacro{title}%
\newunit\newblock
	\printfield{booktitle}\addcomma%
\setunit*{\addspace}%
	\printfield{series}\addcomma%
\setunit*{\addspace}\newblock
	\printlist{publisher}
\setunit*{\addspace}%
	\printfield{volume}%
\setunit*{\addspace}%
	\printfield{pages}
\setunit*{\addspace}
	\usebibmacro{eprint}%
	\printfield{comments}%
\setunit*{\addspace}
	\printfield{doi}
	\usebibmacro{finentry}%
}

\AtEveryBibitem{%
	\clearfield{day}%
	\clearfield{month}%
} 
 
\DeclareNameAlias{sortname}{given-family}
\DeclareNameAlias{default}{given-family}

\newcommand{\acknofootnote}{%
	\blfootnote{%
		\hspace*{1em}%
		Jonathan Hermon%
	\hfill%
		Sam Olesker-Taylor%
		\hspace*{1em}%
	\\%
		\href{mailto:jhermon@math.ubc.ca}{jhermon@math.ubc.ca},
		\href{http://www.math.ubc.ca/~jhermon/}{math.ubc.ca/$\sim$jhermon/}%
	\hfill%
		\href{mailto:sam.ot@posteo.co.uk}{sam.ot@posteo.co.uk},
		\href{https://mathematicalsam.wordpress.com}{mathematicalsam.wordpress.com}%
	\\%
		University of British Columbia, Vancouver, Canada%
	\hfill%
		Department of Mathematical Sciences, University of Bath, UK%
	\\%
		Supported by EPSRC EP/L018896/1 and an NSERC Grant%
	\hfill%
		Supported by EPSRC Grants 1885554 and EP/N004566/1%
	\par\smallskip\par
	\centering%
		The vast majority of this work was undertaken while both authors were at the University of Cambridge%
	}
}

\newcommand*{\mm}{\ensuremath{L}}
\newcommand*{\pp}{\ensuremath{p}}
\newcommand*{\qq}{\ensuremath{q}}

\newcommand*{\LL}{\ensuremath{L}}

\newcommand{\printtoc}[1]{%
	\ifthenelse%
		{\equal{#1}{1}}%
		{\sffamily\boldmath\tableofcontents\unboldmath\normalfont}%
		{\newpage\small\sffamily\boldmath\tableofcontents\unboldmath\normalfont\normalsize}%
	}

\usepackage{indentfirst, environ}
\usepackage[utf8]{inputenc}

\usepackage{empheq}

\usepackage{manyfoot}
\SetFootnoteHook{\hspace*{-1.8em}}
\DeclareNewFootnote{bl}[gobble]
\newcommand{\blfootnote}[1]{\footnotebl{\sffamily#1}}

\usepackage{graphicx}
\graphicspath{{../Figures/}}
\usepackage[labelsep=period]{caption}

\usepackage{chngcntr}
\counterwithin{figure}{section}

\usepackage[UKenglish]{babel}
\usepackage{amsmath, amsthm, amssymb, mathrsfs, bm, mathtools}
\usepackage{wasysym}
\usepackage{mathtools}

\allowdisplaybreaks

\usepackage{cases}

\usepackage{titlesec}
\usepackage{titletoc}

\titleformat{\title}
{\sffamily \Huge \bfseries \scshape \boldmath}{}{}{}
\titleformat{\section}
{\sffamily \Large \bfseries \scshape \boldmath}{\thesection}{1em}{}
\titleformat{\subsection}
{\sffamily \large \bfseries \scshape \boldmath}{\thesubsection}{1em}{}
\titleformat{\subsubsection}
{\sffamily \normalsize \bfseries \scshape \boldmath}{\thesubsubsection}{1em}{}
\titleformat{\paragraph}
{\sffamily \normalsize \bfseries \scshape \boldmath}{\theparagraph}{1em}{}
\titleformat{\subparagraph}[runin]
{\sffamily \normalsize \bfseries \scshape \boldmath}{\thesubparagraph}{1em}{}

\def\IfAmpersandUseAlign#1#2&#3\EndIfAmpersandUseAlign
{%
	\if\relax\detokenize{#3}\relax
	\begin{equation*}%
	#1%
	\end{equation*}%
	\else
	\begin{align*}%
	#1%
	\end{align*}%
	\fi
}
\def\[#1\]%
{%
	\IfAmpersandUseAlign{#1}#1&\EndIfAmpersandUseAlign
}

\usepackage{eqparbox}
\newcommand{\eqmathsbox}[3][\mathrel]{%
	#1{\eqmakebox[#2]{$\displaystyle#3$}}%
}

\usepackage{scrextend}

\usepackage{enumitem}

\newcommand{\numberingroman}{%
	\renewcommand{\labelenumi}{(\roman{enumi})}%
	\renewcommand{\theenumi}{(\roman{enumi})}%
}

\setlist[description]{%
	topsep		= 0pt,		
	noitemsep,				
	font		= {\mdseries\itshape},	
}

\usepackage[dvipsnames,hyperref]{xcolor}

\newcommand{\nt}{\addtocounter{equation}{1}\tag{\theequation}}

\newcommand{\bcdot}{\ensuremath{\bm{\cdot}}}
\let\mod\relax
\DeclareMathOperator{\mod}{\, mod}
\DeclareMathOperator{\MOD}{mod}

\DeclareMathOperator{\Ent}{Ent}

\newcommand{\Quad}[1]{
	\mathchoice
	{\quad\text{#1}\quad}
	{\text{ #1 }}
	{\text{ #1 }}
	{\text{ #1 }}
}

\newcommand{\Qforall}{\Quad{for all}}
\newcommand{\Qfor}{\Quad{for}}
\newcommand{\Qand}{\Quad{and}}
\newcommand{\Qwhere}{\Quad{where}}

\newcommand{\supp}{\mathrm{supp}}

\newcommand{\id}{\mathsf{id}}

\newcommand{\cq}{\coloneqq}
\renewcommand{\epsilon}{\varepsilon}

\newcommand{\eps}{\epsilon}

\usepackage{xifthen}

\newcommand{\binomt}[2]{ \textstyle \binom{#1}{#2} \displaystyle }

\newcommand{\maxt}[1]{ \textstyle \max_{#1} \displaystyle }
\newcommand{\supt}[1]{ \textstyle \sup_{#1} \displaystyle }
\newcommand{\mint}[1]{ \textstyle \min_{#1} \displaystyle }

\usepackage{calc}
\newlength{\halfplusheight}
\DeclareMathOperator*{\sumTT}{\textstyle\sum}
\newcommand{\sumT}[2][]{
	\ifthenelse{\isempty{#1}}
	{\sumTT_{#2}}
	{\sumTT_{#2}^{#1}}
}
\newcommand{\sumt}[2][]{
	\ifthenelse{\isempty{#1}}
	{\textstyle \sum_{#2}      \displaystyle}
	{\textstyle \sum_{#2}^{#1} \displaystyle}
}
\newcommand{\sumd}[2][]{
	\ifthenelse{\isempty{#1}}
	{\displaystyle \sum_{#2}}
	{\displaystyle \sum_{#2}^{#1}}
}

\newcommand{\intt}[2][]{
	\ifthenelse{\isempty{#1}}
	{\textstyle \int_{#2}      \displaystyle}
	{\textstyle \int_{#2}^{#1} \displaystyle}
}

\newcommand{\prodt}[2][]{
	\ifthenelse{\isempty{#1}}
	{\textstyle \prod_{#2}      \displaystyle}
	{\textstyle \prod_{#2}^{#1} \displaystyle}
}
\newcommand{\prodd}[2][]{
	\ifthenelse{\isempty{#1}}
	{\prod_{#2}}
	{\prod_{#2}^{#1}}
}

\usepackage{xspace}

\let\originalexp\exp
\let\exp\relax

\DeclareRobustCommand{\exp} [1]{\originalexp(#1)}
\newcommand{\expb} [1]{\originalexp\bigl( #1 \bigr)}

\newcommand{\expbb}[1]{\originalexp\biggl( #1 \biggr)}

\newcommand{\gap}{}
\newcommand{\GAP}[1]{%
	\renewcommand{\gap}{\hspace{#1em}}}

\newcommand{\abs}  [1]{| #1 |}
\newcommand{\absb} [1]{\big| #1 \bigr|}

\newcommand{\norm}  [1]{\lVert #1 \rVert}
\newcommand{\normb} [1]{\big\lVert #1 \bigr\rVert}

\newcommand{\ipr} [1]{ \langle #1 \rangle }
\newcommand{\ipb} [1]{\bigl\langle #1 \bigr\rangle}

\newcommand{\rbr} [1]{ ( #1 ) }
\newcommand{\rbb} [1]{\bigl( #1 \bigr)}
\newcommand{\rbB} [1]{\Bigl( #1 \Bigr)}

\newcommand{\bra} [1]{ \{ #1 \} }
\newcommand{\brb} [1]{\bigl\{ #1 \bigr\}}

\newcommand{\TV}{\mathrm{TV}}

\newcommand{\DIAM}{\mathrm{diam}}

\newcommand{\st}{{ \ \mathrm{st} \ }}

\DeclareMathOperator{\vol}{\mathrm{vol}}

\DeclareMathOperator{\step}{step}

\newcommand{\Unif}{\mathrm{Unif}}
\newcommand{\iid}{\mathrm{iid}}

\newcommand{\Poisson}{\mathrm{Poisson}}

\newcommand{\Bernoulli}{\mathrm{Bernoulli}}

\newcommand{\LS}{\mathrm{LS}}
\newcommand{\MLS}{\mathrm{MLS}}

\newcommand{\Xs}{(X_s)_{s\ge0}}

\newcommand{\floor}[1]{\lfloor #1 \rfloor}
\newcommand{\floorb}[1]{\bigl \lfloor #1 \bigr \rfloor}

\newcommand{\ceil}[1]{\lceil #1 \rceil}
\newcommand{\ceilb}[1]{\bigl\lceil #1 \bigr\rceil}

\newcommand{\midb}{\bigm\vert}

\newcommand{\one}  [1]{\bm1( #1 )}

\newcommand{\tv} [1]{\lVert #1 \rVert_{\mathrm{TV}}}

\newcommand{\logk}[1][]{
	\ifthenelse{\equal{}{#1}}
	{\log k}
	{(\log k)^{#1}}
}

\newcommand{\logn}[1][]{
	\ifthenelse{\equal{}{#1}}
	{\log n}
	{(\log n)^{#1}}
}

\newcommand{\logm}[1][]{
	\ifthenelse{\equal{}{#1}}
	{\log m}
	{(\log m)^{#1}}
}

\newcommand{\loglogn}[1][]{
	\ifthenelse{\equal{}{#1}}
	{\log\log n}
	{(\log\log n)^{#1}}
}

\newcommand{\relent}[2] {D(#1 \, \Vert \, #2)}

\newcommand{\prt}[2][]{
	\ifthenelse{\equal{}{#1}}
	{\mathbb{P}(#2)}
	{\mathbb{P}_{#1}(#2)}
}
\newcommand{\pr}[2][]{
	\mathchoice
	{\ifthenelse{\isempty{#1}}
		{\mathbb{P}\bigl(#2\bigr)}
		{\mathbb{P}_{#1}\bigl(#2\bigr)}}
	{\ifthenelse{\isempty{#1}}
		{\mathbb{P}(#2)}
		{\mathbb{P}_{#1}(#2)}}
	{\ifthenelse{\isempty{#1}}
		{\mathbb{P}(#2)}
		{\mathbb{P}_{#1}(#2)}}
	{\ifthenelse{\isempty{#1}}
		{\mathbb{P}(#2)}
		{\mathbb{P}_{#1}(#2)}}
}
\newcommand{\prb}[2][]{
	\ifthenelse{\equal{}{#1}}
	{\mathbb{P}\bigl( #2 \bigr)}
	{\mathbb{P}_{#1}\bigl( #2 \bigr)}
}
\newcommand{\prB}[2][]{
	\ifthenelse{\equal{}{#1}}
	{\mathbb{P}\Bigl( #2 \Bigr)}
	{\mathbb{P}_{#1}\Bigl( #2 \Bigr)}
}
\newcommand{\prbb}[2][]{
	\ifthenelse{\equal{}{#1}}
	{\mathbb{P}\biggl( #2 \biggr)}
	{\mathbb{P}_{#1}\biggl( #2 \biggr)}
}
\newcommand{\prBB}[2][]{
	\ifthenelse{\equal{}{#1}}
	{\mathbb{P}\Biggl( #2 \Biggr)}
	{\mathbb{P}_{#1}\Biggl( #2 \Biggr)}
}
\newcommand{\prs}[2][]{
	\ifthenelse{\equal{}{#1}}
	{\mathbb{P}\left( #2 \right)}
	{\mathbb{P}_{#1}\left( #2 \right)}
}

\newcommand{\qr}[2][]{
	\mathchoice
	{\ifthenelse{\isempty{#1}}
		{\mathbb{Q}\bigl(#2\bigr)}
		{\mathbb{Q}_{#1}\bigl(#2\bigr)}}
	{\ifthenelse{\isempty{#1}}
		{\mathbb{Q}(#2)}
		{\mathbb{Q}_{#1}(#2)}}
	{\ifthenelse{\isempty{#1}}
		{\mathbb{Q}(#2)}
		{\mathbb{Q}_{#1}(#2)}}
	{\ifthenelse{\isempty{#1}}
		{\mathbb{Q}(#2)}
		{\mathbb{Q}_{#1}(#2)}}
}
\newcommand{\qrb}[2][]{
	\ifthenelse{\equal{}{#1}}
	{\mathbb{Q}\bigl( #2 \bigr)}
	{\mathbb{Q}_{#1}\bigl( #2 \bigr)}
}
\newcommand{\qrB}[2][]{
	\ifthenelse{\equal{}{#1}}
	{\mathbb{Q}\Bigl( #2 \Bigr)}
	{\mathbb{Q}_{#1}\Bigl( #2 \Bigr)}
}
\newcommand{\qrbb}[2][]{
	\ifthenelse{\equal{}{#1}}
	{\mathbb{Q}\biggl( #2 \biggr)}
	{\mathbb{Q}_{#1}\biggl( #2 \biggr)}
}
\newcommand{\qrBB}[2][]{
	\ifthenelse{\equal{}{#1}}
	{\mathbb{Q}\Biggl( #2 \Biggr)}
	{\mathbb{Q}_{#1}\Biggl( #2 \Biggr)}
}
\newcommand{\qrs}[2][]{
	\ifthenelse{\equal{}{#1}}
	{\mathbb{Q}\left( #2 \right)}
	{\mathbb{Q}_{#1}\left( #2 \right)}
}

\newcommand{\ext}[2][]{
\ifthenelse{\equal{}{#1}}
{\mathbb{E}(#2)}
{\mathbb{E}_{#1}(#2)}
}
\newcommand{\ex}[2][]{
	\mathchoice
	{\ifthenelse{\isempty{#1}}
		{\mathbb{E}\bigl(#2\bigr)}
		{\mathbb{E}_{#1}\bigl(#2\bigr)}}
	{\ifthenelse{\isempty{#1}}
		{\mathbb{E}(#2)}
		{\mathbb{E}_{#1}(#2)}}
	{\ifthenelse{\isempty{#1}}
		{\mathbb{E}(#2)}
		{\mathbb{E}_{#1}(#2)}}
	{\ifthenelse{\isempty{#1}}
		{\mathbb{E}(#2)}
		{\mathbb{E}_{#1}(#2)}}
}
\newcommand{\exb}[2][]{
	\ifthenelse{\equal{}{#1}}
	{\mathbb{E}\bigl( #2 \bigr)}
	{\mathbb{E}_{#1}\bigr( #2 \bigr)}
}
\newcommand{\exB}[2][]{
	\ifthenelse{\equal{}{#1}}
	{\mathbb{E}\Bigl( #2 \Bigr)}
	{\mathbb{E}_{#1}\Bigl( #2 \Bigr)}
}
\newcommand{\exbb}[2][]{
	\ifthenelse{\equal{}{#1}}
	{\mathbb{E}\biggl( #2 \biggr)}
	{\mathbb{E}_{#1}\biggl( #2 \biggr)}
}
\newcommand{\exBB}[2][]{
	\ifthenelse{\equal{}{#1}}
	{\mathbb{E}\Biggl( #2 \Biggr)}
	{\mathbb{E}_{#1}\Biggl( #2 \Biggr)}
}

\newcommand{\fx}[2][]{
	\ifthenelse{\equal{}{#1}}
	{\mathbb{F}(#2)}
	{\mathbb{F}_{#1}(#2)}
}
\newcommand{\fxb}[2][]{
	\ifthenelse{\equal{}{#1}}
	{\mathbb{F}\bigl( #2 \bigr)}
	{\mathbb{F}_{#1}\bigr( #2 \bigr)}
}
\newcommand{\fxB}[2][]{
	\ifthenelse{\equal{}{#1}}
	{\mathbb{F}\Bigl( #2 \Bigr)}
	{\mathbb{F}_{#1}\Bigl( #2 \Bigr)}
}
\newcommand{\fxbb}[2][]{
	\ifthenelse{\equal{}{#1}}
	{\mathbb{F}\biggl( #2 \biggr)}
	{\mathbb{F}_{#1}\biggl( #2 \biggr)}
}
\newcommand{\fxBB}[2][]{
	\ifthenelse{\equal{}{#1}}
	{\mathbb{F}\Biggl( #2 \Biggr)}
	{\mathbb{F}_{#1}\Biggl( #2 \Biggr)}
}

\newcommand{\Var}[1]{\mathbb{V}\mathrm{ar}(#1)}
\newcommand{\Varb}[2][]{
	\ifthenelse{\equal{}{#1}}
	{\mathbb{V}\mathrm{ar} \bigl(#2\bigr)}
	{\mathbb{V}\mathrm{ar}_{#1} \bigl(#2\bigr)}
}
\newcommand{\VAR}[2][]{
	\ifthenelse{\equal{}{#1}}
	{\mathrm{Var}(#2)}
	{\mathrm{Var}_{#1}(#2)}
}

\newcommand{\Oh}  [1]{\mathcal{O}( #1 )}
\newcommand{\Ohb} [1]{\mathcal{O}\bigl( #1 \bigr)}

\newcommand{\Ohbb}[1]{\mathcal{O}\biggl( #1 \biggr)}

\newcommand{\oh}  [1]{o( #1 )}
\newcommand{\ohb} [1]{o\bigl( #1 \bigr)}

\newcommand{\mba}{\mathbb{A}}

\newcommand{\mbn}{\mathbb{N}}

\newcommand{\mbr}{\mathbb{R}}

\newcommand{\mbz}{\mathbb{Z}}

\newcommand{\mce}{\mathcal{E}}

\newcommand{\mci}{\mathcal{I}}

\newcommand{\mcl}{\mathcal{L}}

\newcommand{\mco}{\mathcal{O}}

\newcommand{\mcu}{\mathcal{U}}

\newcommand{\mfgcd}{\mathfrak{g}}

\newcommand{\mfm}{\mathfrak{M}}

\newcommand{\toinf}[1]{\ensuremath{#1\to\infty}}

\newcommand{\tozero}[1]{\ensuremath{#1\to0}}

\newcommand{\kinf}{{k\to\infty}}

\newcommand{\ninf}{{n\to\infty}}

\newcommand{\askinf}{\textnormal{as $k\to\infty$}}

\newcommand{\tinf}{{t\to\infty}}

\makeatletter
\newenvironment{subtheorem}[1]{%
	\def\subtheoremcounter{#1}%
	\refstepcounter{#1}%
	\protected@edef\theparentnumber{\csname the#1\endcsname}%
	\setcounter{parentnumber}{\value{#1}}%
	\setcounter{#1}{0}%
	\expandafter\def\csname the#1\endcsname{\theparentnumber\alph{#1}}%
	\expandafter\def\csname theH#1\endcsname{thm.\theparentnumber\alph{#1}}%
	\unskip\ignorespaces
}{%
	\setcounter{\subtheoremcounter}{\value{parentnumber}}%
	\ignorespacesafterend
}
\makeatother
\newcounter{parentnumber}

\makeatletter
\newenvironment{subtheorem-num}[1]{%
	\def\subtheoremcounter{#1}%
	\refstepcounter{#1}%
	\protected@edef\theparentnumber{\csname the#1\endcsname}%
	\setcounter{parentnumber}{\value{#1}}%
	\setcounter{#1}{0}%
	\expandafter\def\csname the#1\endcsname{\theparentnumber.\arabic{#1}}%
	\expandafter\def\csname theH#1\endcsname{thm.\theparentnumber.\arabic{#1}}%
	\unskip\ignorespaces
}{%
	\setcounter{\subtheoremcounter}{\value{parentnumber}}%
	\ignorespacesafterend
}
\makeatother

\usepackage[colorlinks, pdfauthor = {Jonathan\ Hermon\ and\ Sam\ Olesker-Taylor}, pdfusetitle]{hyperref}
\hypersetup{
	colorlinks,
	linkcolor	= {blue},
	filecolor	= {red},
	urlcolor	= {blue},
	citecolor	= {ForestGreen}
}

\usepackage[capitalise]{cleveref}

\crefformat{equation}{\textup{(#2#1#3)}}
\crefmultiformat{equation}{\textup{(#2#1#3}}{\textup{, #2#1#3)}}{\textup{, #2#1#3}}{\textup{, #2#1#3)}}
\crefrangeformat{equation}{\textup{(#3#1#4--#5#2#6)}}

\newenvironment{Proof}[1][\proofname]{%
	\proof[\upshape\bfseries\sffamily\boldmath#1]
}{\endproof}

\usepackage{etoolbox}
\usepackage{needspace}

\newtheoremstyle{sfsl}
{1\baselineskip}		
{1\baselineskip}		
{\slshape}				
{}						
{\bfseries\sffamily}	
{.}						
{0.5em}					
{\thmname{#1}\thmnumber{ #2}\thmnote{ {\mdseries(#3)}}}

\newtheoremstyle{sfup}
{1\baselineskip}		
{1\baselineskip}		
{\upshape}				
{}						
{\bfseries\sffamily}	
{.}						
{0.5em}					
{\thmname{#1}\thmnumber{ #2}\thmnote{ {\mdseries(#3)}}}

\theoremstyle{sfsl}

\newtheorem*{thm*}{Theorem}
\newtheorem{thm} {Theorem}[section]
\crefname{thm}{Theorem}{Theorems}

\newtheorem*{introthm*}{Theorem}

\crefname{introthm}{Theorem}{Theorems}

\newtheorem*{cor*}{Corollary}
\newtheorem{cor} [thm]{Corollary}
\crefname{cor}{Corollary}{Corollaries}

\newtheorem*{introcor*}{Corollary}

\crefname{introcor}{Corollary}{Corollaries}

\newtheorem*{introconj*}{Conjecture}

\crefname{introconj}{Conjecture}{Conjectures}

\newtheorem*{introques*}{Question}

\crefname{introques}{Question}{Questions}

\newtheorem*{lem*}    {Lemma}
\newtheorem{lem} [thm]{Lemma}
\crefname{lem}{Lemma}{Lemmas}

\newtheorem*{introlem*}{Lemma}

\crefname{introlem}{Lemma}{Lemmas}

\newtheorem*{prop*}    {Proposition}
\newtheorem{prop} [thm]{Proposition}
\crefname{prop}{Proposition}{Propositions}

\newtheorem*{clm*}    {Claim}

\crefname{clm}{Claim}{Claims}

\newtheorem*{defn*}    {Definition}
\newtheorem{defn} [thm]{Definition}
\crefname{defn}{Definition}{Definitions}

\newtheorem*{introdefn*}{Definition}

\crefname{introdefn}{Definition}{Definitions}

\newtheorem{innercustomgenericsl}{\customgenericnamesl}
\providecommand{\customgenericnamesl}{}
\newcommand{\newcustomtheoremsl}[2]{%
	\newenvironment{#1}[1]
	{%
		\renewcommand\customgenericnamesl{#2}%
		\renewcommand\theinnercustomgenericsl{##1}%
		\innercustomgenericsl
	}
	{\endinnercustomgenericsl}
}

\newcustomtheoremsl{customthm}{Theorem}
\newcustomtheoremsl{customprop}{Proposition}
\newcustomtheoremsl{customlemma}{Lemma}
\newcustomtheoremsl{customrmk}{Remark}
\newcustomtheoremsl{customconj}{Conjecture}
\newcustomtheoremsl{customdefn}{Definition}
\newcustomtheoremsl{customquestion}{Question}
\newcustomtheoremsl{customopenproblem}{Open Problem}
\newcustomtheoremsl{customhyp}{Hypothesis}

\newtheorem*{conj*}   {Conjecture}

\crefname{conj}{Conjecture}{Conjectures}

\newenvironment{conj-ind*}
	{\begin{quote}\textsf{\textbf{Conjecture.}}\slshape}
	{\end{quote}}
\newenvironment{conj-ind}
	{\begin{quote}\vspace{-\glueexpr\baselineskip+\topsep}\begin{customconj}}
	{\end{customconj}\end{quote}}

\newenvironment{question-ind*}
	{\begin{quote}\textsf{\textbf{Question.}}\slshape}
	{\end{quote}}
\newenvironment{question-ind}
	{\begin{quote}\vspace{-\glueexpr\baselineskip+\topsep}\begin{customquestion}}
	{\end{customquestion}\end{quote}}

\newenvironment{openproblem-ind*}
	{\begin{quote}\textsf{\textbf{Open Problem.}}\slshape}
	{\end{quote}}
\newenvironment{openproblem-ind}
	{\begin{quote}\vspace{-\glueexpr\baselineskip+\topsep}\begin{customopenproblem}}
	{\end{customopenproblem}\end{quote}}

\newtheorem*{hypothesis*}{Hypothesis}

\newtheorem*{hyp*}{Hypothesis}

\crefname{hyp}{Hypothesis}{Hypotheses}

\newtheorem*{rmk*}{Remark}

\theoremstyle{sfup}

\newtheorem{innercustomgenericup}{\customgenericnameup}
\providecommand{\customgenericnameup}{}
\newcommand{\newcustomtheoremup}[2]{%
	\newenvironment{#1}[1]
	{%
		\renewcommand\customgenericnameup{#2}%
		\renewcommand\theinnercustomgenericup{##1}%
		\innercustomgenericup
	}
	{\endinnercustomgenericup}
}

\crefname{exm} {Example}{Examples}
\crefname{exmT}{Example}{Examples}

\newtheorem{rmkT}[thm]{Remark}
	\newenvironment{rmkt}
	{\pushQED{\qed}\rmkT}
	{\popQED\endrmkT}
\crefname{rmk} {Remark}{Remarks}
\crefname{rmkT}{Remark}{Remarks}

\newtheorem*{rmkTT}{Remark}
\newenvironment{rmkt*}
	{\pushQED{\qed}\rmkTT}
	{\popQED\endrmkTT}

\newcustomtheoremup{customrmkT}{Remark}

\crefname{rmks} {Remarks}{Remarks}
\crefname{rmksT}{Remarks}{Remarks}

\newtheorem*{rmks*} {Remarks}
\newtheorem*{rmksTT}{Remarks}
\newenvironment{rmkst*}
	{\pushQED{\qed}\rmksTT}
	{\popQED\endrmksTT}

\crefname{intrormk} {Remark}{Remarks}
\crefname{intrormkT}{Remark}{Remarks}

\newtheorem*{intrormk*} {Remark}
\newtheorem*{intrormkTT}{Remark}
\newenvironment{intrormkt*}
	{\pushQED{\qed}\intrormkTT}
	{\popQED\endintrormkTT}

\newtheorem*{exm*} {Example}
\newtheorem*{exmTT}{Lemma}
	\newenvironment{exmt*}
	{\pushQED{\qed}\exmTT}
	{\popQED\endexmTT}

\newtheorem*{note*} {Note}
\newtheorem*{noteTT}{Note}
	\newenvironment{notet*}
	{\pushQED{\qed}\noteTT}
	{\popQED\endnoteTT}

\catcode`\^^A=14

\usepackage{xparse}

\newcounter{mixedsubequations}

\ExplSyntaxOn

\int_new:N \g_mixedsubeq_int

\NewDocumentEnvironment{mixedsubequations}{o}
{
	\IfNoValueTF { #1 }
	{
		\addtocounter{equation}{-\g_mixedsubeq_int}
		\stepcounter{mixedsubequations}
	}
	{
		\int_gset:Nn \g_mixedsubeq_int { \clist_count:n { #1 } }
		\clist_map_inline:nn { #1 }
		{
			\refstepcounter{equation}\label{##1}
		}
		\addtocounter{equation}{-\g_mixedsubeq_int}
		\setcounter{mixedsubequations}{1}
	}
	\domixedsubequations
}
{\ignorespacesafterend}

\NewDocumentCommand{\domixedsubequations}{}
{
	\cs_set:Npx \theequation
	{
		\exp_not:o { \theequation }
		\exp_not:n { \alph{mixedsubequations} }
	}
	\ignorespaces
}
\ExplSyntaxOff

\makeatletter
\let\save@mathaccent\mathaccent
\newcommand*\if@single[3]{%
  \setbox0\hbox{${\mathaccent"0362{#1}}^H$}%
  \setbox2\hbox{${\mathaccent"0362{\kern0pt#1}}^H$}%
  \ifdim\ht0=\ht2 #3\else #2\fi
  }
\newcommand*\rel@kern[1]{\kern#1\dimexpr\macc@kerna}
\newcommand*\widebar[1]{\@ifnextchar^{{\wide@bar{#1}{0}}}{\wide@bar{#1}{1}}}
\newcommand*\wide@bar[2]{\if@single{#1}{\wide@bar@{#1}{#2}{1}}{\wide@bar@{#1}{#2}{2}}}
\newcommand*\wide@bar@[3]{%
  \begingroup
  \def\mathaccent##1##2{%
    \let\mathaccent\save@mathaccent
    \if#32 \let\macc@nucleus\first@char \fi
    \setbox\z@\hbox{$\macc@style{\macc@nucleus}_{}$}%
    \setbox\tw@\hbox{$\macc@style{\macc@nucleus}{}_{}$}%
    \dimen@\wd\tw@
    \advance\dimen@-\wd\z@
    \divide\dimen@ 3
    \@tempdima\wd\tw@
    \advance\@tempdima-\scriptspace
    \divide\@tempdima 10
    \advance\dimen@-\@tempdima
    \ifdim\dimen@>\z@ \dimen@0pt\fi
    \rel@kern{0.6}\kern-\dimen@
    \if#31
      \overline{\rel@kern{-0.6}\kern\dimen@\macc@nucleus\rel@kern{0.4}\kern\dimen@}%
      \advance\dimen@0.4\dimexpr\macc@kerna
      \let\final@kern#2%
      \ifdim\dimen@<\z@ \let\final@kern1\fi
      \if\final@kern1 \kern-\dimen@\fi
    \else
      \overline{\rel@kern{-0.6}\kern\dimen@#1}%
    \fi
  }%
  \macc@depth\@ne
  \let\math@bgroup\@empty \let\math@egroup\macc@set@skewchar
  \mathsurround\z@ \frozen@everymath{\mathgroup\macc@group\relax}%
  \macc@set@skewchar\relax
  \let\mathaccentV\macc@nested@a
  \if#31
    \macc@nested@a\relax111{#1}%
  \else
    \def\gobble@till@marker##1\endmarker{}%
    \futurelet\first@char\gobble@till@marker#1\endmarker
    \ifcat\noexpand\first@char A\else
      \def\first@char{}%
    \fi
    \macc@nested@a\relax111{\first@char}%
  \fi
  \endgroup
}
\makeatother

\renewcommand{\thesection}{\Alph{section}}

\numberwithin{equation}{section}

\title{\sffamily Supplementary Material for Random Cayley Graphs Project}

\author{\sffamily Jonathan Hermon\quad Sam Olesker-Taylor}%

\date{}

\usepackage{geometry}
\begin{document}

\maketitle

\acknofootnote

\vspace{-6ex}

\renewcommand{\abstractname}{\sffamily Abstract}
\begin{abstract}
This document contains supplementary material for the main articles, namely
\cite{HOt:rcg:abe:cutoff,HOt:rcg:matrix,HOt:rcg:abe:extra,HOt:rcg:abe:geom}.

\smallskip

We prove refined results about simple random walks on the integers and on the cycle.
See \S\ref{sec-p0:se} and \S\ref{sec-p0:re}, respectively.
We are primarily interested in the entropy of these random walks at certain times and how this entropy changes when the time changes slightly.
Additionally, we prove some large deviation and exit time estimates.
See \S\ref{sec-p0:rp} and \S\ref{sec-p0:gap}, respectively.

We prove some results on the size of discrete lattice balls and how this size changes when the radius changes slightly.
We do this in a general $L_\qq$ norm, with $\qq \in [1,\infty]$.
See \S\ref{sec-p0:balls}.

We also prove some other technical results deferred from \cite{HOt:rcg:abe:cutoff,HOt:rcg:matrix,HOt:rcg:abe:extra,HOt:rcg:abe:geom}.
See \S\ref{sec-p0:deferred}.

\smallskip

We hope that some of the results, particularly the simple random walk estimates, will be useful in their own right for other researchers.
\end{abstract}

\setcounter{tocdepth}{2}
\sffamily
\vfill
\tableofcontents
\vspace*{2ex}
\normalfont

\numberingroman

\newpage
\section*{Notation and Terminology}

While notation will often be recalled later, we list here the majority of what we use below.

\begin{itemize}
	\item 
	The \textit{simple random walk}, abbreviated \textit{SRW}, on $\mbz^k$ (or $\mbz_\mm^k$) is the rate-1 RW which, in each step, chooses a coordinate uniformly at random and adds/subtracts 1 from this value (mod $\mm$) with equal probability.
	The \textit{directed random walk}, abbreviated \textit{DRW}, on $\mbz^k$ (or $\mbz_\mm^k$) is the same except that it only ever adds 1 (mod $\mm$).
	
	When we wish to specify the DRW, we add a $+$-superscript; for the SRW, we add a $-$-superscript.
	When we do not wish to specify to which walk we are referring (as the statement applies for both), we simply speak of the \textit{random walk}, abbreviated \textit{RW}; if we wish to emphasise that the statement applies for both, then we sometimes add a $\pm$-superscript.
	
	\item 
	Rate-1 RW on $\mbz^k$:
	\begin{itemize}
		\item 
		$W = \rbb{ W_i(t) \mid i \in [k], \: t \ge 0 }$ is a rate-1 RW on $\mbz^k$;
		
		\item 
		$\mu_t(\cdot) \cq \mcl\rbb{ W(t) }$, the \textit{law};
		
		\item 
		$Q(t) \cq - \log \mu_t\rbb{ W(t) }$, the \textit{random (Shannon) entropy};
		
		\item 
		$h(t) \cq \exb{Q(t)}$, the \textit{(Shannon) entropy}.
	\end{itemize}
	
	\item 
	Rate-1 RW on $\mbz$:
	\begin{itemize}
		\item 
		$W_i = \rbr{ W_1(t) \mid i \in [k], \: t \ge 0 }$ is a rate-$1/k$ RW on $\mbz$;
		
		\item 
		$X = \rbr{ X(s) = X_s \mid s \ge 0 }$ defined by $X_s \cq W_1(sk)$ is a rate-1 RW on $\mbz$;
		
		\item 
		$\nu_s(\cdot) \cq \mcl\rbb{ W_1(sk) } = \mcl\rbb{ X(s) }$, the \textit{law};
		
		\item 
		$Q_i(t) \cq - \log \nu_{t/k}\rbb{ W_i(t) }$, the \textit{random (Shannon) entropy};
		
		\item 
		$H(s) \cq \exb{ Q_1(sk) } = \exb{ - \log \nu_s\rbb{ X(s) } }$, the \textit{(Shannon) entropy}.
	\end{itemize}
	
	\item 
	Rate-1 RW on $\mbz_\mm^k$:
	\begin{itemize}
		\item 
		$W_\mm = \rbb{ W_{\mm, i}(t) \mid t \ge 0 }$ defined by $W_{\mm, i}(t) \cq W_i(t) \mod \mm$ is a rate-1 RW on $\mbz_\mm$;
		
		\item 
		$\mu_{\mm, t}(\cdot) \cq \mcl\rbb{ W_\gamma(t) }$, the \textit{law};
		
		\item 
		$Q_\mm(t) \cq - \log \mu_{\mm, t}\rbb{ W_\mm(t) }$, the \textit{random (Shannon) entropy};
		
		\item 
		$h_\mm(t) \cq \exb{ Q_\mm(t) }$, the \textit{(Shannon) entropy};
		
		\item 
		$r_\mm(t) \cq \log(\mm^k) - h_\mm(t)$, the \textit{relative entropy} wrt $\Unif(\mbz_\mm^k)$.
	\end{itemize}
	
	\item 
	Rate-1 RW on $\mbz_\mm$:
	\begin{itemize}
		\item 
		$W_{\mm, i} = \rbb{ W_{\mm, i}(t) \mid t \ge 0 }$ is a rate-$1/k$ RW on $\mbz_\mm$;
		
		\item 
		$X_\mm = \rbb{ X_\mm(s) \mid s \ge 0 }$ defined by $X_\mm(s) \cq W_{\mm, 1}(sk)$ is a rate-1 RW on $\mbz_\mm$;
		
		\item 
		$\nu_{\mm, s}(\cdot) \cq \mcl\rbb{ W_{\mm, 1}(sk) } = \mcl\rbb{ X_\mm(s) }$, the \textit{law};
		
		\item 
		$Q_{\mm, i}(t) \cq - \log \nu_{t/k}\rbb{ W_{\mm, i}(t) }$, the \textit{random (Shannon) entropy};
		
		\item 
		$H_\mm(s) \cq \exb{ Q_{\mm, 1}(sk) } = \exb{ - \log \nu_{\mm, s}\rbb{ X_\mm(s) } }$, the \textit{(Shannon) entropy};
		
		\item 
		$R_\mm(s) \cq \log \mm - H_\mm(s)$, the \textit{relative entropy} wrt $\Unif(\mbz_\mm)$.
	\end{itemize}
	
%
%
%
	\item 
	Consider ``mod $\infty$'' to mean no modulation; eg, $W = W_\infty$ or $H = H_\infty$.
	
	\item 
	For $\mm \in \mbn$ and $p \in [1,\infty]$,
	write
	\[
		d_{p,\mm}(s)
	\cq
		\normb{ \pr{ X_\mm(s) \in \cdot } - \Unif(\mbz_\mm) }_{p, \mm}
	=
		\rbb{
			\sumt{x \in \mbz_\mm}
			\tfrac1\mm \absb{ n \pr{ X_\mm(s) = x } - 1 }^p
		}^{1/p};
	\]
	also write $d_{\TV, \mm}(s) \cq \tfrac12 d_{\mm,1}(s)$ for the total variation distance, which can be represented as
	\[
		d_{\TV, \mm}(s)
	=
		\maxt{A \subseteq \mbz_\mm}
		\absb{ \pr{ X_\mm(s) \in A } - \tfrac1\mm \abs A }.
	\]
	For $\mm = \infty$, we usually drop the $\mm$-subscript.
\end{itemize}

\section{Shannon Entropy Estimates and Central Limit Theorem}
\label{sec-p0:se}

%

This part of the appendix (\S\ref{sec-p0:se}) is devoted to properties of the entropic time $t_0$ and cutoff window $t_\alpha - t_0$; this is done through analysis of a CLT for $Q$ (\cref{res-p0:se:CLT}) and variance of $Q_1$ at the entropic time, $\Var{Q_1(t_0)}$.
Accordingly, here we mainly derive properties of the SRW on $\mbz$ evaluated at $t/k$ or of $\Poisson(t/k)$, for $t$ around the entropic time.

\subsection{Key Definitions and Results for Shannon Entropy}

We now define precisely the notion of \textit{entropic times}.
Let $W = (W(t))_{t\ge0}$ be a RW on $\mbz^k$.
Write $\mu_t$, respectively $\nu_s$, for the law of $W(t)$, respectively $W_1(sk)$;
so $\mu_t = \nu_{t/k}^{\otimes k}$.
Define
\[
	Q_i(t) \cq - \log \nu_t\rbb{W_i(t)},
\Quad{and set}
	Q(t) \cq - \log \mu_t\rbb{W(t)} = \sumt[k]{1} Q_i(t).
\]
So $\ex{Q(t)}$ and $\ex{Q_1(t)}$ are the entropies of $W(t)$ and $W_1(t)$, respectively.

\begin{defn}[Entropic and Times]
\label{def-p0:se:t0a}
	For all $k,n \in \mbn$ and all $\alpha \in \mbr$,
	define
	$t_\alpha \cq t_\alpha(k,n)$ so that
	\[
		\ex{ Q_1\rbr{t_\alpha} } = \rbb{ \log n + \alpha \sqrt{v k} }/k
	\Qand
		s_\alpha \cq t_\alpha/k,
	\Qwhere
		v \cq \Varb{Q_1\rbr{t_0}}.
	\]
	We call $t_0$ the \textit{entropic time} and the $\bra{t_\alpha}_{\alpha\in\mbr}$ \textit{cutoff times}.
\end{defn}

The following proposition gives a detailed approximate evaluation of these entropic times.

\begin{prop}[Entropic and Cutoff Times]
\label{res-p0:se:t0a}
Assume that $1 \ll \log k \ll \log n$.
Write $\kappa \cq k/\log n$.
For all $\alpha \in \mbr$ and $\lambda > 0$,
the following relations hold, for some functions $f$ and $g$:
we have $t_\alpha \eqsim t_0$;%
\begin{subequations}
	\label{eq-p0:se:t0a}
\begin{alignat*}{3}
	\text{for } k &\ll \log n,
&\Quad{we have}
	t_0 &\eqsim k \cdot n^{2/k} / (2 \pi e)
&\Qand
	t_\alpha - t_0 &\eqsim \sqrt2 \cdot \alpha t_0 / \sqrt k;
\label{eq-p0:se:t0a:k<<logn}
\nt
\\
	\text{for } k &\eqsim \lambda \log n,
&\Quad{we have}
	t_0 &\eqsim k \cdot f(\lambda)
&\Qand
	t_\alpha - t_0 &\eqsim g(\lambda) \cdot \alpha t_0 / \sqrt k;
\label{eq-p0:se:t0a:k==logn}
\nt
\\
	\text{for } k &\gg \log n,
&\Quad{we have}
	t_0 &\eqsim k \cdot 1/(\kappa \log \kappa)
&\Qand
	t_\alpha - t_0 &\eqsim \sqrt{\kappa \log \kappa} \cdot \alpha t_0 / \sqrt k.
\label{eq-p0:se:t0a:k>>logn}
\nt
\end{alignat*}
\end{subequations}
Moreover, $f, g : (0,\infty) \to (0,\infty)$ are continuous functions, whose value differs between the undirected and directed cases.
In particular, for all $\alpha \in \mbr$, in all cases, we have $t_\alpha \eqsim t_0$.
\end{prop}

Observe that $Q(t) = \sumt[k]{1} Q_i(t)$ is a sum of iid random variables.

\begin{prop}[CLT]
\label{res-p0:se:CLT}
	Assume that $1 \ll k \ll \log n$.
	For all $\alpha \in \mbr$,
	we have
	\[
		\pr{ Q(t_\alpha) \le \log n \pm \omega } \to \Psi(\alpha)
	\Qfor
		\omega \cq \Varb{ Q(t_0) }^{1/4} = (vk)^{1/4}.
	\]
	(There is no specific reason for choosing this $\omega$. We just need some $\omega$ with $1 \ll \omega \ll (vk)^{1/2}$.)
\end{prop}

\subsection{Local CLT for RW on $\mbz$}

We repeatedly use a local CLT for Poisson and simple random walk distributions.
We state it here precisely;
the particular version is given in \cite[Theorem 2.5.6]{LL:RW}.

\begin{thm}[Local CLT, {\cite[Theorem 2.5.6]{LL:RW}}]
\label{res-p0:LCLT}
	Let $\varsigma > 0$ and let $s \in (\varsigma, \infty)$; the implicit constants in the $\mco$-notation notation depend on $\varsigma$.
	Let $X = \rbr{ X_s }_{s\ge0}$ be either a rate-1 SRW or rate-1 DRW on $\mbz$.
	For all $x \in \mbr$ with $x - \ex{X_s} \in \mbz$ and $\abs x \le \tfrac12 s$,
	we have
	\[
		\pr{ X_s - \ext{X_s} = x }
	=
		\frac1{\sqrt{2 \pi s}} \expbb{- \frac{x^2}{2s} } \expbb{ \Ohbb{ \frac1{\sqrt s} + \frac{\abs x^3}{s^2} } }.
	\]
	In particular, if $\abs x \le s^{7/12}$, then
	\[
		\pr{ X_s - \ext{X_s} = x }
	=
		\frac1{\sqrt{2 \pi s}} \expbb{- \frac{x^2}{2s} } \expb{ \Ohb{ s^{-1/4} } }.
	\label{eq-p0:LCLT}
	\nt
	\]
\end{thm}

\begin{Proof}
The result for the SRW is given in \cite[Theorem 2.5.6]{LL:RW}.
For the DRW, observe that $X_s \sim \Poisson(s)$ and use Stirling's approximation.
\end{Proof}

\subsection{Derivation of CLT for $Q$}
\label{sec-p0:se:CLT:proof}

We first justify our CLT application in \cref{res-p0:se:CLT}.
The distribution of $Q_i(t_\alpha)$ depends on $k$ (and $n$), and so we cannot apply the standard CLT. Instead, we apply a CLT for `triangular arrays'; specifically, we now state a special case of the Lindeberg--Feller theorem.

\begin{thm}[CLT for Triangular Arrays; cf {\cite[Theorem~3.4.5]{D:prob-te}}]
\label{res-p0:CLT-triangular}
	For each $k \in \mbn$, let $\{Y_{i,k}\}_{i=1}^k$ be an iid sequence of centralised, normalised random variables, and suppose that $\ex{Y_{1,k}^4} \ll k$.
	Then
	\[
		\sumt[k]{i=1} Y_{i,k}/\sqrt k
	\to^d
		N(0,1)
	\quad
		\askinf,
	\]
	where $N(0,1)$ is a standard Gaussian random variable.
\end{thm}

This version can be deduced easily from the version given in \textcite[Theorem~3.4.5]{D:prob-te}.
Indeed, apply \cite[Theorem~3.4.5]{D:prob-te} to the iid triangular array defined by
\(
	X_{i,k} \cq Y_{i,k} / \sqrt k.
\)
Note~that
\[
	\rbb{ \sumt[k]{1} \ex{ X_{i,k}^2 \one{ \abs{X_{i,k}} \ge \eps} } }^2
&
=
	\ex{ Y_{1,k}^2 \one{ \abs{Y_{1,k}} \ge \eps \sqrt k } }^2
\\&
\le
	\ex{ Y_{1,k}^4 } \pr{ \abs{Y_{1,k}} \ge \eps \sqrt k }
\le
	\ex{ Y_{1,k}^4 } / (\eps^2 k)
\to
	0.
\]

Using this CLT for triangular arrays, we can deduce a CLT for $Q$.

\begin{Proof}[Proof of \cref{res-p0:se:CLT}]
For our application of \cref{res-p0:CLT-triangular},
for each $\alpha \in \mbr$,
we take
\[
	Y_{i,k}
\cq
	Y_{i,k}(\alpha)
\cq
	\frac{ Q_i(t_\alpha) - \ex{ Q_i(t_\alpha) } }{ \sqrt{ \Var{ Q_i(t_0) } } }.
\label{eq-p0:se:CLT:Ydef}
\nt
\]
Observe that $\ex{Y_{i,k}} = 0$ and $\Var{Y_{i,k}} = \ex{Y_{i,k}^2} = 1$.
Assuming that $\ex{Y_{i,k}^4} \ll k$,
we deduce the following result:
for any sequence $(\alpha_n)_{n\ge1}$ which converges to $\alpha$, we deduce that
\[
	\pr{ Q(t) - \ex{Q(t)} \ge \alpha_n \sqrt{\Var{Q(t)}} }
\to
	\Psi(\alpha).
\label{eq-p0:se:CLT:entropic-general}
\nt
\]
(We are also using Slutsky's theorem to allow $\alpha_n$ to depend on $n$, and, of course, the fact that $\kinf$ as $\ninf$.)
We also further rely on the following claim:
\[
	\text{if}
\quad
	t \eqsim t_0,
\Quad{then}
	\Varb{Q_1(t)} \eqsim \Varb{Q_1(t_0)};
\Quad{also}
	\Varb{Q(t_0)} \gg 1.
\label{eq-p0:se:var-sim}
\nt
\]
We prove these two statements in this claim (independently of the proof of the CLT) in \cref{res-p0:se:var:t0} in \S\ref{sec-p0:se:var}.
Now recall \cref{eq-p0:se:t0a}, which says that $t_\alpha \eqsim t_0$ for all $\alpha \in \mbr$.
Taking
\[
	\alpha_n \cq - \alpha \sqrt{ \Var{Q(t_0)} / \Var{Q(t_\alpha)} } \pm \omega / \sqrt{ \Var{ Q(t_\alpha) } }
\Quad{with}
	\omega \cq \Varb{ Q(t_0) }^{1/4} \gg 1,
\]
applying (\ref{eq-p0:se:CLT:entropic-general}, \ref{eq-p0:se:var-sim}) along with the above recollections we obtain the desired result:
\[
	\pr{ Q(t_\alpha) \le \log n \pm \omega }
\to
	\Psi(\alpha).
\]

It remains to verify that $\ex{Y_{i,k}^4} \ll k$.
Roughly, $\abs{ W_1(t) }$ is `well-approximated' by the following:
\begin{alignat*}{2}
	&\absb{ N\rbb{\ext{W_1(t)}, \: t/k} }&
		\Qwhere t/k &\gg 1, \Quad{ie} k \ll \log n;
\\
	&\Bernoulli(t/k)&
		\Qwhere t/k &\ll 1, \Quad{ie} k \gg \log n.
\end{alignat*}
In the interim regime $k \asymp \log n$, we have that $W_1$ behaves like an `order 1' random variable, in the sense that its mean and variance are bounded away from both 0 and $\infty$.
It will actually turn out that the normal approximation is sufficient in the $k \asymp \log n$ regime also.
Below, we abbreviate $Q_1(t_\alpha)$ by $Q_1$, $W_1(t_\alpha)$ by $W_1$ and $t_\alpha$ by $t$.

Write $s \cq t/k$.
We consider separately the cases $s \gtrsim 1$ and $s \ll 1$. When $s \gtrsim 1$, we have $t \gtrsim k \gg 1$; when considering $s \ll 1$, however, we shall only consider $t$ with $1 \ll t \ll k$.
We shall be interested in $t \cq t_\alpha \eqsim t_0$, and \cref{res-p0:se:t0a} says that $t_0 \gg 1$ in all regimes; hence we need only consider $t \gg 1$.
Let $\delta > 0$ be some (arbitrarily) small number.

\medskip

\emph{Consider first $s = t/k$ with $s \ge \delta$.}
In this regime, we approximate $W_1(t)$ by a $N(\ex{W_1},s)$ distribution, where $s = t/k$. Let $Z \sim N(\ex{W_1},s)$, and write $f$ for its density function:
\[
	f(x) \cq (2 \pi s)^{-1/2} \expb{ - \tfrac1{2s} (x-\ext{W_1})^2 }
\Qfor
	x \in \mbr.
\label{eq-p0:se:CLT:f-def}
\nt
\]
Let $R_1$ be a real valued random variable defined so that
\[
	R_1 = - \log f(x)
\Qwhere
	W_1 = x.
\label{eq-p0:se:CLT:R-def}
\nt
\]
Also write $G \cq W_1 + U$, where $U \sim \Unif[-\tfrac12,\tfrac12)$ is independent of $W_1$; then $G$ has density function
\[
	g(x) \cq \pr{ W_1 = [x] }
\Qfor
	x \in \mbr,
\label{eq-p0:se:CLT:g-def}
\nt
\]
where $[x] \in \mbz$ is $x \in \mbr$ rounded to the nearest integer (rounding up when $x \in \mbz + \tfrac12$).
We have
\[
	(a-b)^4
\le
	3^4 \rbb{ (a-a')^4 + (a'-b')^4 + (b' - b)^4 }
\Qforall
	a,a',b,b' \in \mbr.
\]
Applying this inequality with $a = Q_1$, $a' = R_1$, $b = \ext{Q_1}$ and $b' = \ext{R_1}$, we obtain
\[
	3^{-4}
	\ex{ (Q_1 - \ext{Q_1})^4 }
&
\le
	\ex{ (Q_1 - R_1)^4 }
+	\ex{ (R_1 - \ext{R_1})^4 }
+	\ex{R_1 - Q_1}^4
\\&
\le
	\ex{ (R_1 - \ext{R_1})^4 }
+	2 \, \ex{ (Q_1 - R_1)^4 },
\label{eq-p0:se:CLT:norm:ineq}
\nt
\]
with the second inequality following from Jensen (or Cauchy--Schwarz twice).
We study these terms separately. Approximately, the local CLT will say that the second term is small; up to an error term which we control with the local CLT, the first term we can calculate directly using properties of the normal distribution.

We consider first the first term of \cref{eq-p0:se:CLT:norm:ineq}.
In terms of an integral, it is given by
\[
	\ex{ (R_1 - \ext{R_1})^4 }
=
	\intt\mbr g(x) \rbb{ - \log f(x) - \ext{R_1} }^4 \, dx.
\]
The local CLT suggests that we can approximately replace the $g(x)$ factor by $f(x)$, at least for a large range of $x$. So let us first study
\[
	\intt\mbr f(x) \rbb{ - \log f(x) - \ext{R_1} }^4 \, dx
=
	\intt\mbr f(x+\ext{W_1}) \rbb{ - \log f(x+\ext{W_1}) - \ext{R_1} }^4 \, dx.
\]
A direct calculation reveals, remarkably, that the last expression is independent of the mean of $W_1$---this is a property special to the family of normal distributions.
Expanding the fourth power and using moments of $N(0,1)$, one finds that this equals $\tfrac{15}4$; the exact numerical value is unimportant.

Now, by the local CLT \cref{eq-p0:LCLT}, we have
\[
	\intt[s^{7/12}]{-s^{7/12}}
	g(x) \rbb{ - \log f(x) - \ext{R_1} }^4 \, dx
&
=
	\rbb{1 + \Oh{s^{-1/4}}}
	\intt[s^{7/12}]{-s^{7/12}}
	f(x) \rbb{ -\log f(x) - \ext{R_1} }^4 \, dx
\\&
\le
	\rbb{1 + \Oh{s^{-1/4}}} \cdot \tfrac{15}4.
\]
Using bounds on the tail of the SRW and Poisson distribution, as given in \cref{res-p0:rp:lde:poi,res-p0:rp:lde:srw}, it is straightforward to see, in both the undirected and directed cases, that
\[
	\intt{\mbr \setminus [-s^{7/12},s^{7/12}]}
	f(x) \rbb{ -\log f(x) - \ext{R_1} }^4 \, dx
=
	\ohb{s^{-10}}.
\label{eq-p0:se:CLT:tail-integrals}
\nt
\]
(In fact, it is easy to see that it is $\Oh{\exp{-cs^{1/6}}}$ for some sufficiently small constant $c$.)
Hence
\[
	\ex{ (R_1 - \ext{R_1})^4 }
=
	\tfrac{15}4 \rbb{ 1 + \Oh{s^{-1/4}} }
=
	\tfrac{15}4 \rbb{1 + \oh1}.
\label{eq-p0:se:CLT:norm:term1}
\nt
\]

We now turn to the second term of \cref{eq-p0:se:CLT:norm:ineq}.
In terms of an integral, it is given by
\[
	\ex{ (Q_1 - R_1)^4 }
=
	\ex{ \rbb{ \log f(W_1) - \log g(W_1) }^4 }
=
	\intt\mbr g(x) \log\rbb{ f(x)/g(x) }^4 \, dx.
\]
Again by the local CLT \cref{eq-p0:LCLT}, we have
\[
	\intt[s^{7/12}]{-s^{7/12}}
	g(x) \log\rbb{ f(x)/g(x) }^4 \, dx
=
	\Ohb{s^{-1/4}}
	\intt[s^{7/12}]{-s^{7/12}}
	g(x) \, dx
\le
	\Ohb{s^{-1/4}},
\]
and a similar application of the tail bounds in \cref{res-p0:rp:lde:poi,res-p0:rp:lde:srw} shows that
\[
	\intt{\mbr \setminus [-s^{7/12},s^{7/12}]}
	g(x) \log\rbb{ f(x)/g(x) }^4 \, dx
=
	\ohb{s^{-10}}
=
	\Ohb{s^{-1/4}}.
\label{eq-p0:se:CLT:norm:term2}
\nt
\]

Hence, combining (\ref{eq-p0:se:CLT:norm:term1}, \ref{eq-p0:se:CLT:norm:term2}) into \cref{eq-p0:se:CLT:norm:ineq}, we obtain
\[
	\ex{ (Q_1 - \ext{Q_1})^4 }
\le
	\tfrac{15}4 \cdot 3^4 + \oh1
\le
	1000.
\]

\smallskip

Now consider $\Var{Q_1}$.
The arguments used for \cref{res-p0:se:var:s} (in \S\ref{sec-p0:se:var}) below show that
\[
	\ex{ \rbb{ Q_1(sk) - \ex{Q_1(sk)} }{}^4 }
\lesssim
	1
\Qwhere
	s \gtrsim 1.
\]
Since $s_0 = t_0/k \gtrsim 1$ in the regime $k \lesssim \log n$,
we deduce that
\(
	\ex{Y_{1,k}^4}
\lesssim
	1
\ll
	k.
\)

\medskip

\emph{Consider now $s = t/k$ with $s \le \delta$ but $t \gg 1$.}
In this regime, we approximate the number of steps taken by $\Bernoulli(t/k)$. Indeed, we have
\[
	\ex{ W_1 = 0 } = 1 - s + \Oh{bs^2}
\Qand
	\ex{ \abs{W_1} = 1 } = s + \Ohb{s^2}.
\]
We also use the fact that, for both the undirected and directed cases, for $x \ge 0$ we have
\[
	\pr{ W_1 = x }
\ge
	\pr{ \Poisson(s) = x } \cdot 2^{-x}
=
	2^{-x} e^{-s} s^x / x!
\ge
	(s^2/x)^x;
\label{eq-p0:se:CLT:poisson-lower}
\nt
\]
from this one deduces that $-\log\pr{W_1 = x} \le x \log(x/s^2) = x (x + 2 \log(1/s))$.
We use this to show that the terms with $\abs x \ge 2$ contribute subleading order to the expectation
\[
	\ex{ Q_1 }
=
	\sumt{x} \pr{ W_1 = x } \log 1/\pr{ W_1 = x }
=
	s \log(1/s) + \Oh{s}.
\]
Similarly, we can use \cref{eq-p0:se:CLT:poisson-lower} to ignore the terms with $\abs x \ge 2$ in
\[
	\ex{ \absb{ Q_1 - \ext{Q_1} }^r }
&
=
	\sumt{x}
	\pr{ W_1 = x }
	\absb{ - \log\pr{ W_1 = x } - s\log(1/s) + \Oh{s} }^r
\\&
=
	s \log(1/s)^r \rbb{ 1 + \Oh{s} },
\label{eq-p0:se:CLT:Q1^r}
\nt
\]
for any fixed $r \in \mbn$ with $r \ge 2$, say $r \in \{2,3,4\}$.

In particular, this says that $\Var{Q_1} \eqsim s \log(1/s)^2$, and so
\[
	\ex{ Y_{i,k}^4 }
\eqsim
	\rbb{ s \log(1/s)^4 } \big/ \rbb{ s \log(1/s)^2 }^2
=
	1/s
=
	k/t
\ll
	k,
\]
with the final relation holding since while $s \ll 1$ we do have $t \gg 1$.
\end{Proof}

We now have all that we need to get on and calculate the entropic time $t_0$ in the three regimes of $k$. However, in order to find the cutoff times $t_\alpha$, we need to know what the variance of the terms in the sum $Q(t)$, ie $\Var{Q_1(t)}$, is for $t \eqsim t_0$.

\subsection{Variance of $Q_1(t)$}
\label{sec-p0:se:var}

Recall that, for all $t \ge 0$, we have
\[
	Q(t)
=
	-\log \mu(t)
=
	-\sumt[k]{i=1} \log \nu_t\rbb{ W_i(t) }
=
	\sumt[k]{i=1} Q_i(t),
\]
and that the $Q_i(t)$-s are iid (for fixed $t$).
We now determine what its variance is at the entropic time $t_0$, and how the variance changes around this time.
Note that $\Var{Q(t)} = k \Var{Q_1(t)}$.

\begin{prop}
\label{res-p0:se:var:s}
In both the undirected and the directed case,
\begin{subequations}
	\label{eq-p0:se:var:s}
\begin{empheq}[left = {%
	\Varb{Q_1(sk)}
\eqsim
	\empheqlbrace}]%
{alignat=2}
	&1/2
		&&\Quad{as} s \to \infty,
	\label{eq-p0:se:var:s>1}
\\
	&s \log(1/s)^2
		&&\Quad{as} s \to 0;
	\label{eq-p0:se:var:s<1}
\end{empheq}
\end{subequations}
furthermore, the map $s \mapsto \Var{Q_1(sk)} : [0,\infty) \to \mbr_+$ is continuous.
\end{prop}

From this, it is easy to calculate the variance at the entropic time $t_0$. Note that knowledge of the variance \emph{is not} required to calculate $t_0$.
The reader should recall \cref{res-p0:se:t0a}.
(Knowledge of the variance is required in calculating $t_\alpha$ with $\alpha \ne 0$, but not with $\alpha = 0$.)

\begin{cor}
\label{res-p0:se:var:t0}
For all regimes of $k$,
in both the undirected and directed case,
\[
	\text{if}
\quad
	t \eqsim t_0,
\Quad{then}
	\Varb{ Q_1(t) } \eqsim \Varb{ Q_1(t_0) } \gg 1/k.
\label{eq-p0:se:var:sim}
\nt
\]
Moreover,
for all $\lambda > 0$,
we have
\begin{subequations}
	\label{eq-p0:se:var:k}
\begin{empheq}[left = {%
	\Varb{Q_1(t_0)}
\eqsim
	\empheqlbrace}]%
{alignat=2}
	&1/2
		&&\Qwhere k \ll \log n,
	\label{eq-p0:se:var:k<logn}
\\
	&v(\lambda)
		&&\Qwhere k \eqsim \lambda \log n,
	\label{eq-p0:se:var:k=logn}
\\
	&\log n \log(k/\log n) / k
		&&\Qwhere k \gg \log n,
	\label{eq-p0:se:var:k>logn}
\end{empheq}
\end{subequations}
where $v : (0,\infty) \to (0,\infty) : \lambda \mapsto \Var{Q_1(f(\lambda) k)}$ is a continuous function whose value differs between the undirected and directed cases.
\end{cor}

\begin{Proof}[Proof of \cref{res-p0:se:var:t0}]
The first claim is immediate from \cref{res-p0:se:var:s}.
The claim for $k \ll \log n$ also uses \cref{eq-p0:se:t0a}.
For $k \gg \log n$, there is a small amount of work to do.
Set $s_0 \cq t_0/k$, and so
\[
	s_0 = \frac{t_0}{k} \eqsim \frac{\log n/k}{\log(k/\log n)} = \frac1{\kappa\log\kappa}
\Qwhere
	\kappa \cq \frac{k}{\log n} \gg 1.
\]
We then also have
\[
	\log(1/s_0)
=
	\log\log\kappa + \log\kappa + \oh 1
\eqsim
	\log\kappa,
\]
and hence
\[
	s_0 \log(1/s_0)^2
\eqsim
	(\log\kappa)^2 / \rbr{ \kappa \log\kappa }
=
	\log\kappa / \kappa
=
	\log n \log(k/\log n) / k.
\]
Note that while this has $\Var{Q_1(t_0)} \ll 1$, it does have $\Var{Q(t_0)} = k \Var{Q_1(t_0)} \gg 1$.

Finally consider $k \eqsim \lambda \log n$.
Each coordinate runs at rate $1/k$, so for all $s \in \mbr_+$ the map $s \mapsto \Var{Q_1(sk)}$ is continuous.
Hence given $C > 0$ there exists an $M$ so that
\[
	1/M \le \Varb{ Q_1(sk) } \le M
\Qforall
	s \text{ with } 1/C \le s \le C.
\]
By \cref{eq-p0:se:t0a:k==logn}, we have $s = t_0/k \to f(\lambda)$.
Hence $\Var{Q_1(t_0)} \to v$ for some constant $v \in (0,\infty)$ depending only on $\lambda$.
This $v$ is not the same in the directed and undirected cases.
\end{Proof}

\begin{Proof}[Proof of \cref{res-p0:se:var:s}]
Consider first $s \to \infty$.
This proof is similar to the $s \gtrsim 1$ case, in justifying the CLT application.
In particular, if
\[
	g(x) \cq \pr{ W_1(sk) = [x] }
\Qand
	f(x) \cq (2 \pi s)^{-1/2} \expb{ - \tfrac1{2s} (x-\ex{W_1(sk)})^2 },
\]
then the local CLT \cref{eq-p0:LCLT} says, for $s \gtrsim 1$, that
\[
	g(x) = f(x) \, \rbb{ 1 + \Oh{s^{-1/4}} }
\Qfor
	x \in \mbr
\Quad{with}
	\absb{ x - \ex{W_1(sk)} } \le s^{7/12}.
\]
Under the assumption that $W_1(sk)$ is actually distributed as $N(0,s)$, direct calculation as in the previous section shows that the variance is then $\tfrac12$.
Considering the same approximations as before, namely splitting the integration range into $\abs{x - \ex{W_1}} \le s^{7/12}$ and $\abs{x - \ex{W_1}} > s^{7/12}$, and using the local CLT to argue that $\log(g(x)/f(x)) = \Oh{s^{-1/4}}$ for $x$ in the first range, we obtain
\[
	\Varb{ Q_1(sk) } = \tfrac12 + \Ohb{ s^{-1/4} \log s } \eqsim \tfrac12
\Qwhere
	s \gg 1.
\]

Consider next $s \to 0$.
In the CLT justification in the case $s \gtrsim 1$, we showed that
\[
	\ex{ \absb{ Q_1(sk) - \ext{Q_1(sk)} }^r }
=
	s \log(1/s)^r + \Ohb{s^2 \log(1/s)^r},
\tag{\ref{eq-p0:se:CLT:Q1^r}}
\]
and in particular deduced that $\Var{Q_1(sk)} \eqsim s \log(1/s)^2$.
This applies for $s \ll 1$ also.

The continuity of $s \mapsto \Var{Q_1(sk)}$ follows from the dominated convergence theorem.
\end{Proof}

\subsection{Calculating the Entropic and Cutoff Times}
\label{sec-p0:se:entropic-time-calc}

In this section we calculate the entropic time $t_0$, and the cutoff times $t_\alpha$.
Recall that
\[
	h(t) = \ex{ Q(t) }
\Qand
	H(s) \cq \ex{ Q_1(sk) };
\]
note that $H(s)$ is the entropy of $W_1(sk)$, which forms a rate-1 RW on $\mbz$.
The primary purpose of this section is to prove \cref{res-p0:se:t0a}, which the reader should recall.
To prove this, we derive asymptotic expressions for the entropy of rate-1 RW on $\mbz$.
As a consequence of these, one can see by how much the entropy changes when $t_0$ is replaced by $(1 + \xi) t_0$ for a (small)
constant $\xi \in \mbr$.

\begin{lem}
\label{res-p0:se:entropic-time-calc:+omega}
	Let $t_0$ and $t_{\pm2\omega}$ be the times at which the entropy of rate-1 RW on $\mbz^k$ obtains entropy $\log N$ and $\log N \pm 2 \omega$, respectively.
	Then $t_{\pm2\omega} \eqsim t_0$ if $\omega \ll \min\bra{k, \log N}$.
\end{lem}

We give this proof straight away, quoting results which are proved in the upcoming subsections.

\begin{Proof}[Proof of \cref{res-p0:se:entropic-time-calc:+omega}]
We prove the claim for $t_{2\omega}$; the analysis for $t_{-2\omega}$ is identical.

For $s \ge 0$, write $H(s)$ for the entropy rate-1 RW on $\mbz$ evaluated at time $s$.
In \cref{res-p0:se:ent:s>1,res-p0:se:ent:s<1} below,
we establish the following relations:
\begin{subequations}
	\label{eq-p0:ent:approx}
\begin{empheq}%
[ left = { H(s) = \empheqlbrace} ]%
{alignat = 2}
	&\tfrac12 \log(2 \pi e s) + \Oh{ s^{-1/4} }
	&&\Qwhere s \gg 1;
	\label{eq-p0:ent:approx:<}
\\
	&s \log(1/s) + \Oh{s}
	&&\Qwhere s \ll 1;
	\label{eq-p0:ent:approx:>}
\end{empheq}
\end{subequations}
when $s \asymp 1$, we use continuity of $\lambda \mapsto H(1/\lambda)$ (cf \cref{res-p0:se:ent:s=1|k=logn} below).
Let $\delta > 0$; assume that $\delta \downarrow 0$, but more slowly than the error terms in \cref{eq-p0:ent:approx}.
Recall that the entropic time $t_0 = s_0 k$ is defined so that $H(s_0) = \log N / k$.
We need to choose $\delta$ and $\omega$ so that
\[
	H\rbb{ s_0 (1+\delta) }
\ge
	\log N/k + 2\omega / k
\Quad{with}
	\text{$\delta \ll 1$ and $\omega \gg 1$}.
\]
Note that the statement is monotone in $\omega$:
	if it holds for some $\omega$, then it holds for any $0 \le \omega' \le \omega$, since then $t_0 \le t_{2\omega'} \le t_{2\omega}$.
Hence we may assume lower bounds on $\omega$, if desired.

\smallskip

\emph{Regime $k \ll \log N$.}
	%
We have $s_0 \gg 1$.
By \cref{eq-p0:ent:approx:<},
we have
\[
	H\rbb{ s_0(1+\delta) }
&
=
	\tfrac12 \log(2 \pi e s_0) + \tfrac12 \log(1+\delta) + \Ohb{s_0^{-1/4}}
\\&
=
	\log N / k + \tfrac12 \delta + \Ohb{ \min\brb{ \delta^2, s_0^{-1/4} } }.
\]
We take $\delta \cq 5 \omega/k$, and so need $s_0^{-1/4} \ll \omega/k \ll 1$.
Hence $\omega \ll k$ suffices.

\smallskip

\emph{Regime $k \gg \log N$.}
	%
We have $s_0 \ll 1$.
By \cref{eq-p0:ent:approx:>},
we have
\[
	H\rbb{ s_0(1+\delta) }
&
=
	(1+\delta) \cdot s_0 \log(1/s_0)
-	s_0 (1+\delta) \log(1+\delta)
+	\Oh{s_0}
\\&
=
	\log N / k + \delta \log N / k + \Oh{s_0}.
\]
We take $\delta \cq 2 \omega / \log N$, and so need $s_0 \ll \omega / \log N \ll 1$.
Hence $\omega \ll \log N$ suffices.

\smallskip

\emph{Regime $k \asymp \log N$.}
By continuity and the strict increasing property of the entropy, all we require is that $\log N / k + 2\omega / k = (1+\oh1) \log N / k$, and hence only require $\omega \ll \log N \asymp k$.
\end{Proof}

\subsubsection{Regime $k \ll \log n$}

We first consider the regime $k \ll \log n$, which corresponds to $s_0 = t_0/k \gg 1$.

\begin{prop}
\label{res-p0:se:ent:s>1}
	For $s \gtrsim 1$,
	the entropy $H$ of a rate-1 SRW or DRW on $\mbz$ satisfies
	\[
		H(s) = \tfrac12 \log(2 \pi e s) + \Ohb{ s^{-1/4} }.
	\label{eq-p0:se:ent:s>1}
	\nt
	\]
\end{prop}

\begin{Proof}
We consider both the directed and undirected cases together.
Write $t \cq sk$.
Define $f$, $R_1$ and $g$ as in
(\ref{eq-p0:se:CLT:f-def}, \ref{eq-p0:se:CLT:R-def}, \ref{eq-p0:se:CLT:g-def}),
respectively.
By \cref{eq-p0:se:CLT:norm:term2}, we have
\[
	\absb{ \ext{Q_1} - \ext{R_1} }
\le
	\ex{ (Q_1 - R_1)^4 }^{1/4}
=
	\ohb{s^{-5/2}}
=
	\Ohb{s^{-1/4}}
\Qwhere
	s \gtrsim 1.
\]
Direct calculation with its pdf shows that the entropy of $N(0,s)$ is precisely $\tfrac12 \log(2 \pi e s)$.
Using this along with a similar calculation as used for \cref{eq-p0:se:CLT:norm:term1} gives
\[
	\ex{ R_1 }
=
	\rbb{1 + \Oh{s^{-1/4}}} \cdot \log(2 \pi e s).
\]
Hence we obtain our desired expression, namely \cref{eq-p0:se:ent:s>1}.
\end{Proof}

We now calculate the derivative of this entropy.

\begin{prop}
\label{res-p0:se:der:s>1}
	For $s \gtrsim 1$,
	the entropy $H$ of a rate-1 SRW or DRW on $\mbz$ satisfies
	\[
		H'(s) = (2s)^{-1} \rbb{ 1 + \Ohb{ s^{-10} } }.
	\label{eq-p0:se:der:s>1}
	\nt
	\]
	By the chain rule,
	for $t \gtrsim k$,
	the entropy $h$ of rate-1 SRW or DRW on $\mbz^k$ then satisfies
	\[
		h'(t) = H'(t/k) = (2t/k)^{-1} \rbb{ 1 + \Ohb{ (t/k)^{-10} } }.
	\]
\end{prop}

\begin{Proof}
Write $t \cq sk$.
Define $f$, $R_1$ and $g$ as in
(\ref{eq-p0:se:CLT:f-def}, \ref{eq-p0:se:CLT:R-def}, \ref{eq-p0:se:CLT:g-def})
respectively.
We have
\[
	H(s) = - \sumt{x \in \mbz} \pr{ X_s = x } \log\pr{ X_s = x }.
\]
Differentiating this with respect to $t$ we obtain
\[
	H'(s)
=
	-\sumt{x \in \mbz} \tfrac{d}{ds} \pr{ X_s = x } \rbb{ \log\pr{ X_s = x } + 1 }
=
	-\sumt{x \in \mbz} \tfrac{d}{ds} \pr{ X_s = x } \cdot \log \pr{ X_s = x }.
\]

Consider first the SRW.
Using the Kolmogorov backward equations for the SRW,
we obtain
\[
	\tfrac{d}{ds} \pr{ X_s = x }
=
	\tfrac12 \pr{ X_s = x+1 }
+	\tfrac12 \pr{ X_s = x-1 }
-	\pr{ X_s = x }.
\]
Recall that $\nu_s(x) \cq \pr{ X_s = x }$; write $g_s(x) \cq \nu_s([x])$.
Since $\sum_{x \in \mbz} \nu_s(x) = 1$, we obtain
\begin{subequations}
	\label{eq-p0:se:der:s>1:srw:integrals}
\[
	H'(s)
&
=
	\sumt{x \in \mbz}
	\rbb{ \nu_s(x) - \tfrac12( \nu_s(x+1) + \nu_s(x-1) ) } \log \nu_s(x)
\\&
=
	\intt\mbr
	\rbb{ g_s(x) - \tfrac12( g_s(x+1) + g_s(x-1) ) } \log g(x).
\\&
=
	\intt\mbr
	\rbb{ g_s(x) - \tfrac12( g_s(x+1) + g_s(x-1) ) } \log f_s(x) \, dx
\label{eq-p0:se:der:s>1:srw:main}
\nt
\\&\qquad+
	\intt\mbr
	\rbb{ g_s(x) - \tfrac12( g_s(x+1) + g_s(x-1) ) } \log\rbb{ g_s(x)/f_s(x) } \, dx,
\label{eq-p0:se:der:s>1:srw:error}
\nt
\]
\end{subequations}
where $f_s(x) \cq (2 \pi s)^{-1/2} \exp{-x^2/(2s)}$.
The same arguments as used for \cref{eq-p0:se:CLT:tail-integrals} show that the integral in \cref{eq-p0:se:der:s>1:srw:error} is $\oh{s^{-10}}$.
Now consider the integral in \cref{eq-p0:se:der:s>1:srw:main}.
Using a simple shift, we have
\[
	\intt\mbr g_s(x+1) \log f_s(x) \, dx
=
	\intt\mbr g_s(x) \log f_s(x) \, dx
-	\intt\mbr g_s(x) \log\rbb{ f_s(x-1)/f_s(x) } \, dx,
\]
and we consider $\int_\mbr g_s(x-1) \log f_s(x) \, dx$ similarly; hence we have
\GAP{3}
\[
&	\intt\mbr \rbb{ g_s(x) - \tfrac12( g_s(x+1) + g_s(x-1) ) } \log f_s(x) \, dx
\\&\gap
=
	\tfrac12 \intt\mbr g_s(x) \rbb{ \log\rbb{ f_s(x-1)/f_s(x) } + \log\rbb{ f_s(x+1)/f_s(x) } } \, dx
\\&\gap
=
	\tfrac12 \intt\mbr g_s(x) \log\rbb{ f_s(x-1) f_s(x+1) / f_s(x)^2 } \, dx.
\]
Since $f_s(x) = (2 \pi s)^{-1/2} \exp{-x^2/(2s)}$, this log is precisely $1/s$ (independent of $x$). Since it is a distribution, $g_s$ integrates to 1, so the integral equals $1/(2s)$.
This proves the SRW case.

Now consider the DRW.
Here the backward Kolmogorov equations read
\[
	\tfrac{d}{ds} \pr{ X_s = x }
=
	\pr{ X_s = x-1 } - \pr{ X_s = x }
\Qfor
	x \in \mbn
\]
and
\(
	\tfrac{d}{ds} \pr{ X_s = 0 } = - \pr{ X_s = 0 } = - e^{-s}.
\)
Hence, as above, we have
\begin{subequations}
	\label{eq-p0:se:der:s>1:poi:integrals}
\[
	H'(s) - s e^{-s}
&
=
	\sumt{x \in \mbn}
	\rbb{ \nu_s(x) - \nu_s(x-1) } \log \nu_s(x)
\\&
=
	\intt[\infty]{1/2} \rbb{ g_s(x) - g_s(x-1) } \log f_s(x) \, dx
\label{eq-p0:se:der:s>1:poi:main}
\nt
\\&\qquad+
	\intt[\infty]{1/2} \rbb{ g_s(x) - g_s(x-1) } \log\rbb{ g_s(x)/f_s(x) } \, dx.
\label{eq-p0:se:der:s>1:poi:error}
\nt
\]
\end{subequations}
As for \cref{eq-p0:se:der:s>1:srw:error} above, the same arguments as used for \cref{eq-p0:se:CLT:tail-integrals} show that the integral in \cref{eq-p0:se:der:s>1:poi:error} is $\oh{s^{-10}}$.
Note also that $s e^{-s} = \oh{s^{-10}}$ as $s \to \infty$.
Now consider the integral in \cref{eq-p0:se:der:s>1:poi:main}.
Using a simple shift as before, we have
\[
	\intt[\infty]{1/2} \rbb{ g_s(x) - g_s(x-1) } \log f_s(x) \, dx
&
=
	-\intt[\infty]{1/2} g_s(x) \log\rbb{ f_s(x+1) / f_s(x) } \, dx
\\&
=
	\intt[\infty]{1/2} g_s(x) \rbb{ (x-s)/s + 1/(2s) } \, dx
=
	1/(2s),
\]
recalling that here $f_s(x) = (2 \pi s)^{-1/2} \exp{ - (x-s)^2/(2s) }$, $\ex{X_s} = s$ and $g_s$ integrates to 1.
In the same way as for the SRW, this proves the DRW case.
\end{Proof}

We wish to find the times $s_\alpha = t_\alpha/k$ defined so that, recalling \cref{eq-p0:se:var:k<logn},
\[
	H(s_\alpha) = \rbb{ \log n + \alpha \sqrt{v k} }/k
\Qwhere
	v \cq \Varb{ Q_1(t_0) } \eqsim \tfrac12.
\]

\begin{prop}
\label{res-p0:se:ent:k<logn}
	For $k \ll \log n$, we have
	\[
		s_0
	=
		t_0/k
	\eqsim
		n^{2/k} / (2 \pi e),
	\tag{\ref{eq-p0:se:t0a:k<<logn}}
	\]
	and,
	for each $\alpha \in \mbr$,
	we have
	$s_\alpha \eqsim s_0$,
	and furthermore
	\[
		\rbr{ s_\alpha - s_0 }/s_0
	=
		\rbr{ t_\alpha - t_0 }/t_0
	\eqsim
		\alpha \sqrt{2/k} = \oh1.
	\tag{\ref{eq-p0:se:t0a:k<<logn}}
	\]
\end{prop}

\begin{Proof}
We consider the directed and undirected cases simultaneously.
By directly manipulating \cref{eq-p0:se:ent:s>1}, we see that
if $H(s_0) = \log n/k$ then
\[
	s_0
=
	n^{2/k} / (2 \pi e) \cdot \rbb{ 1 + \Ohb{ s_0^{-1/4} } }
\eqsim
	n^{2/k} / (2 \pi e),
\]
noting that $k \ll \log n$ and so $n^{2/k} \gg 1$.
This proves the first part of \cref{eq-p0:se:t0a:k<<logn}.

We now turn to finding $t_\alpha$.
Fix $\alpha \in \mbr$.
Note that $H$ is increasing and $\alpha \sqrt{v/k} = \oh1$.
So from the form of $H(s)$ given in \cref{eq-p0:se:ent:s>1} we see that, for all $\eps > 0$, we have $(1 - \eps) s_0 \le s_\alpha \le (1 + \eps) s_0$ for $n$ sufficiently large (depending on $\alpha$);
hence $s_\alpha \eqsim s_0$ for all $\alpha \in \mbr$.

By definition of $s_\alpha$, we have
\[
	H(s_\alpha) - H(s_0) = \alpha \sqrt{v/k},
\Quad{and hence}
	\tfrac{ds_\alpha}{d\alpha} \, H'(s_\alpha) = \sqrt{v/k}.
\]
Hence we have
\[
	s_\alpha - s_0
=
	\intt[\alpha]{0} \tfrac{ds_a}{da} \, da
=
	\sqrt{v/k} \intt[\alpha]{0} 1/H'(s_a) \, da.
\]
But, by \cref{res-p0:se:der:s>1}, we may write $H'(s) = (2s)^{-1} \rbr{ 1 + \oh1 }$ with $\oh1$ term uniform over $s \in [\tfrac12 s_0, 2 s_0]$, which is an interval containing the cutoff window.
Hence, recalling from \cref{eq-p0:se:var:k<logn} that $v \eqsim \tfrac12$ in this regime, the second part of \cref{eq-p0:se:t0a:k<<logn} follows:
\[
	s_\alpha - s_0
\eqsim
	2 \alpha s_0 \sqrt{1/2} / \sqrt k
\eqsim
	\alpha s_0 \sqrt{2/k}.
\qedhere
\]
\end{Proof}

\subsubsection{Regime $k \asymp \log n$}

We next consider the regime $k \eqsim \lambda \log n$ with $\lambda \in (0,\infty)$, which corresponds to $s_0 = t_0/k \asymp 1$.

\begin{prop}
\label{res-p0:se:ent:s=1|k=logn}
	There exists a decreasing, continuous bijection $f : (0,\infty) \to (0,\infty)$, whose value differs between the undirected and directed cases, so that,
	for all $\lambda > 0$,
	for $k \eqsim \lambda \log n$,
	we have
	\[
		s_0
	=
		t_0/k
	\eqsim
		f(\lambda)
	\Qwhere
		f(\lambda) \cq H^{-1}(1/\lambda),
	\tag{\ref{eq-p0:se:t0a:k==logn}}
	\]
	and,
	for each $\alpha \in \mbr$, we have
	$s_\alpha \eqsim s_0$,
	and furthermore
	\[
	\begin{gathered}
		\rbr{ s_\alpha - s_0 }/s_0
	=
		\rbr{ t_\alpha - t_0 }/t_0
	\eqsim
		\alpha g(\lambda) / \sqrt k
	=
		\oh1
	\\
		\text{where}
	\quad
		g(\lambda) \cq \sqrt{ \Var{ Q_1( f(\lambda) \, k ) } } / \rbb{ f(\lambda) H'\rbr{ f(\lambda) } }.
	\end{gathered}
	\tag{\ref{eq-p0:se:t0a:k==logn}}
	\]
	(Note that, for $s \in \mbr_+$, the law of $Q_1(sk)$ is independent of $n$ and $k$, so $g$ is a continuous function.)
\end{prop}

\begin{Proof}
Since $\log n/k \eqsim 1/\lambda$, we must choose $s_0$ so that $H(s_0) \eqsim 1/\lambda \in \mbr$.
Since $H$ is strictly increasing and continuous,
we thus deduce that $t_0/k = s_0 \eqsim H^{-1}(1/\lambda) \eqqcolon f(\lambda)$.
So $f$ is a decreasing, continuous bijection from $(0,\infty)$ to itself.
This proves the first part of \cref{eq-p0:se:t0a:k==logn}.


We wish to find times $s_\alpha$ defined so that, recalling \cref{eq-p0:se:var:k=logn},
\[
	H(s_\alpha) = \rbb{ \log n + \alpha \sqrt{v k} } / k
\Qwhere
	v \cq \Varb{Q_1(t_0)} \eqsim \Varb{Q_1(f(\lambda)k)},
\]
which is a constant whose value differs between the undirected and directed cases.

We now turn to finding $s_\alpha$.
Fix $\alpha \in \mbr$.
Note that $G$ is increasing and $\alpha \sqrt{v/k} = \oh1$.
So from the continuity of $H$ we see that, for all $\eps > 0$, we have $(1 - \eps) s_0 \le s_\alpha \le (1 + \eps) s_0$ for $n$ sufficiently large (depending on $\alpha$);
hence $s_\alpha \eqsim s_0$ for all $\alpha \in \mbr$.

Using the same arguments as in the previous derivative proof,
we have
\[
	s_\alpha - s_0
=
	\intt[\alpha]{0} \tfrac{ds_a}{da} \, da
=
	\sqrt{v/k} \intt[\alpha]{0} 1/H'(s_a) \, da.
\]
Using continuity of $H'$ along with $s_\alpha \eqsim s_0 \eqsim f(\lambda)$,
the second part of \cref{eq-p0:se:t0a:k==logn} follows:
\[
	s_\alpha - s_0
&
\eqsim
	\alpha \sqrt{v_*/k} / H'(s_0)
\eqsim
	\alpha \sqrt{v_*/k} / H'(f(\lambda))
\\&
\eqsim
	\alpha s_0 \rbb{ \sqrt{ \Var{Q_1( f(\lambda) \, k } } / \rbb{ f(\lambda) H'(f(\lambda))} } / \sqrt k.
\qedhere
\]
\end{Proof}

\subsubsection{Regime $k \gg \log n$}

Finally we consider the regime $k \gg \log n$, which corresponds to $s_0 = t_0/k \ll 1$ but $t_0 \gg 1$.
We have to handle the directed and undirected cases slightly differently here.
The entropic time $t_0$ and cutoff times $t_\alpha$ will be the same (up to smaller order terms), but the technical details of the proofs will differ ever so slightly.

\begin{prop}
\label{res-p0:se:ent:s<1}
	For $s \ll 1$,
	the entropy $H$ of a rate-1 SRW or DRW on $\mbz$ satisfies
	\[
		H(s) = s \log(1/s) + \Oh{s}.
	\label{eq-p0:se:ent:s<1}
	\nt
	\]
\end{prop}

\begin{Proof}
This follows immediately from \cref{eq-p0:se:CLT:Q1^r} given in the justification of the CLT when $s \ll 1$.
\end{Proof}

\begin{prop}
\label{res-p0:se:der:>}
	For $s \ll 1$,
	the entropy $H$ of a rate-1 SRW or DRW on $\mbz$ satisfies
	\[
		H'(s) = \log(1/s) + \Oh1.
	\label{eq-p0:se:der:<}
	\nt
	\]
	(For SRW, this $\Oh1$ is $\log2 + \Oh{s}$; for DRW, it is $\Oh{s \log(1/s)}$.)
\end{prop}

\begin{Proof}
We proceed as in the previous derivative proof, ie the proof of \cref{res-p0:se:der:s>1}.

Consider first the undirected case.
Using the Kolmogorov backward equations,
we obtain
\[
	H'(s)
=
	\sumt{x \in \mbz} \rbb{ \nu_s(x) - \tfrac12( \nu_s(x+1) + \nu_s(x-1) ) } \log \nu_s(x).
\]
Recall that we have
\[
	\pr{ X_s = 0 } = 1 - s + \Ohb{s^2}
\Qand
	\pr{ X_s = x } = \tfrac12 s + \Ohb{s^2}
\text{ for }
	x \in \bra{\pm1},
\]
and hence $\pr{ X_s = x } = \Oh{s^2}$ for $x \notin \bra{0,\pm1}$.
Also, as previously,
in the above sum we may ignore the $x$ with $x \notin \bra{0,\pm1}$ to give an error $\Oh{s \log(1/s)}$.
(Note that it is not $\Oh{s^2 \log(1/s)}$, since the $x$-th term of the sum contains $\nu_s(x+1)$ and $\nu_s(x-1)$.)
Direct calculation then gives
\[
	H'(s)
=
	\log(1/s) + \log2 + \Oh{s}
=
	\log(1/s) + \Oh1.
\]
This proves the undirected case.

Now consider the directed case.
Here, $X_s \sim \Poisson(s)$, and so $\pr{ X_s = x } = e^{-s} x^s / x!$. Then $W_1(t) \sim X_s$.
Direct differentiation shows that
\[
	\tfrac{d}{ds} \pr{ X_s = x } = \pr{ X_s = x-1 } - \pr{ X_s = x } = e^{-s} s^{x-1} (x-s) / x!
\Qfor
	x \in \mbn,
\]
and $\tfrac{d}{ds} \pr{ X_s = 0 } = -\pr{ X_s = 0 } = -e^{-s}$,
as in the previous derivative proof.
As there, we have
\[
	H'(s)
=
	-\sumt{x \in \mbz_+} \tfrac{d}{ds} \pr{ X_s = x } \rbb{ \log\pr{ X_s = x } + 1 }.
\]
As previously,
we may ignore the terms with $x \notin \bra{0, \pm1}$, giving an error $\Oh{s \log(1/s)}$.
Plugging in the derivative, we obtain
\[
	H'(s)
&
=
	- e^{-s} \log\rbb{e^{-s}} -e^{-s}(1-s) \log\rbb{s e^{-s}} + \Ohb{s \log(1/s)}
\\&
=
	s\rbb{ 1 - s + \Oh{s^2} } - (1-s)\rbb{1 - s + \Oh{s^2}}\rbb{ \log s - s } + \Ohb{ s \log(1/s) }
\\&
=
	\log(1/s) + \Ohb{ s \log(1/s) }
=
	\log(1/s) + \Oh1.
\]
This proves the directed case.
\end{Proof}

We wish to find the times $s_\alpha = t_\alpha/k$ defined so that, recalling \cref{eq-p0:se:var:k>logn},
\[
	H(s_\alpha) = \rbb{ \log n + \alpha \sqrt{v k} }/k
\Qwhere
	v \cq \Varb{ Q_1(t_0) } \eqsim (\log n/k) \log(k/\log n).
\]

\begin{prop}
\label{res-p0:se:ent:k>logn}
	For $k \gg \log n$,
	we have
	\[
		s_0
	=
		t_0/k
	\eqsim
		k^{-1} \log n / \log(k/\log n),
	\tag{\ref{eq-p0:se:t0a:k>>logn}}
	\]
	and,
	for each $\alpha \in \mbr$,
	we have
	$s_\alpha \eqsim s_0$,
	and furthermore
	\[
		\rbr{ s_\alpha - s_0 }/s_0
	=
		\rbr{ t_\alpha - t_0 }/t_0
	\eqsim
		\alpha \sqrt{ \log(k/\log n)/ \log n } = \oh1.
	\tag{\ref{eq-p0:se:t0a:k>>logn}}
	\]
\end{prop}

\begin{Proof}
We consider the directed and undirected cases simultaneously.
By directly manipulating \cref{eq-p0:se:ent:s<1},
we see that
if $H(s_0) = \log n/k$ then
\[
	s_0 \log(1/s_0) \eqsim \log n/k.
\Quad{and hence}
	\log(1/s_0) \eqsim \log(k/\log n),
\]
with the final relation holding since $k \gg \log n$ and so $\log(k/\log n) \gg 1$;
this implies that
\[
	s_0
=
	t_0/k
\eqsim
	k^{-1} \log n / \log(k/\log n).
\]

We now turn to finding $s_\alpha$.
Fix $\alpha \in \mbr$.
From the form \cref{eq-p0:se:ent:s<1} of $H$, observe that
\[
	H\rbb{ s_0 (1 \pm \eps) }
=
	(1 \pm \eps) H(s_0) + \Oh{s_0}
=
	(1 \pm \eps) H(s_0) \cdot \rbb{ 1 + \oh1 },
\]
noting that $s_0 \ll 1$ and so $H(s_0) \eqsim s_0 \log(1/s_0) \gg s_0$.
Note also that
\[
	\sqrt{ v k } \eqsim \sqrt{ \log n \log(k/\log n) } \ll \log n,
\]
since $\log k \ll \log n$.
Hence $H(s_\alpha) = h(t_0) \cdot \rbr{ 1 + \oh1 }$.
Hence, for all $\eps > 0$, we have $(1 - \eps) s_0 \le s_\alpha \le (1 + \eps) s_0$ for $n$ sufficiently large (depending on $\alpha$);
hence $s_\alpha \eqsim s_0$ for all $\alpha \in \mbr$.

As in the previous derivative proofs, we have
\[
	s_\alpha - s_0
=
	\intt[\alpha]{0} \tfrac{ds_a}{da} \, da
=
	\sqrt{v/k} \intt[\alpha]{0} 1/H'(s_a) \, da.
\]
But, by \cref{res-p0:se:der:>}, we may write $H'(s) = \log(1/s) \rbr{ 1 + \oh1 }$ with $\oh1$ term uniform over $t \in [\tfrac12 s_0, 2 s_0]$, which is an interval containing the cutoff window.
Hence, recalling the expressions for $v$ from \cref{eq-p0:se:var:s<1} and $s_0$ from above, the second part of \cref{eq-p0:se:t0a:k>>logn} follows:
\[
	s_\alpha - s_0
\eqsim
	\alpha \sqrt{ s_0 /k }
=
	\alpha s_0 / \sqrt{ s_0 k }
\eqsim
	\alpha s_0 / \sqrt{ \log n / \log(k/\log n) }.
\]
Note that $\log k \ll \log n$, and so $\log n / \log(k/\log n) \gg 1$.
So we do indeed have $\abs{s_\alpha - s_0} = \oh{s_0}$.
\end{Proof}

\begin{rmkt*}
In the directed case, we can actually find an explicit closed-form solution for the entropy:
\[
	H(s)
=
	s \rbb{ \log(1/s) + 1 + e^{-s} \sumt[\infty]{\ell=2} s^{\ell-1} \log(\ell!)/\ell! }.
\]
From this explicit expression, one can derive an approximation to the entropy when $s \gg 1$; see \cite{EBBJ:poisson-entropy}.
An analogous result for $s \ll 1$ is easy to obtain.
For $s \asymp 1$, no simple closed form is known.
\end{rmkt*}

\section{Relative Entropy Estimates, Growth Rates and Concentration}
\label{sec-p0:re}

Let $X^+_\mm \cq (X^+_\mm(s))_{s\ge0}$ be a DRW on $\mbz_\mm$ and $X^-_\mm \cq (X^-_\mm(s))_{s\ge0}$ be a SRW on $\mbz_\mm$.
Throughout this section, we use $+$-superscript to indicate DRW, eg $X^+_\mm(s)$, and $-$-superscript to indicate SRW, eg $X^-_\mm(s)$; when the result holds for both the SRW and DRW, we use either $\pm$-superscript or none at all, eg $X^\pm_\mm(s)$ or $X_\mm(s)$.
Write $\nu_{\mm, s}(\cdot)$ for the law of $X_\mm(s)$, adding $+$/$-$-superscripts as appropriate.


\medskip

The aim of this section is to derive some estimates on relative entropy.
In the first two subsections (\S\ref{sec-p0:re:s>=mix2} and \S\ref{sec-p0:re:s<<mix2}), the results will be for general times $s$.
In the final subsection (\S\ref{sec-p0:re:app}), we are interested in the behaviour of the relative entropy around the so-called \textit{entropic times}; see \S\ref{sec-p0:re:app} for the definition, namely \cref{def-p0:re:app:ent-times}.

\medskip

First, we prove some general estimates on the relative entropy for RW on $\mbz_\mm$.
We then specialise to $s \gtrsim \mm^2$ (in \S\ref{sec-p0:re:s>=mix2}) and then to $s \ll \mm^2$ (in \S\ref{sec-p0:re:s<<mix2}).

\begin{lem}
\label{res-p0:re:prelim-est}
	There exists an absolute constant $c > 0$ so that, for all $\mm \ge 2$ and all $s \ge c$, we have
	\[
		R_\mm(s)
	\ge
		c \log\rbr{\mm/\sqrt s}.
	\label{eq-p0:re:s>=mix2:r>log}
	\nt
	\]
	Moreover, for all $p \ge 2$ and $s \ge 0$, we have
	\[
		\tfrac12 e^{- 2 \gamma_2 s }
	\le
		2 \, d_{\TV, \mm}(s)^2
	\le
		R_\mm(s)
	\le
		d_{\infty, \mm}(2s)
	\le
		\sumt[\mm]{\ell=2} e^{-2 \gamma_\ell s}
	\label{eq-p0:re:s>=mix2:<>-evals}
	\nt
	\]
	where $\gamma_\ell \cq 1 - \cos(2\pi(\ell-1)/\mm)$ for $\ell \in [\mm]$.
	In particular,
	the following hold:
	\begin{alignat*}{2}
		R_\mm(s) &\eqmathsbox{re:prelim}{\ll} 1
	&\Quad{if and only if}
		s &\eqmathsbox{re:prelim}{\gg} \mm^2;
	\\
		R_\mm(s) &\eqmathsbox{re:prelim}{\asymp} 1
	&\Quad{if and only if}
		s &\eqmathsbox{re:prelim}{\asymp} \mm^2;
	\\
		R_\mm(s) &\eqmathsbox{re:prelim}{\gg} 1
	&\Quad{if and only if}
		s &\eqmathsbox{re:prelim}{\ll} \mm^2.
	\end{alignat*}
\end{lem}

\begin{Proof}
The first claim is an immediate consequence of \cite[Proposition~4.1]{HP:hypercont}; it applies for both SRW and DRW.
In particular, in the notation of \cite[Proposition~4.1]{HP:hypercont}, the set $A$ is chosen to be an interval of width $2 \sqrt s$ around the mode of the RW location.


\smallskip

For the lower bound, recall Pinsker's inequality, which says that
\[
	R_\mm(s)
\ge
	2 \, d_{\TV, \mm}(s)^2.
\]
Recall the standard fact that, for any eigenvalue $\psi$ of the transition matrix and $s \ge 0$, we have
\[
	d_{\TV, \mm}(s)
\ge
	\tfrac12 e^{-s \Re(1-\psi)};
\]
see \cite[(12.15)]{LPW:markov-mixing} for the discrete-time analogue.
Write $q \cq e^{-2 \pi i / \mm}$, where $i$ is the imaginary unit (not an index).
The eigenvalues of the transition matrix for DRW, respectively SRW, are given by
\[
	\rbb{ \lambda^+_\ell \cq q^{\ell-1} \mid \ell \in [\mm]},
\Quad{respectively}
	\rbb{ \lambda^-_\ell \cq \Re(q^{\ell-1}) = \Re(\lambda^+_\ell) \mid \ell \in [\mm] }.
\]
Apply this with $\psi \cq \lambda^\pm_2$.
As $\lambda^-_2 = \Re(\lambda^+_2)$, this proves the lower bound for both SRW and DRW.

\smallskip

We turn to the upper bounds.
By Jensen's inequality, for two measures $\mu$ and $\pi$, we have
\[
	\relent \mu\pi
&
=
	\intt{} \mu(x) \log\rbb{ \mu(x)/\pi(x) } \, dx
\le
	\log\rbb{ \intt{} \mu(x)^2 / \pi(x) \, dx}
\\&\qquad
=
	\log\rbb{ 1 + \intt{} \pi(x) \abs{ \mu(x)/\pi(x) - 1 }^2 \, dx }
=
	\log\rbb{ 1 + \norm{ \mu - \pi }_{L^2(\pi)}^2 }.
\]
Applying this and using the inequality $\log(1 + x) \le x$ for $x > -1$,
we obtain
\[
	R_\mm(s)
\le
	\log\rbb{ 1 + d_{2, \mm}(s)^2 }
\le
	d_{2, \mm}(s)^2.
\]
For reversible chains, it is well-known that $d_{2, \mm}(s)^2 = d_{\infty, \mm}(2s)$.
For the DRW, we have $d_{2, \mm}(s)^2 \le d_{\infty, \mm}(2s)$.
Indeed, by symmetry, the $L_2$ mixing profile for the DRW and its time reversal are identical.
The claim then follows from $L_2$--$L_\infty$ mixing time relations in \cite[Appendix]{MT:mixing-times}.

For the SRW, by transitivity, and since we are working in continuous time,
\[
	d^-_{\infty, \mm}(2s)
=
	\mm \, \nu^-_{\mm, 2s}(0) - 1
=
	\mathrm{trace}(P^-_{\mm, 2s}) - 1
=
	\sumt[\mm]{\ell=2} e^{-2 \gamma_\ell s},
\]
where $P^-_{\mm, \cdot}$ is the transition kernel for rate-1 SRW on $\mbz_\mm$.
(See \cite[Lemma~3.20, (3.60)]{AF:book} for justification of the first equality.)
This establishes the upper bound for SRW.

For the DRW, we use the spectral decomposition.
Let $\rbr{ f^+_\ell \mid \ell \in [\mm] }$ be the orthonormal eigenbasis corresponding to $\rbr{ \lambda^+_\ell \mid \ell \in [\mm] }$.
We have $f^+_\ell(x) \cq \exp{ - 2 \pi i (\ell - 1) x / \mm}$ for $x \in \mbz_\mm$.
By the spectral decomposition, for all $s \ge 0$ and all $x,y \in \mbz_\mm$, we have
\[
	P^+_{\mm, 2s}(x,y) - \tfrac1\mm
=
	\tfrac1\mm
	\sumt[\mm]{\ell=2}
	f^+_\ell(x) \overline{f^+_\ell(y)}
	\expb{ - 2s \rbr{ 1 - \lambda^+_\ell } }
\le
	\tfrac1\mm
	\sumt[\mm]{\ell=2}
	\expb{ - 2s \rbb{ 1 - \Re(\lambda^+_\ell) } },
\]
where $P^+_{\mm, \cdot}$ is the transition kernel for rate-1 DRW on $\mbz_\mm$ and we have used the fact that $\abs{f^+_\ell(z)} = 1$ for all $z \in \mbz_\mm$.
As $\Re(\lambda^+_\ell) = \lambda^-_\ell$ for all $\ell \in [\mm]$, this establishes the upper bound for DRW.
\end{Proof}

\subsection{Estimates for $s \gtrsim \mm^2$}
\label{sec-p0:re:s>=mix2}

This subsection is devoted to analysing the regime $s \gtrsim \mm^2$.
(Recall that centred RW is diffusive, and so $\mm^2$ is the order of the mixing and maximal hitting time of the RW.)

\cref{res-p0:re:prelim-est} has the following simple, but extremely useful, corollary.

\begin{cor}
\label{res-p0:re:s>=mix2:asymp}
	For all $\mm \ge 2$, if $s \gtrsim \mm^2$, then 
	\[
		d_{\TV, \mm}(s)^2
	\asymp
		R_\mm(s)
	\asymp
		d_{\infty, \mm}(2s)
	\asymp
		d_{\infty, \mm}(s)^2
	\asymp
		e^{-2 \gamma_2 s}.
	\]
\end{cor}

\begin{Proof}
Note that $\gamma_2 = \gamma_m$.
Hence from \cref{eq-p0:re:s>=mix2:<>-evals} we deduce that
\[
	\tfrac12 e^{-2 \gamma_2 s}
\le
	R_\mm(s)
\le
	e^{-2 \gamma_2 s} \rbb{ 2 + \sumt[\mm-1]{\ell=3} e^{-2( \gamma_\ell - \gamma_2 ) s} }.
\]
Since $s \gtrsim \mm^2$ and $\gamma_\ell - \gamma_2 \gtrsim \min\bra{\ell, \mm-\ell}^2/\mm^2$, the sum above is $\Oh1$.
\end{Proof}

\begin{lem}
\label{res-p0:re:s>=mix2:<nu<}
	For all $c > 0$,
	there exists a constant $\sigma \in (0,\infty)$ so that,
	for all $s \ge c \mm^2$,
	we have
	\[
		1 / \rbb{ 1 + \sigma \sqrt{ R_\mm(s) } }
	\le
		\mm \mint{x \in \mbz_\mm} \nu_{\mm, s}(x)
	\le
		\mm \maxt{x \in \mbz_\mm} \nu_{\mm, s}(x)
	\le
		1 + \sigma \sqrt{ R_\mm(s) }.
	\]
\end{lem}

\begin{Proof}
By \cref{res-p0:re:prelim-est,res-p0:re:s>=mix2:asymp},
there exists a constant $\sigma_+ \in (0,\infty)$ so that
\[
	d_{\infty, \mm}(s)
\le
	\sigma_+ \sqrt{ R_\mm(s) },
\Quad{and hence}
	p \maxt{x} \nu_{\mm, s}(x) \le 1 + \sigma_+ \sqrt{ R_\mm(s) }.
\]
If
\(
	R_\mm(s) \le (2 \sigma_+)^{-2},
\)
then
\(
	1 - \sigma_+ \sqrt{ R_\mm(s) }
\ge
	1 / \rbr{ 1 + 2 \sigma_+ \sqrt{ R_\mm(s) } };
\)
the claim follows with $\sigma \cq 2 \sigma_+$.

It remains to prove the lower bound under the assumption that
\(
	R_\mm(s) \ge (2 \sigma_+)^{-2}.
\)
It then suffices to show that $\mint{x} \nu_{\mm, s}(x) \gtrsim 1/\mm$.
This follows from a relatively simple application of the local CLT (see \cref{res-p0:LCLT}), for either SRW or DRW, noting that $s \gtrsim L^2$.
\end{Proof}


\begin{cor}
\label{res-p0:re:s>=mix2:var}
	For all $c > 0$,
	there exists a constant $\sigma > 0$ so that,
	for all $s \ge c \mm^2$,
	we have
	\[
		\Varb{ Q_{\mm, 1}(sk) } \le \sigma^2 R_\mm(s).
	\]
\end{cor}

\begin{Proof}
Abbreviate $\rho_x \cq \nu_{\mm, s}(x) \cdot \mm$ for $x \in \mbz_\mm$.
Since $Q_{\mm, 1}(sk)$ takes the value
$- \log \nu_{\mm, s}(x)$
with probability $\nu_{\mm, s}(x)$ for each $x \in \mbz_\mm$,
if we define the random variable $Y$ to take the value $\log \rho_x$ with probability $\nu_{\mm, s}(x)$ for each $x \in \mbz_\mm$, then
\[
	\Varb{Q_{\mm, 1}(s k)}
=
	\Var{Y}
\le
	\ex{Y^2}.
\]
Applying \cref{res-p0:re:s>=mix2:<nu<},
we deduce the corollary:
\[
	\ex{Y^2}
=
	\sumt{x} \nu_{\mm, s}(x) \rbr{ \log\rho_x }^2
\le
	\maxt{x} \rbr{ \log\rho_x }^2
\le
	\log\rbb{ 1 + \sigma R_L(s)^{1/2} }^2
\le
	\sigma^2 R_L(s).
\qedhere
\]
\end{Proof}

\subsection{Estimates for $s \ll \mm^2$}
\label{sec-p0:re:s<<mix2}

This subsection is devoted to analysing the regime $s \ll \mm^2$;
however, we only consider $s \ge \varsigma$, for some absolute constant $\varsigma$.
Many of the constants below will depend on the choice of $\varsigma$; however, since $\varsigma$ should be thought of as fixed throughout this whole section, we do not restate this dependence.

\begin{prop}
\label{res-p0:re:s<<mix2:ent-approx}
	Uniformly in all $\mm \in \mbn$,
	we have
	\[
		\maxt{ s \in [r, c \mm^2] }
		\absb{ H_\mm(s) - \tfrac12 \log(2 \pi e s) }
	=
		o_{\toinf r}(1) + o_{\tozero c}(1).
	\]
	Equivalently,
	uniformly in all $\mm \in \mbn$,
	we have
	\[
		\maxt{ s \in [r, c \mm^2] }
		\absb{ R_\mm(s) - \tfrac12 \log(\mm^2/s) - \tfrac12 \log(2 \pi e) }
	=
		o_{\toinf r}(1) + o_{\tozero c}(1).
	\]
\end{prop}

\begin{Proof}
We can uniquely write
\(
	X_\infty(s)
=
	\widetilde X_\mm(s) + \mm M_\mm(s) + m_s
\)
with the following definitions:
\begin{itemize}[noitemsep, topsep = \smallskipamount, label = \bcdot]
	\item 
	$\widetilde X_\mm(\cdot)$ is the RW on $[-\tfrac12 \mm, \tfrac12 \mm]$ centred to have mean 0;
	
	\item 
	$M_\mm(s)$ indicates in which interval of width $L$ the RW on $\mbz$ lives;
	
	\item 
	$m_s$ is the mode of $X_\infty(s)$.
\end{itemize}
As $X_\infty(s)$ determines $\rbr{ X_\mm(s), \: M_\mm(s) }$ and vice versa,
by standard properties of entropy,
we have
\[
	\Ent\rbb{ X_\mm(s) }
\le
	\Ent\rbb{ X_\mm(s), \: M_\mm(s) }
=
	\Ent\rbb{ X_\infty(s) }
\le
	\Ent\rbb{ X_\mm(s) } + \Ent\rbb{ M_\mm(s) }.
\]
The upper bound on $H_\mm(s)$ now follows immediately from \cref{res-p0:se:ent:s>1}.

We now turn to the lower bound.
Using large deviations estimates for the SRW and the Poisson distribution from \cref{res-p0:rp:lde:poi,res-p0:rp:lde:srw},
it is routine to show that
\[
	- \log \Ent\rbb{ M_\mm(s) }
\asymp
	\mm^2 / s
\ge
	1/c.
\qedhere
\]
\end{Proof}

The above proof actually quantifies the errors, in the way described below.

\begin{cor}
\label{res-p0:re:s<<mix2:ent-exp}
	There exists a constant $c$ so that,
	for	all $\mm \in \mbn$
	and
		all $\varsigma \le s \le c \mm^2$,
	we have
	\[
		0 \le H_\infty(s) - H_\mm(s) \le e^{-c \mm^2/s}/c.
	\]
\end{cor}

\begin{lem}
\label{res-p0:re:s<<mix2:var}
	There exist positive constants $c$ and $C$ so that,
	for all $\mm \in \mbn$,
	if $\varsigma \le s \le c \mm^2$, then
	\[
		\Varb{ Q_{\mm, 1}(sk) } \le C.
	\]
\end{lem}

\begin{Proof}
We may assume that $\mm$ is larger than any constant which we desire, otherwise all random variables are order 1 and so the statement holds easily.
%
For $\delta \in (0, \tfrac12)$, consider the set
\[
	A_\delta \cq \brb{x \in \mbz_\mm \midb \nu_{\mm, s}(x) \ge \delta / \sqrt s };
\Quad{write}
	B_\delta \cq \mbz_\mm \setminus A_\delta.
\]
By \cref{res-p0:re:s<<mix2:ent-approx} and the local CLT (\cref{res-p0:LCLT}),
since $\varsigma \le s \le c \mm^2$,
we have
\[
	\ex{ Q_{\mm, 1}(sk) }
=
	H_\mm(s)
=
	\tfrac12 \log s + \Oh1
\Qand
	\maxt{x} \nu_{\mm, s}(x)
=
	\nu_{\mm, s}(0)
\gtrsim
	1/\sqrt s.
\]
From this
and the definition of $A_\delta$,
we deduce the following relations:
\[
	\alpha
&
\cq
	\sumt{x \in A_{\delta}}
	\nu_{\mm, s}(x) \rbb{ \log(1/\nu_{\mm, s}(x)) - \ex{Q_{\mm,1}(sk)} }^2
\lesssim
	1/\delta;
\\
	\beta
&
\cq
	\sumt{x \in B_\delta}
	\nu_{\mm, s}(x) \rbb{ \log(1/\nu_{\mm, s}(x)) - \ex{Q_{\mm,1}(sk)} }^2
\lesssim
	1
+	\sumt{x \in B_{\delta}}
	\nu_{\mm, s}(x) \log \rbb{ \sqrt s \nu_{\mm, s}(x)}^2.
\]

To analyse the sum over $x \in B_\delta$, note that $\delta \le \tfrac 12$.
Under this assumption,
\[
	\nu_{\mm, s}(x) \le s^{-1/2} e^{-\sqrt u}
\Quad{iff}
	\log\rbb{ \sqrt s \nu_{\mm, s}(x) } \le - \sqrt u
\Quad{iff}
	\log\rbb{ \sqrt s \nu_{\mm, s}(x) }^2 \ge u.
\]
From this, using a simple change of variables,
taking $\delta \cq \exp{-\sqrt{10}} \in (0, \tfrac12)$,
we find that
\[
&	\sumt{x \in B_{\delta}}
	\nu_{\mm, s}(x) \log \rbb{ \sqrt s \nu_{\mm, s}(x)}^2
=
	\intt[\infty]{0}
	\pr{ \log\rbb{ \sqrt s \nu_{\mm, s}(X_s) }^2 > u \mid X_s \in B_\delta }
	\pr{ X_s \in B_\delta} du
\\&\qquad
=
	\intt[\infty]{0}
	\pr{ \nu_{\mm, s}(X_s) \le s^{-1/2} \min\bra{e^{-\sqrt u}, \delta} } du
\\&\qquad
\le
	\intt[\infty]{10}
	\nu_{\infty, s}\rbb{ x \in \mbz \mid \nu_{\infty, s}(x) \le s^{-1/2} e^{-\sqrt u} } du
+	10.
\]
It is easy to verify that the last integral is bounded from above, uniformly in $s \ge \varsigma$.

The result now follows, since $\Var{Q_{\mm,1}(sk)} = \alpha + \beta$, and $\delta = \exp{-\sqrt{10}}$.
\end{Proof}

\subsection{Variations Around the Entropic Time: General Abelian Groups}
\label{sec-p0:re:app}

For rate-1 RW on $\mbz_\mm^k$, the entropy function is denoted $h_\mm(\cdot)$.
For rate-1 RW on $\mbz_\mm$, the Shannon, respectively relative, entropy function is denoted $H_\mm(\cdot)$, respectively $R_\mm(\cdot)$; recall that $h_\mm(\cdot) = \log \mm - R_\mm(\cdot)$.
The Shannon entropy functions are strictly increasing bijections with
\[
	h_\mm : [0, \infty) \to \bigl[ 0, \log(\mm^k) \bigr) = [0, k \log \mm)
\Quad{and}
	H_\mm : [0, \infty) \to [0, \log \mm).
\]

\subsubsection{Entropic Time Definitions and Preliminaries}
\label{sec-p0:re:app:ent-def}

We are primarily interested in a target entropy of $N \cq \log\abs{G/\mm G}$, where $G$ is an arbitrary Abelian group $G$.
For certain $\mm$, the time $H_\mm^{-1}\rbr{ \log \abs{G/\mm G} / k }$ may well be $\oh1$, but since $k \lesssim \log n$ the maximum over $\mm$ is at least order 1.
In the definition below, we take a maximum with $\varsigma$.

\begin{defn}
\label{def-p0:re:app:ent-times}
	For $\mm, N \in \mbn$, the \textit{entropic time} is defined by
	\[
		s_0(\mm, N)
	\cq
		H_\mm^{-1} \rbb{ (\log N)/k }
	\Qand
		t_0(\mm, N)
	\cq
		s_0(\mm, N) k
	=
		h_\mm^{-1} \rbr{ \log N }.
	\]
	For the special case $N \cq \abs{G/\mm G}$, write
	\(
		t_\mm \cq s_\mm k
	\Quad{and}
		t_* \cq s_* k
	\)
	where
	\[
		s_\mm
	\cq
		s_0\rbr{ \mm, \abs{G/\mm G} } \vee \varsigma
	\Qand
		s_*
	\cq
		\maxt{\mm \in \mbn} s_\mm.
	\]
	For an Abelian group $G$, write $d(G)$ for the minimal size of a generating subset of $G$.
	Abbreviate
	\[
		\zeta_\mm
	\cq
		\tfrac1k \rbb{ k - d(G) } \log \mm.
	\]
\end{defn}

\begin{lem}
\label{res-p0:re:app:G/LG}
	For all $\mm \in \mbn$, we have
	\(
		\abs{\mm G} \ge \mm^{-d(G)} \abs G,
	\)
	and in particular
	\(
		\abs{G/\mm G} \le \mm^{d(G)}.
	\)
\end{lem}

\begin{Proof}
	Decompose $G$ as $\oplus_1^d \: \mbz_{m_j}$.
	Then $\mm G$ can be decomposed as $\oplus_1^d \: \mbz_{m_j/\gcd(\mm,m_j)}$.
	Thus
	\(
		\abs{\mm G}
	=
		\prodt[d]{1} m_j/\gcd(\mm,m_j)
	\ge
		\prodt[d]{1} m_j/\mm
	=
		\abs G / \mm^d.
	\)
	The second part follows from Lagrange's theorem.
\end{Proof}

\begin{cor}
\label{res-p0:re:app:re>zeta}
	For all $\mm \ge 2$,
	we have
	\[
		R_\mm\rbb{ s_0(\mm, \abs{G/\mm G}) }
	=
		\log \mm - \rbr{ \log \abs{G/\mm G} }/k
	\ge
		\tfrac1k \rbb{ k - d(G) } \log \mm
	=
		\zeta_\mm.
	\]
\end{cor}

We first determine the asymptotic behaviour of $s_*$.
Afterwards, we determine the rate of growth of the entropy around the entropy time.
For both investigations, the following is useful.

Use the usual functional inner product:
\(
	\ipr{f,g}_\pi
\cq
	\sumt{z} f(z) g(z).
\)

\begin{defn}
	For
		all transition matrices $P$
	and
		all functions $f$ and $g$,
	define the \textit{Dirichlet form}
	\[
		\mce_P(f,g)
	\cq
		\ipb{ f, (I - P)g }_\pi
	=
		\sumt{x,y}
		f(x) \rbb{ g(x) - g(y) } P(x,y) \pi(x).
	\]
\end{defn}

For a transition matrix $P$, write $P^*$ for its time reversal; then $P^\times \cq \tfrac12 (P + P^*)$ is its additive symmetrisation.
Observe that, for all functions $f$ and $g$, we have the following:
\[
	\mce_P(f,f)
=
	\tfrac12
	\sumt{x,y}
	\rbb{ f(x) - f(y) }^2 P(x,y) \pi(x)
\Quad{and}
	\mce_P(f,f) = \mce_{P^*}(f,f) = \mce_{P^\times}(f,f);
\]
if $P$ is reversible, then also $\mce(f,g) = \mce(g,f)$.
We now define logarithmic-Sobolev constants.

We write
\(
	\Ent_\pi(g) \cq E_\pi\rbr{ g \log(g / E_\pi(g)) } = \relent{g/E_\pi(g)}{\pi} E_\pi(g)
\)
for a function~$g \ge 0$.

\begin{defn}
\label{def-p0:re:app:logsob}
	Define the \textit{usual}, respectively \textit{modified}, \textit{log-Sobolev constants} by
	\[
		c_{\LS,P}
	\cq
		\inf_{f : f \ne 0}
		\frac{\mce_P(f,f)}{\Ent_\pi(f^2)}
	\Quad{and}
		c_{\MLS,P}
	\cq
		\inf_{f : f > 0}
		\frac{\mce_P(f, \log f)}{\Ent_\pi(f)}.
	\]
	Observe that $c_{\LS,P} = c_{\LS,P^*} = c_{\LS,P^\times}$, ie this is the same for the reversal and the symmetrisation.
\end{defn}

\begin{lem}
\label{res-p0:re:app:logsob:comp:gen}
	For every irreducible transition matrix $P$,
	we have
	\(
		2 c_{\LS,P}
	\le
		c_{\MLS,P}.
	\)
\end{lem}

\begin{Proof}
Using the inequality $\log c \ge 1 - 1/c$ for $c > 0$,
it is straightforward to show that
\(
	\mce(f, \log f)
\ge
	2 \, \mce(\sqrt f, \sqrt f)
\)
for $f > 0$; see, eg, \cite[Lemma~2.8]{MT:mixing-times}.
From this and the definitions, the claim follows.
\end{Proof}

Simple direct calculations establish the following lemma; see, eg, \cite[Lemma~2.4]{MT:mixing-times}.

\begin{lem}[{\cite[Lemma~2.4]{MT:mixing-times}}]
\label{res-p0:re:app:logsob:deriv}
	Let $\Omega$ be a state space and let $s \ge 0$.
	Let $P$ be an irreducible transition matrix with invariant distribution $\pi$ and write $P_s \cq e^{s(P-I)}$ for its heat kernel.
	Let $\mu$ be a distribution on $\Omega$ and write $\mu_s \cq \mu P_s$, ie the law of the chain started from $\mu$ and run for time $s$.
	For $x \in \Omega$, write $h_s(x) \cq \mu_s(x)/\pi(x)$, ie the density with respect to $\pi$.
	Then
	\[
		\tfrac{d}{ds}
		\relent{\mu_s}{\pi}
	=
		\tfrac{d}{ds}
		\Ent_\pi(h_s)
	=
		- \mce( h_s, \log h_s )
	\le
		- c_{\MLS,P} \Ent_\pi(h_s)
	=
		- c_{\MLS,P} \relent{\mu_s}{\pi}.
	\]
\end{lem}

\begin{cor}
\label{res-p0:re:app:logsob:final}
	In the set-up of \cref{res-p0:re:app:logsob:deriv},
	we have
	\[
		\relent{\mu_s}{\pi}
	\le
		\relent{\mu}{\pi} e^{-c_{\MLS,P} s}
	\le
		\relent{\mu}{\pi} e^{-2 c_{\LS,P} s}.
	\]
\end{cor}

\begin{Proof}
This follows immediately from \cref{res-p0:re:app:logsob:comp:gen,res-p0:re:app:logsob:deriv} and Gronwall's lemma.
\end{Proof}

It remains to estimate the log-Sobolev constant for the random walks on $\mbz_\mm$---recall that this is the same for the SRW and DRW, as the SRW is the additive symmetrisation of the DRW.

\begin{lem}
\label{res-p0:re:app:logsob:comp:cycle}
	For all $\mm \in \mbn$,
	the log-Sobolev constants of the RW on $\mbz_\mm$ satisfy
	\[
		c_{\MLS,\mm} \ge 2 c_{\LS,\mm} \gtrsim 1/\mm^2.
	\]
\end{lem}

\begin{Proof}
In \cite[Corollary~3.11]{DSc:log-sobolev},
it is shown that the $L_2$ mixing time is bounded below by $1/(2 c_{\LS,P})$.
For the SRW on $\mbz_\LL$, the $L_2$ mixing time is well-known to be order $\mm^2$.
The claim follows.
\end{Proof}

\subsubsection{Asymptotic Evaluation of Entropic Time}
\label{sec-p0:re:app:ent-eval}

The precise definitions of $s_0(\mm, N)$ and $t_0(\mm, N)$
Here we asymptotically evaluate the entropic time.
We first determine its order when $1 \ll k \lesssim \log \abs G$;
second we evaluate it asymptotically when $1 \ll k \ll \log \abs G$ and $k - d(G) \asymp k$;
finally we evaluate it asymptotically when $k \gg \log \abs G$.

\begin{subtheorem}{thm}
	\label{res-p0:re:app:s*:order}

\begin{prop}
\label{res-p0:re:app:s*:order:asymp}
	Let $d,n \in \mbn$.
	Suppose that $1 \ll k \lesssim \log n$ and $k - d \asymp k$.
	Then,
	with implicit constant uniform over all Abelian groups $G$ with $\abs G = n$ and $d(G) = d$,
	we have
	\[
		\maxt{\mm \in \mbn}
		t_0\rbr{ \mm, \abs{G/\mm G} }
	\asymp
		k \abs G^{2/k}.
	\]
\end{prop}

\begin{prop}
\label{res-p0:re:app:s*:order:>>}
	Let $d,n \in \mbn$.
	Suppose that $1 \ll k \lesssim \log n$ and $k > d$.
	Then,
	with implicit constant uniform over all Abelian groups $G$ with $\abs G = n$ and $d(G) = d$,
	we have
	\[
		k \abs G^{2/k}
	\lesssim
		\maxt{\mm \in \mbn}
		t_0\rbr{ \mm, \abs{G/\mm G} }
	\lesssim
		k \abs G^{2/k} \log k.
	\]
\end{prop}

\end{subtheorem}

When $d \ll \log n$ and $k - d(G) \asymp k$,
the entropic time $t_*(k, G)$ is asymptotically equivalent to $t_0(\infty, \abs G)$, ie the time at which the entropy of the RW on $\mbz^k$ reaches $\log \abs G$.
This entropic time $t_0(\infty, \abs G)$ is evaluated, in different regimes, in \cref{res-p0:se:t0a}; there it is denoted $t_0(k, \abs G)$.

\begin{prop}
\label{res-p0:re:app:s*:eqsim:d<<log}
	Let $d,n \in \mbn$.
	Suppose that $d \ll \log n$ and $k - d \asymp k$.
	Then,
	with implicit constant uniform over all Abelian groups $G$ with $\abs G = n$ and $d(G) = d$,
	we have
	\[
		\maxt{\mm \in \mbn}
		t_0\rbr{ \mm, \abs{G/\mm G} }
	\eqsim
		t_0\rbr{ \infty, \abs G }.
	\]
\end{prop}

\begin{prop}
\label{res-p0:re:app:s*:eqsim:k>>log}
	Suppose that $k \gg \log \abs G$.
	Write $\rho \cq \log k / \log \log \abs G$.
	Then
	\[
		\maxt{\mm \in \mbn}
		t_0\rbr{ \mm, \abs{G/\mm G} }
	\eqsim
		t_0(\infty, \abs G)
	\eqsim
		\tfrac{\rho}{\rho-1} \log_k \abs G
	=
		\log(\abs G) / \log(k/\log \abs G).
	\]
\end{prop}

\begin{Proof}[Proof of \cref{res-p0:re:app:s*:order:asymp}]
For the lower bound,
take $\mm \cq \infty$ and use \cref{res-p0:se:t0a}
\[
	\widetilde s_*
\cq
	\maxt{\mm \in \mbn}
	s_0\rbr{ \mm, \abs{G/\mm G} }
\ge
	s_0(\infty, \abs G)
\asymp
	\abs G^{2/k}.
\]
In particular, this says that the maximising $\mm$ satisfies
\(
	s_0(\mm, \abs{G/\mm G}) \gtrsim \abs G^{2/k} \gtrsim 1.
\)

We now turn to the upper bound.
Assume that $k - d \ge \eps k$ with $\eps \in (0,1)$.
By \cref{res-p0:re:app:G/LG},
the entropy of a single-coordinate at $t \cq t_0(\mm, \abs{G/\mm G})$,
which we denote $x$,
satisfies
\[
	x
=
	\tfrac1k \log \abs{G/\mm G}
\le
	\rbb{ \tfrac1k \log n } \wedge \rbb{ \tfrac1k \log(\mm^d) }
\le
	\rbb{ \tfrac1k \log n } \wedge \rbb{ (1 - \eps) \log \mm }.
\]
The relative entropy, which we denote $\xi$, thus satisfies
\(
	\xi = \log \mm - x \ge \eps \log \mm \gtrsim 1.
\)
By \cref{res-p0:re:prelim-est,res-p0:re:s<<mix2:ent-approx},
there exists a constant $C$ so that if $\xi \ge C$ then
\(
	t/k \le C e^{2x} \le C n^{2/k}.
\)
On the other hand, if $\xi \le C$, then necessarily $\xi \asymp 1$ and $\eps \log \mm \le C$, ie $\mm \le e^{C/\eps}$.
\cref{res-p0:re:prelim-est} then implies that
\(
	t/k \asymp \mm^2 \le e^{2C/\eps} \le e^{2C/\eps} n^{2/k},
\)
with implicit constant uniform over the compact interval $[\eps \log 2, C]$ in which $\xi$ lies.
Hence, in either case, we have $t \lesssim C' k n^{2/k}$.
\end{Proof}

\begin{Proof}[Proof of \cref{res-p0:re:app:s*:order:>>}]
The lower bound of $s_0(\infty, \abs G) \asymp \abs G^{2/k}$ holds using \cref{res-p0:se:t0a}.

The relative entropy of a single coordinate at time $t \cq t_0(\mm, \abs{G/\mm G})$,
denoted $\xi$,
satisfies
\[
	\xi
=
	\log \mm - \tfrac1k \log \abs{G/\mm G}
\ge
	\log \mm - \rbb{ \tfrac1k \log n } \wedge \rbb{ \tfrac1k \log(\mm^d) }
=
	\rbb{ \tfrac1k \rbr{k-d} \log \mm } \vee \rbb{ \log \mm - \tfrac1k \log n}.
\]

Consider first the case that $\xi = \tfrac1k(k-d) \log \mm \eqqcolon \zeta$, ie $\mm \le n^{1/d}$.
This is precisely the entropic time studied in \S\ref{sec-p0:re:p} below, specifically \cref{def-p0:re:p:entropic-time,res-p0:re:p:t0} with $\pp \cq \mm$ and $\alpha \cq 0$ (in the notation there).
These references show that
\(
	t/k
\lesssim
	(\mm^d)^{2/k} \lvert \log(\zeta \wedge \tfrac12) \rvert.
\)
Since here we are considering $L \le n^{1/d}$, the proof is completed in this case.

Now suppose that $\xi = \log \mm - \tfrac1k \log n$, ie $L \ge n^{1/d}$.
Then
\(
	\xi
\ge
	k^{-2} \rbr{ k - d } \log n
\ge
	1/k^2.
\)
By \cref{res-p0:re:prelim-est,res-p0:re:s<<mix2:ent-approx},
there exists a constant $C$ so that if $\log \mm - \tfrac1k \log n \ge C$ then
\(
	t/k \le C n^{2/k}
\)
as in the proof of \cref{res-p0:re:app:s*:order:asymp}.
On the other hand, if $\log \mm - \tfrac1k \log n \le C$, then necessarily $\xi \lesssim 1$ and $n^{1/d} \le \mm \lesssim n^{1/k}$.
\cref{res-p0:re:s>=mix2:asymp} implies that $R_\mm\rbr{ C' \mm^2 \log k } \le 1/k^2 \le \xi$ for a sufficiently large constant $C'$.
Since $C' \mm^2 \log k \lesssim n^{2/k} \log k$ the proof is completed in this case.
\end{Proof}

\begin{Proof}[Proof of \cref{res-p0:re:app:s*:eqsim:d<<log}]
For the lower bound,
take $\mm \cq \infty$ and use \cref{res-p0:se:t0a}:
\[
	\widetilde s_*
\cq
	\maxt{\mm \in \mbn}
	s_0\rbr{ \mm, \abs{G/\mm G} }
\ge
	s_0(\infty, n).
\]
When $k \ll \log n$, we have
\(
	s_0(\infty, n)
\eqsim
	\tfrac1{2 \pi e} n^{2/k}
\gg
	1
\)
while
\(
	s_0(\infty, n)
\asymp
	1
\)
when $k \asymp \log n$.

We now turn to the upper bound.
Choose $\mm$ to be an optimiser,
ie with
\(
	\widetilde s_*
=
	s_0\rbr{ \mm, \abs{G/\mm G} }.
\)
By \cref{res-p0:re:app:G/LG},
we have
\(
	\widetilde s_*
\le
	s_0\rbr{ \mm, \mm^d }
\eqqcolon
	\widetilde s_\mm.
\)
We consider first the case that $1 \ll k \ll \log n$ with $k - d \asymp k$; so $\widetilde s_* \gg 1$.
If $k \gtrsim \log(\mm^d)$, then $\widetilde s_\mm \asymp 1$.
But $\widetilde s_\mm \gg 1$, so we must have
\(
	1
\ll
	k
\ll
	\log(\mm^d).
\)
But $k \ge d$, so we must have $\mm \gg 1$.
Hence $\zeta \cq \tfrac1k (k - d) \log \mm \gg 1$.
By \cref{res-p0:re:prelim-est,res-p0:re:app:re>zeta}, we thus have
\(
	R_\mm(\widetilde s_\mm)
\ge
	\zeta
\gg
	1
\)
and
\(
	\widetilde s_\mm
\ll
	\mm^2.
\)
By \cref{res-p0:re:s<<mix2:ent-approx},
\[
	H_\mm(\widetilde s_\mm)
=
	\tfrac12 \log\rbr{ 2 \pi e \widetilde s_\mm } + \oh1.
\]
The target entropy is
\(
	\tfrac1k \log \abs{G/\mm G}
\le
	\tfrac1k \log n.
\)
Thus
\(
	\widetilde s_\mm
\eqsim
	\tfrac1{2 \pi e} n^{2/k},
\)
completing the upper bound.

Consider now $k \asymp \log n$ with $d \ll \log n$.
Choose $m \gg 1$ such that $d \log m \ll \log n \asymp k$.
Consider first $\mm \le m$.
Since $\log(m^d) / k = d \log m / k \ll 1$, we have
\(
	\widetilde s_\mm
\le
	s_0(\mm, m^d)
\ll
	1.
\)
But $\widetilde s_* \asymp 1$, so we must have $\mm \ge m \gg 1$.
Since $\widetilde s_\mm = \widetilde s_* \asymp 1$, for $\mm \ge m$, we have $\widetilde s_\mm \ll \mm^2$.
Thus we may apply \cref{res-p0:re:s<<mix2:ent-exp} to deduce that
\(
	0
\le
	H_\infty(\widetilde s_*) - H_\mm(\widetilde s_*)
=
	\oh1.
\)
Further, $H_\infty(\widetilde s_*) \asymp 1$.
We thus deduce that the entropic times for the RW on $\mbz^k$ and $\mbz_\mm^k$ are asymptotically equivalent.
\end{Proof}

\begin{Proof}[Proof of \cref{res-p0:re:app:s*:eqsim:k>>log}]
We start with the lower bound.
Clearly $t_0(\mm, \abs{G/\mm G}) \ge t_0(\infty, \abs G)$.
By \cref{eq-p0:se:t0a:k==logn} and some simple algebraic manipulations,
we have
\(
	t_0(\infty, \abs G)
\eqsim
	\tfrac{\rho}{\rho-1} \log_k \abs G.
\)

We turn to the upper bound.
Clearly $t_0(\mm, \abs{G/\mm G}) \le t_0(2, \abs G)$ for all $\mm \in \mbn$.
In the regime $k \gg \log \abs G$,
in \S\ref{sec-p0:se},
to prove \cref{eq-p0:se:t0a:k==logn},
we approximated the rate-1 RW run for time $s \ll 1$ on $\mbz$ by one on $\mbz_2$.
Thus the arguments for \cref{eq-p0:se:t0a:k==logn} imply the upper bound here.
\end{Proof}
\color{black}

\subsubsection{Rate of Change of Entropy Around the Entropic Time}

We now move onto determining the rate of growth of the entropy.
The following lemma is valid for any $s \ge \varsigma$, but we are particularly interested in applying it at an entropic time $s_\mm$.
(This is one place in which we need the bound $s \ge \varsigma$, and so need to deal with $s_\mm$, rather than $s_0\rbr{ \mm, \abs{G/\mm G} }$.)

\begin{lem}
\label{res-p0:re:app:ent-growth-rate}
	There exists a continuous function $\widetilde c : (0,1) \to (0,1)$ so that,
	for all $\mm \ge 2$, all $\xi \in (-1,1) \setminus \bra{0}$ and all $s \ge \varsigma$,
	we have
	\[
		\absb{ H_\mm\rbb{ s(1 + \xi) } - H_\mm(s) }
	\ge
		2 \widetilde c_{\abs \xi} \rbb{ R_\mm(s) \wedge 1 }.
	\]
\end{lem}

\begin{Proof}
If $s \asymp 1$, then the claim is immediate, noting that $R_\mm(s) \asymp 1$.
Now assume that $s \gg 1$.

Consider first the case where $s/\mm^2$ is small; in particular, $\mm$ is large.
By \cref{res-p0:re:s<<mix2:ent-approx},
there exists constants $\mm_0 \in \mbn$ and $\alpha, c \in (0,\infty)$ so that,
for all $\mm \ge \mm_0$ and all $s \in [\varsigma, 2 \alpha \mm^2]$,
the difference in entropy is
\(
	\tfrac12 \log(1 + \xi) + \oh1.
\)
The claim thus follows in this case.

Now suppose that $s \ge \alpha \mm^2$.
By \cref{res-p0:re:app:logsob:final},
for all $\mm \ge 2$ and all $s \ge \alpha \mm^2$,
we have
\[
	H_\mm\rbb{ s (1 + \eps) } - H_\mm\rbr{ s }
=
	R_\mm\rbr{ s } - R_\mm\rbb{ s (1 + \eps) }
\ge
	(1 - e^{-2 c_{\LS,\mm} s \eps}) R_\mm\rbr{ s }
\ge
	\delta_\eps R_\mm\rbr{ s },
\]
where
\(
	\delta_\eps
\cq
	\liminf_\mm \bra{ 1 - e^{-2 \alpha \eps c_{\LS,\mm} \mm^2 } }
\in
	(0,1)
\)
by \cref{res-p0:re:app:logsob:comp:cycle}.
This completes the proof.
	%
\end{Proof}

Abbreviate
\(
	\rho_\mm
\cq
	R_\mm(s_\mm).
\)
By \cref{res-p0:re:app:re>zeta}, if $s_\mm = s_0(\mm, \abs{G/\mm G})$, then
\(
	\rho_\mm \ge \zeta_\mm.
\)

\begin{prop}
\label{res-p0:re:app:conc}
	Suppose that $k - d(G) \gg 1$
	There exists a continuous function $c : (0,1) \to (0,1)$ so that,
	for all $\mm \ge 2$ with $\mm \wr \abs G$ and all $\eps \in (0,1)$,
	the following hold:
	\[
		\pr{
			Q_\mm\rbb{ t_* (1 + \eps) }
		\le
			\log \abs{G/\mm G} + c_\eps \rbr{ \zeta_\mm \wedge 1 } k
		}
	&\le
		\expb{ - c_\eps \rbr{ \zeta_\mm \wedge 1 } k };
	\\
		\pr{
			Q_\mm\rbb{ t (1 - \eps) }
		\ge
			\log \abs{G/\mm G} - c_\eps \rbr{ \zeta_\mm \wedge 1 } k
		}
	&=
		\oh1
	\Quad{for all}
		t \le t_0(\mm, \abs{G/\mm G}).
	\]
\end{prop}

The outline is of the proof is relatively straightforward.
Replace $s_\mm$ with $s_0(\mm, \abs{G/\mm G})$.
Consider $k - d \asymp k$, so that $\zeta_\mm \asymp \zeta_\mm \wedge 1 \asymp 1$.
Both parts use the entropy growth rate lemma, \cref{res-p0:re:app:ent-growth-rate}.
The non-quantitative part is then an application of Chebyshev's inequality, once one has shown that the variance $\Var{Q_{\mm,1}(sk)}$ is uniformly bounded over $s \ge \varsigma k$.
The quantitative part requires a (one-sided) large deviations estimate given below in \cref{res-p0:re:app:lde:thm}.
We are not exactly sure who proved this originally; the earliest reference that we can find is in a survey by \textcite[Theorem~2.7]{Mcd:concentration};
we use the version given in the very nice survey paper by \textcite[Theorem~3.4]{CL:concentration-survey}.

\begin{thm}
\label{res-p0:re:app:lde:thm}
	Let $(\xi_i)_{i=1}^k$ be a sequence of iid, mean-0 random variables with $\xi_1 \ge - M$ (deterministically), for some $M$.
	Set $\sigma^2 \cq \Var{\xi_1} = \ex{\xi_1^2}$.
	For all $x > 0$, we have
	\[
		\prb{ \sumt[k]{1} \xi_i \le -x }
	\le
		B(x,M,k\sigma^2)
	\Qwhere
		B(x,M,v^2)
	\cq
		\expb{ - \tfrac12 x^2 /(v^2 + xM/3) }.
	\]
\end{thm}

Recall the definition of the random variable $Q$ and the entropies $h$ and $H$.
Define
\[
	\xi_i
\cq
	Q_{\mm,i}(t) - H_\mm(t/k);
\Quad{in particular recall that}
	\ex{Q_{\mm,i}(t)} = H_\mm(t/k).
\]

To apply the large deviations estimate to $\sumt[k]{1} \xi_i$, we wish to find an $M \in \mbr$ so that $\xi_1 \ge -M$ deterministically.
We also need to bound the variance.
The following auxiliary lemmas do these.

\begin{subtheorem}{thm}
	\label{res-p0:re:app:lde:app}

\begin{lem}
\label{res-p0:re:app:lde:app:M}
	There exists an absolute constant $M$ so that,
	for all $\mm \ge 2$ and $s \ge \varsigma$,
	we have
	\[
		Q_{\mm,1}(sk) - \ex{ Q_{\mm,1}(sk) }
	\ge
		- M \rbb{ \sqrt{ R_\mm(s) } \wedge 1 }.
	\]
\end{lem}

\begin{lem}
\label{res-p0:re:app:lde:app:var}
	There exist an absolute constant $\sigma^2$ so that,
	for all $\mm \ge 2$ and $s \ge \varsigma$,
	we have
	\[
		\Varb{ Q_{\mm,1}(sk) }
	\le
		\sigma^2 \rbb{ R_\mm(s) \wedge 1 }.
	\]
\end{lem}

\end{subtheorem}

%

We combine these lemmas to get our own large deviations estimate on $Q_\mm(\cdot)$.

\begin{Proof}[Proof of \cref{res-p0:re:app:conc}]
Let $\eps \in (0,1)$.
We are interested at looking at time $t_*$, which satisfies
\[
	\maxt{\mm \in \mbn} t_\mm
=
	t_*
=
	\maxt{\mm \in \mbn} t_0\rbr{ \mm, \abs{G/\mm G} }.
\]
Let $\mm \in \mbn$ and set $t_0 \cq t_0\rbr{ \mm, \abs{G/\mm G} }$ and $s_0 \cq t_0/k$.
Abbreviate $r_\mm \cq R_\mm(s_0)$ and $\hat r_\mm \cq r_\mm \wedge 1$.
By \cref{res-p0:re:app:re>zeta}, we have
\(
	r_\mm
\ge
	\zeta_\mm.
\)
By \cref{res-p0:re:app:ent-growth-rate},
there exists a constant $\widetilde c_\eps > 0$ so that
\[
	h_\mm\rbb{ t_0( 1 + \eps) } - h_\mm\rbr{ t_0 }
\ge
	2 \widetilde c_\eps \hat r_\mm k
\ge
	2 \widetilde c_\eps \hat \zeta_\mm k,
\]
where $\zeta_\mm \cq \zeta_\mm \wedge 1$.
Recall that $\ex{Q_\mm(t')} = h(t')$ for all $t' \ge 0$.
For each $i \in [k]$,
set
\[
	\xi_i
\cq
	Q_{\mm,i}\rbb{ t_0(1 + \eps) } - \ex{ Q_{\mm,i}\rbb{ t_0(1 + \eps) } }.
\]
Altogether, these relations imply that
\[
	\brb{ Q_\mm\rbb{ t_0(1 + \eps) } \le \log \abs{G/\mm G} + \widetilde c_\eps \hat \zeta_\mm k }
\subseteq
	\brb{ \sumt[k]{1} \xi_i \le - \widetilde c_\eps \hat \zeta_\mm k }.
\]
Applying the large deviations estimate \cref{res-p0:re:app:lde:thm}, with parameters controlled by Lemma~\ref{res-p0:re:app:lde:app},
setting
\(
	c_\eps
\cq
	\tfrac12 \widetilde c_\eps^2 / \rbr{ \sigma^2 + \tfrac13 \widetilde c_\eps M },
\)
a little algebra shows that
\[
	\pr{ Q_\mm\rbb{ t_0(1 + \eps) } \le \log \abs{G/\mm G} + \widetilde c_\eps \hat \zeta_\mm k }
\le
	\expb{ - c_\eps \hat \zeta_\mm k }.
\]

We want to apply this for $\max_{\mm \in \mbn} t_0\rbr{ \mm, \abs{G/\mm G} }$ instead of each individual $t_0\rbr{ \mm, \abs{G/\mm G} }$.
This follows from the above analysis, due to the fact that $t' \mapsto Q_\mm(t')$ is stochastically increasing.

\smallskip

For the lower bound, observe that the growth rate lemma \cref{res-p0:re:app:ent-growth-rate} has the same form for time $1 - \eps$ as for $1 + \eps$.
As for the upper bound, it suffices to prove the result for $t \cq t_0(\mm, \abs{G/\mm G})$.
We apply Chebyshev's inequality, with the variance controlled by \cref{res-p0:re:app:lde:app:var}.
The standard deviation is order $\rbr{ \hat \zeta_\mm k }^{1/2}$ and the displacement order $\hat \zeta_\mm k$.
So to deduce the result from Chebyshev's inequality, we need $\hat \zeta_\mm k \gg 1$.
This is immediate: $\hat \zeta_\mm k = (k - d) \log \mm \gg 1$ as $k - d(G) \gg 1$.
\end{Proof}

It remains to prove Lemma~\ref{res-p0:re:app:lde:app}, which has two parts.

\begin{Proof}[Proof of \cref{res-p0:re:app:lde:app:M}]
	%
Recall that $R_\mm(s) \lesssim 1$ if $s \gtrsim \mm^2$ and $R_\mm(s) \gtrsim 1$ if $s \lesssim \mm^2$.
We have
\[
	\ex{ Q_{\mm,1}(sk) } - Q_{\mm,1}(sk)
\le
	\log \mm + \log\rbb{ \maxt{x \in \mbz_\mm} \nu_{\mm,s}(x) }
=
	\log\rbb{ 1 + d_{\infty, \mm}(s) }
\le
	d_{\infty, \mm}(s).
\]
By \cref{res-p0:re:prelim-est},
writing $\gamma_\ell \cq 1 - \cos\rbr{ 2 \pi (\ell-1)/ \mm }$ for $\ell \in [\mm]$,
for both SRW and DRW,
we have
\[
	d_{\infty, \mm}(s)
\le
	\sumt[\mm]{\ell=2} e^{-\gamma_\ell s}
=
	2 e^{ - \gamma_2 s } \rbb{ 1 + \tfrac12 \sumt[\mm-1]{\ell=3} e^{- (\gamma_\ell - \gamma_2) s} }.
\]
Since $\gamma_\ell - \gamma_2 \asymp \min\bra{\ell, \mm - \ell}^2 / \mm^2$,
there exists a constant $\beta$ so that,
if $s \ge \beta \mm^2$, then
\[
	\ex{ Q_{\mm,1}(sk) } - Q_{\mm,1}(sk)
\le
	d_{\infty, \mm}(s)
\le
	5 \sqrt{ R_\mm(s) }.
\]

On the other hand, if $\varsigma \le s \le \beta \mm^2$, then we upper bound
\[
	\ex{ Q_{\mm,1}(sk) }
=
	H_\mm(s)
\le
	H_\infty(s)
=
	\tfrac12 \log\rbr{ 2 \pi e s } + \Ohb{ s^{-1/4} },
\]
with the last relation following from \cref{res-p0:se:ent:s>1}.
By the local CLT (see, eg, \cite[Theorem~2.5.6]{LL:RW} or \cref{res-p0:LCLT}),
the mode has probability order $1/\sqrt s$ in this regime.
Hence
\[
	\ex{ Q_{\mm,1}(sk) } - Q_{\mm,1}(sk)
\le
	\Oh1.
\qedhere
\]
	%
\end{Proof}

\begin{Proof}[Proof of \cref{res-p0:re:app:lde:app:var}]
This is an immediate consequence of \cref{res-p0:re:s>=mix2:var,res-p0:re:s<<mix2:var}.
\end{Proof}

\begin{Proof}[Proof of \cref{res-p0:re:app:lde:app:var}]
Recall that $R_\mm(s) \lesssim 1$ if $s \gtrsim \mm^2$ and $R_\mm(s) \gtrsim 1$ if $s \lesssim \mm^2$.

By \cref{res-p0:re:s<<mix2:var}, there exist positive constants $\beta$ and $C$ so that,
if $\varsigma \le s \le \beta \mm^2$, then
\[
	\Varb{ Q_{\mm,1}(sk) } \le C.
\]
On the other hand,
if $s \ge \beta \mm^2$, then 
by \cref{res-p0:re:s>=mix2:var} we have
\[
	\Varb{ Q_{\mm,1}(sk) } \lesssim R_\mm(s).
\qedhere
\]
	%
\end{Proof}

%

\subsection{Variations Around the Entropic Time: The Special Case of $\mathbb Z_p^d$}
\label{sec-p0:re:p}

In this section we specialise to the group $\mbz_\pp^d$; these entropic results do not require $\pp$ to be prime.
(Note that $d(\mbz_\pp^d) = d$.)
Here we not only establish cutoff, but also get a bound on the order of the cutoff window when $(k - d)\pp \gg 1$.
This section has two main propositions.

Use the notation from the previous sections, but drop the $\mm$-subscripts: we only consider RWs on $\mbz_\pp$ or $\mbz_\pp^k$. This is because the maximiser $\mm$ is clearly given by $\mm = \pp$.


\smallskip

We first define precisely the entropic times under consideration here.

\begin{defn}
\label{def-p0:re:p:entropic-time}
Recall that
\(
	\zeta
=
	\tfrac1k (k - d) \log \pp.
\)
For $\alpha \in \mbr$,
define
\[
	t_\alpha
\cq
	h^{-1}\rbb{ d \log p + 2 \alpha \sqrt{ k(\zeta \wedge 1) } }.
\]
Equivalently,
\(
	t_\alpha \cq s_\alpha k
\)
where
\[
	\zeta_\alpha
\cq
	\zeta - 2 \alpha \sqrt{(\zeta \wedge 1)/k}
=
	\zeta \rbb{ 1 - 2 \alpha / \sqrt{\zeta k (\zeta \vee 1)} }
\Qand
	s_\alpha
\cq
	R^{-1}(\zeta_\alpha).
\]
We call $t_0$ the \textit{entropic time} and $\bra{t_\alpha}_{\alpha \in \mbr}$ \textit{cutoff times}.
Note that $\zeta_0 = \zeta$.
\end{defn}

The next proposition estimates these entropic times.

\begin{subtheorem}{thm}
	\label{res-p0:re:p:t0a}

\begin{prop}[Entropic Times]
\label{res-p0:re:p:t0}
	Suppose that $1 \ll k \lesssim d \log \pp$.
	The following hold:
	\begin{alignat*}{2}
		\text{if}
	\quad
		\zeta &\eqmathsbox{re:p:t0}{\ll} 1,
	&\Quad{then}
		t_0/k = s_0 &\eqsim \tfrac12 \log(1/\zeta) / \rbb{ 1 - \cos(2\pi/\pp) };
	\\
		\text{if}
	\quad
		\zeta &\eqmathsbox{re:p:t0}{\gtrsim} 1,
	&\Quad{then}
		t_0/k = s_0 &\asymp \pp^2 e^{-2\zeta} = (\pp^d)^{2/k};
	\intertext{further, if in fact $1 \ll k \ll d \log \pp$, then}
		\text{if}
	\quad
		\zeta &\eqmathsbox{re:p:t0}{\gg} 1,
	&\Quad{then}
		t_0/k = s_0 &\eqsim \pp^2 e^{-2\zeta} / (2 \pi e) = (\pp^d)^{2/k} / (2 \pi e).
	\end{alignat*}
	Note that
	\(
		1 - \cos(2\pi/\pp)
	\eqsim_{\toinf \pp}
		2 \pi^2 / \pp^2
	=
		2 \pi^2 \pp^{-2d/k} e^{2\zeta}.
	\)
\end{prop}

\begin{prop}[Cutoff Times]
\label{res-p0:re:p:ta}
	Suppose that $1 \ll k \lesssim d \log \pp$ and $(k - d) \pp \gg 1$, ie $\zeta \gg 1/k$.
	Then,
	for all $\alpha \in \mbr$,
	we have
	\(
		t_\alpha \eqsim t_0
	\)
	and furthermore the following hold:
	\begin{alignat*}{3}
		\text{if}
	\quad
		\zeta &\eqmathsbox{re:p:ta}{\lesssim} 1,
	\Quad{then}&
		(t_\alpha - t_0) / t_0 &\lesssim 1 / \rbb{ \sqrt{\zeta k} \log( (1/\zeta) \vee e ) };
	\\
		\text{if}
	\quad
		\zeta &\eqmathsbox{re:p:ta}{\gg} 1,
	\Quad{then}&
		(t_\alpha - t_0) / t_0 &\lesssim 1 / \sqrt k
	\quad
		\text{for the SRW}.
	\end{alignat*}
\end{prop}

\end{subtheorem}

\begin{rmkt}
\label{rmk-p0:re:p:ta}
	We strongly believe that the last result also holds for the DRW; see \cref{rmk-p0:re:p:deriv:rmk} for justification of this belief.
	In short, $\zeta \gg 1$ implies that $s_0 \ll \pp^2$, and so the RW on $\mbz_\pp$ should look almost the same as the RW on $\mbz$, once recentred to have mean 0.
	In particular, the growth of the entropy (as a function of time) should be similar.
\end{rmkt}

The next result is a concentration result.
For $\alpha \in \mbr$, define
\[
	Q_\alpha^+ \cq \bra{ Q(t_{ \alpha}) \ge \log n + \alpha \sqrt{ k (\zeta \wedge 1) } }
\Qand
	Q_\alpha^- \cq \bra{ Q(t_{-\alpha}) \le \log n - \alpha \sqrt{ k (\zeta \wedge 1) } };
\]

\begin{prop}[Concentration]
\label{res-p0:re:p:conc}
	For all $\alpha \in (0,\infty)$ with $\abs{ \zeta_\alpha - \zeta_0 } \le \tfrac12 \zeta_0$,
	we have
	\[
		\pr{ (Q_\alpha^\pm)^c } \lesssim \alpha^{-2}.
	\]
\end{prop}

This proposition is an easy consequence of relative entropy results proved earlier in \S\ref{sec-p0:re}.

\begin{Proof}[Proof of \cref{res-p0:re:p:conc}]
Using the definition of $\zeta_\alpha$ and $\zeta$,
we have
\[
	h(t_\alpha)
=
	\log n + 2 \alpha \sqrt k (\zeta \wedge 1);
\]
recall that $n = d \log p$.
Note that
\(
	\Var{Q} = k \Var{Q_1}.
\)
By Chebyshev's inequality,
we have
\[
	\pr{ \absb{ Q(t_\alpha) - \rbb{ \log n + 2 \alpha \sqrt{ k (\zeta \wedge 1) } } } \ge \abs\alpha \sqrt{ k (\zeta \wedge 1) } }
\le
	\alpha^{-2} \, \Varb{ Q_1(t_\alpha) } / (\zeta \wedge 1).
\]
Consider $\zeta \lesssim 1$.
\cref{res-p0:re:prelim-est} implies that $s_\alpha \gtrsim \pp^2$.
Then,
by \cref{res-p0:re:s>=mix2:var},
we have
\(
	\Var{Q(t_\alpha)} \lesssim \zeta_\alpha \asymp \zeta.
\)
From this and the definition of $Q_\alpha^\pm$, the proposition follows.
When ${\zeta \gg 1}$, \cref{res-p0:re:prelim-est} gives $s_\alpha \ll \pp^2$.
The argument proceeds as before, replacing \cref{res-p0:re:s>=mix2:var} with \cref{res-p0:re:s<<mix2:var}.
	%
\end{Proof}

We separate the proof of Proposition~\ref{res-p0:re:p:t0a} into multiple parts.

\begin{Proof}[Proof of \ref{res-p0:re:p:t0}]
Apply \cref{res-p0:re:prelim-est,res-p0:re:s>=mix2:asymp} together:
	for $\zeta \ll 1$, they imply that
	\(
		s_0
	\eqsim
		\tfrac12 \gamma_2^{-1} \log(1/\zeta);
	\)
	for $\zeta \asymp 1$, they imply that
	\(
		s_0
	\asymp
		\gamma_2^{-1}.
	\)
Also, note that
\(
	\gamma_2
=
	1 - \cos(2\pi/\pp)
\eqsim_{\toinf \pp}
	2 \pi^2 / \pp^2.
\)
For $\zeta \gg 1$, \cref{res-p0:re:s<<mix2:ent-approx} implies that
\(
	s_0
\asymp
	n^{2/k}
=
	p^{2d/k}.
\)
Further, if $k \ll d \log \pp = \log n$, then $s_0 \gg 1$, and so in fact \cref{res-p0:re:s<<mix2:ent-approx} implies that
\(
	s_0
\eqsim
	n^{2/k} / (2 \pi e).
\)
	%
\end{Proof}

We first show cutoff, namely $t_\alpha \eqsim t_0$ for all $\alpha \in \mbr$.
This just uses the entropy growth rate.

\begin{Proof}[Proof of \ref{res-p0:re:p:ta}: Cutoff]
Since $\zeta k \gg 1$,
we have
\[
	R(s_\alpha)
=
	\zeta_\alpha
=
	\zeta \rbb{ 1 + \oh1 }
=
	R(s_0) \rbb{ 1 + \oh1 }.
\]
But by \cref{res-p0:re:app:ent-growth-rate},
replacing $s_0$ by $s_0(1 + \xi)$ changes the entropy by at least order $\zeta \wedge 1$.
Case analysis gives
\(
	s_0(1 - \eps) \le s_\alpha \le s_0(1 + \eps)
\)
asymptotically for all $\alpha \in \mbr$ for all $\eps \in (0,1)$.
\end{Proof}

We next bound the window.
Note that $(k - d)\pp \gg 1$ implies $\zeta k (\zeta \vee 1) \gg 1$.

\begin{Proof}[Proof of \ref{res-p0:re:p:ta}: Window when $\zeta \lesssim 1$]
Recall \cref{res-p0:re:app:logsob:final,res-p0:re:app:logsob:comp:cycle}:
\[
	R(u+v) \le R(v) e^{-2 c_{\LS,\pp} u}
\Quad{for all}
	u,v \ge 0
\Qand
	c_{\LS,\pp} \gtrsim 1/\pp^2.
\]

Consider first $\alpha > 0$.
Applying this
with $v \cq s_0$ and $u \cq s_\alpha - s_0$
gives
\[
	\zeta_\alpha
=
	R(s_\alpha)
\le
	e^{-c_{\LS,\pp} (s_\alpha - s_0)} R(s_0)
=
	e^{-c_{\LS,\pp} (s_\alpha - s_0)} \zeta_0.
\]
We hence deduce that
\(
	s_\alpha - s_0
\le
	c_{\LS,\pp}^{-1} \log(\zeta_0/\zeta_\alpha),
\)
ie
\[
	s_\alpha - s_0
\le
	- c_{\LS,\pp}^{-1} \log\rbb{ 1 - \alpha / \sqrt{ \zeta k (\zeta \vee 1) } }
\asymp
	\alpha \pp^2 / \sqrt{ \zeta k (\zeta \vee 1) }
\asymp
	\alpha \pp^2 / \sqrt{ \zeta k }.
\]
For $\alpha < 0$,
use $v \cq s_\alpha$ and $u \cq s_0 - s_\alpha$
to deduce the analogous result.
Hence
\[
	\abs{ s_\alpha - s_0 }
\asymp
	\abs \alpha \pp^2 / \sqrt{ \zeta k }
\asymp
	\abs \alpha s_0 \big/ \rbb{ \sqrt{ \zeta k } \log( (1/\zeta) \vee e ) }.
\qedhere
\]
\end{Proof}

To analyse the window,
we need to use the derivative of the entropy.

\begin{lem}
\label{res-p0:re:p:deriv:srw}
	There exist positive constants $c$ and $c'$ so that,
	for all $\mm \in \mbn$,
	if $\varsigma \le s \le c \mm^2$, then
	\[
		\tfrac{d}{ds} H^-_\mm(s)
	=
		- \tfrac{d}{ds} R^-_\mm(s)
	\ge
		c' / s.
	\]
\end{lem}

\begin{rmkt}
\label{rmk-p0:re:p:deriv:rmk}
	This intuition for this claim is simple.
	For $s \ll \mm^2$, the RW on $\mbz_\mm$ and $\mbz$ look almost the same.
	This is quantified by \cref{res-p0:re:s<<mix2:ent-exp}.
	We showed the the RW on $\mbz$ that the derivative satisfies
	\(
		H'_\infty(s) \eqsim 1/(2s).
	\)
	We thus expect the RW on $\mbz_\mm$ to exhibit the same property when $s \ll \mm^2$.
	However, due to a technical hurdle, we have only been able to show this for the SRW.
	
	The claim is somewhat analogous to the usual log-Sobolev inequality (see \cref{res-p0:re:app:logsob:final} above):
	\[
		R_\mm(u+v)
	\le
		R_\mm(v)
		e^{-c_{\LS, \mm} v}
	\Quad{for all}
		u,v \ge 0.
	\qedhere
	\]
\end{rmkt}

\begin{Proof}[Proof of \ref{res-p0:re:p:ta}: Window when $\zeta \gg 1$ for SRW]
In \cref{res-p0:se:ent:k<logn} we performed an analogous calculation.
Exactly as there,
using \cref{res-p0:re:p:deriv:srw},
we deduce that
\(
	\abs{s_\alpha - s_0} / s_0
\lesssim
	\abs \alpha / \sqrt k.
\)
\end{Proof}

Finally, we prove \cref{res-p0:re:p:deriv:srw}.

\begin{Proof}[Proof of \cref{res-p0:re:p:deriv:srw}]
We may assume that $s$, and hence $\mm$, is larger than any constant which we desire, otherwise all terms on the left-hand side are order 1 and so the statement holds easily.

There is a constant $C$ sufficiently large so that,
writing $K \cq C \sqrt s$,
we~have
\[
	\nu^-_{\mm, s}([-K, K])
\cq
	\sumt{x \in [-K, K]} \nu^-_{\mm, s}(x)
\ge
	\tfrac89.
\]
Moreover, provided $c$ is sufficiently small, we can choose $C$ so that $K \le \mm/10$. Define $\xi^-$ by
\[
	\xi^-_x
\cq
	\nu^-_{\mm, s}(x) \one{x \in [-K, K]} / \nu^-_{\mm, s}([-K, K]).
\]

Let $\mcu$ denote the uniform distribution on $\mbz_{2K+1}$.
For $g : \mbz_{2K+1} \to (0,\infty)$, define
\[
	\Ent_\mcu(g)
\cq
	E_\mcu\rbb{ g \log\rbb{ g / E_\mcu\rbr{g} } }.
\]
By definition of the \textit{modified log-Sobolev} constant (of RW on $\mbz_{2K+1}$), denoted $c_{\MLS, 2K+1}$,
we have
\[
	\sumt{x \in \mbz_{2K+1}}
	\tfrac12 \mcu(x) \rbb{ g(x) - g(x+1) }
	\log\rbb{ g(x) / g(x+1) }
\ge
	c_{\MLS, 2K+1} \Ent_\mcu(g);
\label{eq-p0:re:s<<mix2:deriv:mls}
\nt
\]
see \cref{def-p0:re:app:logsob} below.
Further, it is well-known that,
\(
	c_{\MLS, 2K+1}
\gtrsim
	c_{\LS, 2K+1}
\gtrsim
	1/K^2;
\)
see \cref{res-p0:re:app:logsob:comp:cycle} below.
This holds for both the SRW and the DRW.

Using the backward Kolmogorov equations,
an elementary calculation for the SRW gives
\[
	- \tfrac{d}{ds} R^-_\mm(s)
&
=
	\sumt{x \in \mbz_\mm}
	\tfrac{d}{ds} \nu^-_{\mm, s}(x) \cdot \log \nu^-_{\mm, s}(x)
\\&
=
	\sumt{x \in \mbz_\mm}
	\tfrac12
	\rbb{ \nu^-_{\mm, s}(x) - \nu^-_{\mm, s}(x+1) }
	\log\rbb{ \nu^-_{\mm, s}(x) / \nu^-_{\mm, s}(x+1) }.
\label{eq-p0:re:s<<mix2:deriv:kolmo}
\nt
\]
This actually holds for the DRW too; we explain this at the end.
Note that $(a - b)\log(a/b) \ge 0$ whenever $a,b > 0$.
Hence all terms above are non-negative,
and so,
recalling the definition of $\xi^-$ above,
we have
\[
	- \tfrac{d}{ds} R^-_\mm(s)
\ge
	\nu^-_{\mm, s}([-K, K]) \sumt{x \in [-K, K]} \tfrac12 \rbb{ \xi^-_x - \xi^-_{x+1} } \log\rbb{ \xi^-_x/\xi^-_{x+1} },
\]
where we identify $\mbz_{2K+1}$ with $[-K, K] \cap \mbz$ and $K+1 \equiv -K$, and used the fact that $\xi^-_K = \xi^-_{-K}$.

Combining this with \cref{eq-p0:re:s<<mix2:deriv:mls},
noting that $\nu_{\mm, s}([-K, K]) \ge \tfrac89$,
we see that 
\[
	- \tfrac{d}{ds} R^-_\mm(s)
&
\ge
	\nu_{\mm, s}([-K, K]) \sumt{x \in \mbz_{2K+1}} \tfrac12 \rbb{ \xi^-_x - \xi^-_{x+1} } \log\rbb{ \tfrac{\xi^-_x}{\xi^-_{x+1}} }
\\&
=
	\nu_{\mm, s}([-K, K]) \sumt{x \in \mbz_{2K+1}} \mcu(x) \tfrac12 \rbb{ (2K+1) \xi^-_x - (2K+1) \xi^-_{x+1} } \log\rbb{ \tfrac{(2K+1) \xi^-_x}{(2K+1) \xi^-_{x+1}} }
\\&
\gtrsim
	K^{-2} \Ent_\mcu\rbb{(2K+1) \xi^-}
=
	K^{-2} \relent{\xi^-}{\mcu},
\]
where the final expression is the the relative entropy of $\xi^-$ with respect to $\mcu$ on $\mbz_{2K+1}$.

We argue that $\relent{\xi^-}{\mcu} \gtrsim 1$; the lemma then follows, since $K \asymp \sqrt s$.
By standard exit time estimates for SRW on a cycle,
\(
	\nu^-_{\mm,s}([-K,K])
\asymp
	1
\)
since $K \asymp \sqrt s$.
Since $K \le \mm/10$, the support of $\xi$, namely $[-K,K] \subseteq \mbz_\mm$, contains fewer than half the vertices of $\mbz_\mm$.
Hence by Pinsker's inequality and then the triangle inequality,
we have
\(
	\tfrac12 \relent{\xi^-}{\mcu}^2
\ge
	\tv{ \xi^- - \mcu }
\ge
	\tfrac12.
\)
\end{Proof}

\subsection{Comparison of Entropic Times Between Related Groups}
\label{sec-p0:re:abe-com}

Here we compare the entropic times $t_*(k, A)$ and $t_*(k, A \oplus B)$ where $A$ and $B$ are Abelian groups with $\log \abs B \ll \log \abs A$.
We derive conditions under which these are asymptotically equivalent.

\begin{prop}
\label{res-p0:re:abe-com:eqsim}
	Let $A$ and $B$ be finite, Abelian groups and $k$ be such that $1 \ll \log k \ll \log \abs A$.
	\begin{itemize}[noitemsep, topsep = \smallskipamount, label = \bcdot]
		\item 
		If $k \eqmathsbox{ent-com}{\lesssim} \log \abs{A \oplus B}$, then
		suppose that $k \gg d(B) \log \abs B$ and $k - d(A) \gg d(B)$.
		
		\item 
		If $k \eqmathsbox{ent-com}{\gg} \log \abs{A \oplus B}$, then
		suppose only that $\log \abs B \ll \log \abs A$.
	\end{itemize}
	(In either case, $\log \abs B \ll \log \abs A$.)
	Then
	\(
		t_*(k, A \oplus B)
	\eqsim
		t_*(k, A).
	\)
\end{prop}

\begin{Proof}
For any Abelian group $H$, if $k \gg \log \abs H$ then
\(
	t_*(k, \abs H)
\eqsim
	\log(\abs H) / \log(k / \log \abs H);
\)
see \cref{res-p0:re:app:s*:eqsim:k>>log}.
From this the claim when $k \gg \log \abs{A \oplus B}$ follows immediately.

Now consider $1 \ll k \lesssim \log \abs A$.
Recall the definition of the entropic times $t_*$ and $t_\mm$ for $\mm \in \mbn \cup \bra{\infty}$ from \cref{def-p0:re:app:ent-times}.
(Recall that we use the convention $\mbz_\infty = \mbz$ and $\infty H = \bra{\id}$ for any Abelian group $H$.)
In particular, for an Abelian group $H$, for $t_\mm$ the target entropy for the RW on $\mbz_\mm^k$ is $\log \abs{H/\mm H}$.
The lower bound
\(
	t_\mm(k, A \oplus B)
\ge
	t_\mm(k, A)
\)
is immediate for all $\mm$.

Regarding the upper bound,
the difference in target entropies for $t_\mm(A \oplus B)$ vs $t_\mm(A)$ is
\[
	\log \rbb{ \abs{A \oplus B} / \abs{\mm(A \oplus B)} }
-	\log \rbb{ \abs{A} / \abs{\mm A} }
=
	\log \rbb{ \abs{B} / \abs{\mm B} }.
\]

Consider first $\mm = \infty$.
The difference is then simply $\log \abs B$.
By \cref{res-p0:se:entropic-time-calc:+omega}, at $t_\infty(A)$, increasing the time by a factor $1 + \eps$, for arbitrarily small but constant $\eps > 0$, increases the entropy by an additive term of order $k$.
Hence if $k \gg \log \abs B$ then this increase swallows up the difference in the target entropies.
Thus
\(
	t_\infty(A \oplus B)
\eqsim
	t_\infty(A).
\)

Consider now $\mm \in \mbn$.
Write $\mfgcd \cq \gcd(\mm, \abs B)$.
We have $\mm B = \mfgcd B$ and $\abs{B} / \abs{\mfgcd B} \le \mfgcd^{d(B)}$ by \cref{res-p0:re:app:G/LG}.
Hence the above difference is at most
\(
	d(B) \log \mfgcd.
\)
By \cref{res-p0:re:app:ent-growth-rate,res-p0:re:app:re>zeta}, at $t_\mm(A)$, increasing the time by a factor $1 + \eps$, for arbitrarily small but constant $\eps > 0$, increases the entropy by an additive term of order $k \wedge \rbr{ (k - d(A)) \log \mm }$.
We thus desire
\[
	d(B) \log \mfgcd
\ll
	\rbb{ k - d(A) } \log \mm
\Qand
	d(B) \log \mfgcd
\ll
	k.
\]
Since $\mfgcd \le \min\bra{\mm, \abs B}$,
these hold if
	$k - d(A) \gg d(B)$
and
	$k \gg d(B) \log \abs B$.
\end{Proof}

\section{Large Deviation Estimates for Random Walk on $\mathbb Z$}
\label{sec-p0:rp}

The aim of this section is to prove a large deviations result for the RW on $\mbz$.
First we must define the times at which we wish to evaluate the RW.
Roughly, we look at the time at which the entropy of the RW on $\mbz$ is $\log n / k$.
This corresponds to roughly the time in the following definition.

\begin{defn}
\label{def-p0:rp:s0-def}
	Abbreviate $\kappa \cq k / \log n$.
	Let $s_0 \cq s_0(k,n)$ be any time satisfying
	\[
		s_0 \lesssim n^{2/k} \log k
	\text{ when }
		k \lesssim \log n
	\Qand
		s_0 \eqsim 1/(\kappa\log\kappa)
	\text{ when }
		k \gg \log n.
	\]
\end{defn}

When we say ``RW'', we mean either a SRW or a DRW.

\begin{defn}
\label{def-p0:rp:rp-def}
	Let $X = \Xs$ be a rate-1 RW on $\mbz$.
	Define $r(k,n)$ and $p(k,n)$ as follows:
	\[
		r(k,n)
	&\cq
		\min
		\brb{ r \in \mbz \midb \pr{ \absb{ X_{s_0} - \ext{X_{s_0}} } > r } \le 1/k^{3/2} };
	\\
		p(k,n)
	&\cq
		\min
		\brb{ \pr{ X_{s_0} - \ext{X_{s_0}} = j } \midb \abs{j} \le r(k,n) }.
	\]
	Also define
	\(
		r_*(k,n) \cq \tfrac12 n^{1/k} \logk[2]
	\Quad{and}
		p_*(k,n) \cq n^{-1/k} k^{-2}.
	\)
\end{defn}

\begin{prop}
\label{res-p0:rp:res}
	We have
	\(
		r(k,n) \ge r_*(k,n)
	\Quad{and}
		p(k,n) \ge p_*(k,n).
	\)
\end{prop}

\begin{rmkt}
\label{rmk-p0:rp:k-large}
	An easy inspection of the proof reveals that we may replace $r_*$ and $p_*$ by $\tfrac12 \log k$ and $k^{-2}$, respectively, in the regime $k \gtrsim \log n$.
\end{rmkt}

This proposition will follow from standard large deviation theory, but the details are non-trivial.
The exponent 2 in $\logk[2]$ is not optimal, but is chosen for convenience of proof and to enable us to deal with all regimes of $k$ simultaneously.

The following propositions provide asymptotic estimates for tails of the Poisson distribution and for continuous-time SRW on $\mbz$, as well as for the ratio between the `tail' and `point' probabilities.
We note that in the regime $r \in [\sqrt s, s^{2/3}]$ stronger assertions can be made via the local CLT \cref{eq-p0:LCLT}.

Below, for $a,b \in \mbr$, we write $a \vee b \cq \max\bra{a,b}$ and $a \wedge b \cq \min\bra{a,b}$.

\begin{prop}[Poisson Bounds]
\label{res-p0:rp:lde:poi}
	For $s \in (0,\infty)$, let $X_s \sim \Poisson(s)$.
	Then, uniformly in $s \in (0,\infty)$ and in $r$ with $r \ge \sqrt s$ and $s+r \in \mbz$,
	we have the following relations:
	\begin{mixedsubequations}[eq-p0:rp:lde:poi:log,eq-p0:rp:lde:poi:rat]
		^^A\label{eq-p0:rp:lde:poi:log}\label{eq-p0:rp:lde:poi:rat}
	\[
		- \log \pr{ X_s \ge s + r }
	&\asymp
		r \rbb{(r/s) \wedge 1} \log \rbb{ (r/s) \vee e};
	\label{eq-p0:rp:lde:poi:log:>}
	\nt
	\\
		\pr{ X_s \ge s + r } / \pr{ X_s = s + r }
	&\asymp
		(s/r) \vee 1.
	\label{eq-p0:rp:lde:poi:rat:>}
	\nt
	\]
	\end{mixedsubequations}
	Moreover, uniformly in $s \in (0,\infty)$ and in $r \in [\sqrt s, s]$ with $s-r \in \mbz$
	we have the following relations:
	\begin{mixedsubequations}
	\[
		- \log \pr{ X_s \le s - r }
	&\asymp
		r \rbb{(r/s) \wedge 1} \log \rbb{ (r/s) \vee e};
	\label{eq-p0:rp:lde:poi:log:<}
	\nt
	\\
		\pr{ X_s \le s - r } / \pr{ X_s = s - r }
	&\asymp
		(s/r) \vee 1.
	\label{eq-p0:rp:lde:poi:rat:<}
	\nt
	\]
	\end{mixedsubequations}
\end{prop}

\begin{prop}[SRW Bounds]
\label{res-p0:rp:lde:srw}
	Let $X = (X_s)_{s\ge0}$ be a rate-1 SRW on $\mbz$ started at 0.
	Then, uniformly in $s \in (0,\infty)$ and in $r$ with $r \ge \sqrt s$ and $r \in \mbz$,
	we have the following relations:
	\[
		- \log \pr{ X_s \ge r }
	&\asymp
		r \rbb{(r/s) \wedge 1} \log \rbb{ (r/s) \vee e};
	\label{eq-p0:rp:lde:srw:log}
	\nt
	\\
		\pr{ X_s \ge r } / \pr{ X_s = r }
	&\asymp
		(s/r) \vee 1.
	\label{eq-p0:rp:lde:srw:rat}
	\nt
	\]
\end{prop}

From these, we can deduce the proof of \cref{res-p0:rp:res}.

\begin{Proof}[Proof of \cref{res-p0:rp:res}]
Recall that $\kappa = k / \log n$ and that the time $s$ being considered satisfies
\[
	s \lesssim n^{2/k} \log k
\text{ when }
	k \lesssim \log n
\Qand
	s \eqsim 1/(\kappa\log\kappa)
\text{ when }
	k \gg \log n.
\]

\smallskip

Consider the SRW.
Equations (\ref{eq-p0:rp:lde:poi:log}--\ref{eq-p0:rp:lde:srw:rat}) are all ``$f \asymp g$''-type statements; let $c > 0$ be a universal constant such that $c$ is a lower and $C \cq 1/c$ an upper bound for these relations.

For $r$, it is enough to find an $\widetilde r$ so that
\[
	- \log \pr{ X_s \ge \widetilde r } \ge 2 \log k.
\]
For $p$, since we only consider $j$ with $\abs j \le r$, and $r$ is defined as a minimum, we have
\(
	\pr{ X_s \ge \abs j } \ge k^{-3/2}
\)
for all such $j$.
We split into two regimes, namely $s \ge 2C \log k$ and $s < 2C \log k$.

First suppose that $s \ge 2C \log k$.
Set $\widetilde r \cq \sqrt{ 2 C s \log k }$. Then $\widetilde r \le s$, and so, by \cref{eq-p0:rp:lde:srw:log}, we have
\[
	- \log \pr{ X_s \ge \widetilde r }
\ge
	c \widetilde r \rbb{ (\widetilde r/s) \wedge 1 } \log\rbb{ (\widetilde r/s) \vee e }
=
	c \widetilde r^2/s
\ge
	2 \log k.
\]
For $p_*$, since $\widetilde r \le s$, by \cref{eq-p0:rp:lde:srw:rat}, we have
\[
	\pr{ X_s = j }
\gtrsim
	(s/r) \pr{ X_s \ge j }
\gtrsim
	\logk[1/2] n^{-1/k} \cdot k^{-3/2}
\gg
	n^{-1/k} k^{-2}.
\]

Suppose now that $s < 2C \log k$.
Set $\widetilde r \cq 2C \log k$. Then $\widetilde r \ge s$, and so, by \cref{eq-p0:rp:lde:srw:log}, we have
\[
	- \log \pr{ X_s \ge \widetilde r }
\ge
	c \widetilde r \rbb{ (\widetilde r/s) \wedge 1 } \log\rbb{ (\widetilde r/s) \vee e }
\ge
	c \widetilde r
=
	2 \log k.
\]
For $p_*$, since $\widetilde r \ge s$, by \cref{eq-p0:rp:lde:srw:rat}, we have
\[
	\pr{ X_s = j }
\gtrsim
	\pr{ X_s \ge j }
\ge
	k^{-3/2}
\gg
	k^{-2}
\ge
	n^{-1/k} k^{-2}.
\]

Observe that, in either regime, we have $\widetilde r \le r_*$, with $r_*$ defined in \cref{def-p0:rp:rp-def}.
This completes the proof of \cref{res-p0:rp:res} in the undirected case.

\smallskip

The DRW case, using Poisson bounds, is in essence the same, due to the similarity of \cref{res-p0:rp:lde:poi,res-p0:rp:lde:srw}.
It is slightly messier to write down, since one must take care that $s+r \ge 0$.
\end{Proof}

\begin{Proof}[Proof of \cref{res-p0:rp:lde:poi} \textnormal{(Poisson)}]
For $s \le 10$, all that is needed is the observation that
\[
	\pr{ X_s \ge r }
\asymp
	\pr{ X_s = r }
\asymp
	s^r / r!
\asymp
	\rbb{es/r}^r / \sqrt r.
\]

We now consider the case $s \ge 1$.
First we state that, for all $r \ge 0$, we have
\[
	\max\brb{ \pr{ X_s \ge s + r }, \pr{ X_s \le s - r } }
\le
	\expb{ - \tfrac12 r^2 / (s + r/3) };
\label{eq-p0:rp:lde:proof:poi:bern}
\nt
\]
this follows from Bernstein's inequality, by taking an appropriate limit.

A direct calculation involving Stirling's approximation shows, uniformly in $s$ and in $r$ with $r \ge \tfrac12 s$ and $s+r \in \mbz$, respectively $\tfrac12 s \le r \le s$, the following relations:
\[
	\pr{ X_s \ge s + r }
&\asymp
	\pr{ X_s = s + r }
\asymp
	\frac{ e^r \rbb{ s / (s+r) }^{s+r} }{ \sqrt{ 2 \pi(s+r) } },
\\
	\pr{ X_s \le s - r }
&\asymp
	\pr{ X_s = s - r }
\asymp
	\frac{ e^r \rbb{ s / (s-r) }^{s-r} }{ \sqrt{ 2 \pi(s-r) } };
\]
from these, one can verify (\ref{eq-p0:rp:lde:poi:rat:>}, \ref{eq-p0:rp:lde:poi:rat:<}) for such $r$.

We can obtain lower bounds on $\pr{ X_s \ge s + r }$ and $\pr{ X_s \le s - r }$ for $r \le \tfrac12 s$, from which, together with \cref{eq-p0:rp:lde:proof:poi:bern}, we can verify (\ref{eq-p0:rp:lde:poi:rat:>}, \ref{eq-p0:rp:lde:poi:rat:<}) for such $r$:
\[
	\pr{ X_s = s + r } \sqrt{ 2 \pi (s+r) }
\asymp
	e^r \rbB{ \frac{s}{s+r} }^{s+r} 
&\asymp
	\expbb{ - \frac{r^2}{2(s+r)} - \Ohbb{ \frac{r^3}{(s+r)^2} } },
\\
	\pr{ X_s = s - r } \sqrt{ 2 \pi (s-r) }
\asymp
	e^{-r} \rbB{ \frac{s}{s-r} }^{s-r} 
&\asymp
	\expbb{ - \frac{r^2}{2(s-r)} - \Ohbb{ \frac{r^3}{(s-r)^2} } };
\]
these are found using Stirling's approximation, and both hold uniformly for $r \le \tfrac12 s$.

\smallskip

We now prove \cref{eq-p0:rp:lde:poi:log:>}; the proof of \cref{eq-p0:rp:lde:poi:log:<} is similar and is omitted.
We consider $s \ge 10$, having already considered $s \le 10$ initially.
Observe that $r \mapsto \pr{ X_s = s \pm r }$ is decreasing on $r \ge 0$ with $s \pm r \in \mbz$.
Using the formula for $\pr{\Poisson(\lambda) = k}$, we have
\[
	\frac{ \pr{ X_s = s+r } }{ \pr{ X_s = s+r+1 } }
=
	\frac{ s+r+1 }{ s }.
\]

If $r \ge \tfrac14s$, then this ratio is at least $11/9$, when $s \ge 10$, from which one can readily see that \cref{eq-p0:rp:lde:poi:log:>} holds.
Now suppose that $r \in [\sqrt s, \tfrac14 s]$.
To conclude the proof, we show that there exist universal constants $c_1,c_2 \in (0,1)$ so that, for such $r$, we have
\[
	c_1 \, \pr{ X_s = s + r }
\le
	\pr{ X_s = s + r + \ceil{s/(2r)} }
\le
	c_2 \, \pr{ X_s = s + r }.
\label{eq-p0:rp:lde:proof:poi-suff}
\nt
\]
This, together with the decreasing statement above, can easily be seen to imply \cref{eq-p0:rp:lde:poi:log:>}.
We now prove \cref{eq-p0:rp:lde:proof:poi-suff}.
If $\sqrt s \le r \le \tfrac14 s$, then
\[
	\frac{ \pr{ X_s = s+r } }{ \pr{ X_s = r+r+j } }
& =
	\prod_{i=1}^j
	\frac{ s+r+i }{ s }
=
	\prod_{i=1}^j
	\rbb{1 + (r+i)/s }
\\&
\le
	\expb{ \sumt[j]{i=1} (r+i)/s }
=
	\expb{ \tfrac12 j (j+2r+1) / s }.
\]
If in addition $j \le \tfrac12 s/r$, then the last estimate is tight up to a constant factor.
Indeed, in this case we have $\exp{ \tfrac12 j(j+2r+1)/s } \le e^3$.
Conversely, using the fact that $1 + \theta \ge \exp{\theta - 2\theta^2}$ for $\theta \in [0,\tfrac12]$, we find some universal constant $c_0 > 1$ so that $\exp{ \tfrac12 j(j+2r+1)/s } \ge c_0$.
\end{Proof}

\begin{Proof}[Proof of \cref{res-p0:rp:lde:srw} \textnormal{(SRW)}]
Fix an $s \in (0,\infty)$; without loss of generality, assume $r \ge 0$.
Recall that $X$ has the same law as $Y_N \cq \sum_1^N \xi_i$, where $(\xi_i)_{i\in\mbn}$ is an iid sequence of random variables with $\pr{\xi_1 = +1} = \tfrac12 = \pr{\xi_1 = -1}$ and $N \sim \Poisson(s)$, independent of $(\xi_i)_{i\in\mbn}$.
Then $(Y_k \cq \sum_1^k \xi_i)_{k \in \mbz_+}$ is a discrete-time SRW on $\mbz$ started at the origin.

\smallskip

We first prove \cref{eq-p0:rp:lde:srw:log}.
Observe that $\ex{e^{\lambda \xi_1}} = \tfrac12 e^\lambda + \tfrac12 e^{-\lambda} \le e^{\lambda^2/2}$,
and so $\ex{e^{\lambda Y_k}} \le e^{\lambda^2 k/2}$,
and hence $\pr{Y_k \ge r} \le \exp{-r^2/(2k)}$, by taking $\lambda \cq r/k$.
Further, an elementary calculation involving Stirling's approximation shows, uniformly over $r$ with $\sqrt{k \log k} < r \le k$, that
\[
	- \log \pr{ Y_k \ge r }
\le
	- \log \pr{ Y_k \in \bra{r,r+1} }
\asymp
	r^2/k;
\]
for $\sqrt k \le r \le \sqrt{ k \log k }$ one can use the local CLT (see \cref{res-p0:LCLT}) to verify that
\[
	- \log \pr{ Y_k \ge r }
\asymp
	r^2/k.
\]

For $r \le \sqrt2 s$, we average over $N$ and use the above bounds on $Y_k$.
In particular, we have
\[
	\ex{ e^{\lambda X_s} }
\le
	\sumt[\infty]{r=0}
	\pr{ N = k } \, e^{\lambda^2 k/2}
=
	\ex{ e^{\lambda^2 N/2} }
=
	\expb{ s (e^{\lambda^2/2} - 1) }
\le
	\expb{ s \rbr{ \lambda^2/2 + (\lambda^2/2)^2 } },
\]
with the final inequality holding when $\lambda^2 \le 2$, applying the inequality $e^\theta - 1 \le \theta + \theta^2$ valid for $\theta \in [-1,1]$.
We now set $\lambda \cq r/s$ and use Chernoff to deduce that
\[
	\pr{ X_s \ge r }
\le
	\expb{ - \tfrac12 (r^2/s) \rbr{ 1 - \tfrac12 (r/s) } }
\le
	\expb{ - \tfrac18 (r^2/s) }.
\]

For $r \ge \sqrt2 s$, we use the inequalities
\[
	\pr{ X_s \ge r } \le \pr{ \Poisson(s) \ge r }
\Qand
	\pr{ X_s \ge r } \ge \pr{ N = 2r } \, \pr{ Y_{2r} \ge r }.
\]
This case is completed by applying (\ref{eq-p0:rp:lde:poi:log}, \ref{eq-p0:rp:lde:poi:rat}), ie \cref{res-p0:rp:lde:poi}.

\smallskip

We now prove \cref{eq-p0:rp:lde:srw:rat}.
For $r \ge \tfrac12 s$, it follows from the fact that $r \mapsto \pr{X_s = r}$ is decreasing and
\[
	\supt{s, \, r \st r \ge s/2}
	\pr{ X_s = r+2 } / \pr{ X_s = r }
<
	1,
\]
which can be verified via a direct calculation involving averaging over $N$ and applying Stirling's approximation; we omit the details.
For $r \le \tfrac12 s$, it suffices to prove the following corresponding result for $(Y_k)_{k\in\mbz_+}$:
	uniformly in $k > 0$ and $r \in [\sqrt k, \tfrac12k]$ with $r \in \mbz$, we have
	\[
		\frac{ \pr{ Y_{2k} \ge 2r } }{ \pr{ Y_{2k} = 2r } }
	\asymp
		\frac{ k }{ r }
	\asymp
		\frac{ \pr{ Y_{2r+1} \ge 2r+1 } } { \pr{ Y_{2r+1} = 2r+1 } };
	\label{eq-p0:rp:lde:proof:srw-suff}
	\nt
	\]
	from this, the original claim follows by averaging over $N$.
Using Stirling's approximation, it is not hard to verify for $r \in [\sqrt k, \tfrac12 k]$ that there exist universal constants $c_1, c_2 \in (0,1)$ such that the following hold:
\begin{gather*}
	c_1 \, \pr{ Y_{2k} = 2r }
\le
	\pr{ Y_{2k} = 2(r + \ceil{k/r}) }
\le
	c_2 \, \pr{ Y_{2k} = 2r };
\\
	c_1 \, \pr{ Y_{2k+1} = 2r+1 }
\le
	\pr{ Y_{2k+1} = 2(r + \ceil{k/r})+1 }
\le
	c_2 \, \pr{ Y_{2k+1} = 2r+1 }.
\end{gather*}
This, together with the fact that both $r \mapsto \pr{ Y_{2k} = 2r }$ and $r \mapsto \pr{ Y_{2k+1} = 2r+1 }$ are decreasing on $[0,k]$, is easily seen to imply \cref{eq-p0:rp:lde:proof:srw-suff}.
\end{Proof}

\section{Simple Random Walk Exit Times Estimates}
\label{sec-p0:gap}

In this section,
we prove some estimates on exit times for SRW on the integers.
These results were used in the spectral gap analysis of \cite{HOt:rcg:abe:geom}.
The following auxiliary lemma is needed.

\begin{lem}
\label{res-p0:gap:costaylor}
	For $\varphi \in [-\tfrac12 , \tfrac12]$,
	we have
	\[
		2(\pi \varphi)^2
	\ge
		1 - \cos(2 \pi \varphi)
	\ge
		\tfrac23 (\pi \varphi)^2.
	\]
\end{lem}


\begin{Proof}[Proof of \cref{res-p0:gap:costaylor}]
Let $\varphi \in [- \tfrac12 , \tfrac12 ]$. Then, using the fact
	that $\log(1-x) \ge - x - x^2$ for $\abs x < 1$,
	that $\sum_1^\infty 1/i^2 = \tfrac16 \pi^2$,
	that $\sum_1^\infty 1/i^4 = \tfrac1{90} \pi^4$
and
	that $\varphi \in [-\tfrac12,\tfrac12]$,
we can calculate directly:
\[
	1
\ge
	\frac{1 - \cos(2\pi\varphi)}{2 (\pi\varphi)^2}
=
	\rbB{ \frac{\sin(\pi\varphi)}{\pi\varphi} }^2
=
	\prod_{\ell=1}^\infty
	\rbB{ 1 - \frac{\varphi^2}{\ell^2} }^2
\ge
	\expbb{ - 2 \sum_{\ell=1}^\infty \rbB{ \frac{\varphi^2}{\ell^2} + \frac{\varphi^4}{\ell^4 }} }
\ge
	0.383
\ge
	\tfrac13.
\qedhere
\]
\end{Proof}

\begin{lem}
\label{res-p0:gap:exitinterval}
	Let $\ell \in \mbn$ and $\tau \cq \inf\bra{ s \ge 0 \mid \abs{Y_s} = \ell }$, where $(Y_s)_{s\ge0}$ is a continuous-time rate-1 SRW on $\mbz$.
	Let $\theta \cq \tfrac12 \pi / \ell$ and $\lambda \cq 1 - \cos\theta$.
	Then, for all $s \ge 0$, we have
	\[
		\pr[0]{ \tau > s }
	\ge
		e^{- \lambda s}
	\ge
		\expb{- \tfrac18 s (\pi/\ell)^2 }.
	\]
\end{lem}

\begin{Proof}[Proof of \cref{res-p0:gap:exitinterval}]
The second inequality follows from \cref{res-p0:gap:costaylor}.

For the first inequality, we first note that
\[
	\mu
:
	x \mapsto \cos(\theta x) \big/ \sumt[\ell]{j=-\ell} \cos(\theta j)
:
	\bra{-\ell, ..., \ell} \to [0,1]
\]
is a distribution satisfying $\mu(\pm \ell) = 0$ and
\[
	(\mu \widehat P)(x) = \mu(x) \, \cos\theta
\Qfor
	x \in J = \bra{ -\ell+1, ...,\ell-1 },
\]
where $\widehat P$ is the transition matrix of discrete-time SRW on $\{-\ell, ...,\ell \} $ with absorption at the boundary.
	Indeed, using $\mu(\pm \ell) = 0$ we have
	\(
		(\mu \widehat P)(x) = \tfrac12 \rbr{ \mu(x+1)+\mu(x-1) } = \mu(x) \cos( \pi/(2\ell) ),
	\)
	where we have used $\cos(a+b) + \cos(a-b) = 2 \cos a \cos b$.
If follows that starting from initial distribution $\mu$ we have
$\mu \widehat P^i(J) = (1- \lambda)^i$,
where $\widehat P^i$ is the matrix $\widehat P$ raised to the power $i$, and so $\mu \widehat P^i(J)$ is the probability of not getting absorbed at the boundary by the $i$-th step when the initial distribution is $\mu$.
It follows that
\[
	\pr[\mu]{ \tau > s }
=
	\sumt[\infty]{i=0}
	\mu \widehat P^i(J) \, \pr{ \Poisson(t) = i }
=
	e^{-\lambda s}.
\]
By considering the continuous-time chain with jump-matrix $\widehat P$, we obtain,
for all $s \ge 0$,
that 
\[
	\pr[0]{ \tau > s }
=
	\maxt{j \in J}
	\pr[j]{ \tau > s },
\]
as can be seen by a simple coupling argument; cf \cite[Example 5.1]{LPW:markov-mixing}.
This concludes the proof.
\end{Proof}

\begin{defn}
\label{def-p0:gap:diric}
	For a transition matrix $P$ and a set $A$, let $\lambda_A$ be the \textit{minimal Dirichlet eigenvalue}, namely the minimal eigenvalue of minus the generator of the chain killed upon exiting $A$, ie of
	\[
		I_A - P_A
	\Qwhere
		(I_A - P_A)(x,y) \cq \one{x,y \in A} \rbb{ \one{x = y} - P(x,y) }.
	\]
	Also, for a set $A$, write $\tau_{A^c}$ for the (first) exit time of this set by the chain.
\end{defn}

\begin{lem}
\label{res-p0:gap:diric}
	Consider a rate-1, continuous-time, reversible Markov chain with transition matrix $P$.
	Let $A$ be a connected set, and let $\lambda_A$ and $\tau_{A^c}$ be as in \cref{def-p0:gap:diric}.
	For all $a \in A$,
	we have
	\[
		- \tfrac1t \log \pr[a]{ \tau_{A^c} > t }
	\to
		\lambda_A
	\quad
		\text{as $\tinf$}.
	\]
\end{lem}

\begin{Proof}[Proof of \cref{res-p0:gap:diric}]
For connected $A$,
by the Perron--Frobenius theorem,
the quasi-stationary distribution of $A$, which we denote by $\mu = (\mu_a)_{a \in A}$, is positive everywhere on $A$. (See \cite[\S 3.6.5]{AF:book} for the definition of quasi-stationarity.)
Since
\(
	\pr[\mu]{ \tau_{A^c} > t }
=
	\sumt{a \in A}
	\mu_a \,
	\pr[a]{ \tau_{A^c} > t },
\)
we have
\[
	\pr[a]{ \tau_{A^c} > t }
\le
	\mu_a^{-1} \, \pr[\mu]{ \tau_{A^c} > t }
=
	\mu_a^{-1} \exp{ - \lambda_A t },
\]
since the exit time starting from the quasi-stationary distribution is exponential with rate $\lambda_A$, as shown in the equation proceeding (3.83) in \cite{AF:book}.
This proves the upper bound, taking $\tinf$.

For the other direction, we claim that there exists a constant $c$, independent of $a$ and $t$, so that
\[
	\mint{a \in A}
	\pr[a]{ \tau_{A^c} > t }
\ge
	c \,
	\maxt{a \in A}
	\pr[a]{ \tau_{A^c} > t + 1 }.
\]
Indeed,
	let $a'$ be an element of $A$ attaining the maximum at time $t+1$.
	Using the connectedness of $A$, for any other $a \in A$ there exists a path from $a'$ to $a$ consisting of states belonging to $A$.
	The probability that the walk traverses this path, and does so in time less than $1$, is at least $c$, for some $c$ independent of $t$.
From this we deduce that
\[
	\mint{a \in A}
	\pr[a]{ \tau_{A^c} > t }
\ge
	c \,
	\pr[\mu]{ \tau_{A^c} > t + 1 }
=
	c \, \expb{ - \lambda_A (t + 1) }.
\]
This proves the lower bound, taking $\tinf$, and hence proves the lemma.
\end{Proof}

\section{Size of Discrete Lattice Ball Estimates}
\label{sec-p0:balls}

We wish to determine the size of the $L_\qq$ balls in $\mbz^k$.
In particular, we desire $R_{k,\qq}$ so that
\[
	\abs{B_{k,\qq}(R_{k,p})} \approx n
\Qwhere
	B_{k,\qq}(R) \cq \brb{ a \in \mbz^k \midb \sumt[k]{1} \abs{a_i}^\qq \le R^\qq }.
\]
This is done by Lemmas \ref{res-p0:balls:size} and \ref{res-p0:balls:R}.
First we need a definition and preliminary lemma.

\begin{defn}
\label{def-p0:balls:R}
	Set $\omega \cq \max\bra{ \logk[2], k/n^{1/(2k)} }$, and choose $R_{k,\qq}$ to be the minimal integer satisfying $\abs{ B_{k,\qq}(R_{k,\qq}) } \ge n e^\omega$.
	Note that $\omega$ satisfies $1 \ll \omega \ll k$ if $k \ll \log n$.
\end{defn}

For $\qq \in [1,\infty)$, write $V_{k,\qq}(R)$ for the (Lebesgue) volume of the $L_\qq$ ball of radius $R$ in $\mbr^k$, ie
\[
	V_{k,\qq}(R) \cq \vol\brb{ x \in \mbr^k \midb \norm{x}_\qq \le R };
\]
also write $V_{k,\qq} \cq V_{k,\qq}(1)$ and note that $V_{k,\qq}(R) = R^k V_{k,\qq}$.
It is known (see \cite{W:Lp-balls}) that
\[
	V_{\ell,\qq} = 2^\ell \Gamma(1/\qq+1)^\ell / \Gamma(\ell/\qq+1).
\label{eq-p0:balls:vol}
\nt
\]
We can use this, along with \cref{res-p0:balls:size:p} below, to well-approximate $\abs{B_{k,\qq}(R)}$ when $\qq \notin \bra{1,\infty}$;
for $\qq = 1$ we directly bound $\abs{B_{k,1}(\cdot)}$, while for $\qq = \infty$ we have an exact expression.

\begin{subequations}
	\label{eq-p0:balls:size}
\begin{subtheorem}{thm}
	\label{res-p0:balls:size}
	
\begin{lem}
\label{res-p0:balls:size:1}
	For $\qq = 1$ and all $R \ge 0$,
	we have
	\[
		2^{k \wedge R} \binomt{\floor{R}}{k} \one{R \ge k}
	\le
		\absb{ B_{k,1}(R) }
	\le
		2^{k \wedge R} \binomt{\floor{R}+k}{k}.
	\label{eq-p0:balls:size:1}
	\nt
	\]
\end{lem}

\begin{lem}
\label{res-p0:balls:size:p}
	For $\qq \in (1,\infty)$ and all $R \ge k^{1+1/\qq}$,
	we have
	\[
		\absb{ B_{k,\qq}(R) }
	=
		V_{k,\qq}(R) \, \rbb{ 1 + \Ohb{ k^{1+1/\qq}/R } }.
	\label{eq-p0:balls:size:p}
	\nt
	\]
\end{lem}

\begin{lem}
\label{res-p0:balls:size:inf}
	For $\qq = \infty$ and all $R \ge 0$,
	we have
	\[
		\absb{ B_{k,\infty}(R) }
	=
		\rbb{2\floor{R} + 1}^k.
	\label{eq-p0:balls:size:inf}
	\nt
	\]
\end{lem}

\end{subtheorem}
\end{subequations}

\begin{Proof}[Proof of \cref{res-p0:balls:size:1}]
Assume $R \in \mbn$.
Observe that
\[
	\absb{ B_{k,1}(R) }
=
	\absb{ \brb{ a \in \mbz^k \mid \sumt[k]{i=1} \abs{a_i} \le R } }.
\]
Moreover, it is a standard combinatorial identity that
\[
	\absb{ \brb{ \alpha \in \mbz_+^k \mid \sumt[k]{i=1} \alpha_i \le R } }
=
	\binomt{R+k}{k}.
\]
The upper and lower bounds will follow easily from this view point, setting $\alpha_i \cq \abs{a_i}$.

For the upper bound, note that $\alpha_i = \abs{\pm a_i}$, and so given the value of $\alpha_i$, there are two choices for $a_i$ if $\alpha_i > 0$, otherwise there is only one (since $0 = -0$). Hence
\[
	\absb{ \brb{ a \in \mbz^k \mid \sumt[k]{i=1} \abs{a_i} \le R } }
\le
	2^{k \wedge R} \absb{ \brb{ a \in \mbz_+^k \mid \sumt[k]{i=1} \alpha_i \le R } }
=
	2^{k \wedge R} \binomt{R+k}{k},
\]
noting that there are at most $k \wedge R$ non-zero coordinates for which a sign can be chosen.

For the lower bound, we get the factor of $2^{k \wedge R}$ by only considering $a \in \mbz^k$ with $\abs{a_i} > 0$ for all $i$, and then setting $\beta_i \cq \alpha_i - 1$. Concretely, for $R \ge k$, we have
\[
&	\absb{ \brb{ a \in \mbz^k \midb \sumt[k]{i=1} \abs{a_i} \le R } }
\ge
	\absb{ \brb{ a \in \mbz^k \midb \sumt[k]{i=1} \abs{a_i} \le R, \: a_i \neq 0 \: \forall \, i = 1, ..., k } }
\\&\qquad
=
	2^{k \wedge R} \,
	\absb{ \brb{ \alpha \in \mbz^k \midb \sumt[k]{i=1} \alpha_i \le R, \: \alpha_i > 0 \: \forall \, i = 1, ..., k } }
\\&\qquad
=
	2^{k \wedge R} \,
	\absb{ \brb{ \beta \in \mbz^k \midb \sumt[k]{i=1} \beta_i \le R - k, \: \beta_i \ge 0 \: \forall \, i = 1, ..., k } }
=
	2^{k \wedge R} \binomt{R}{k}.
\qedhere
\]
\end{Proof}

\begin{Proof}[Proof of \cref{res-p0:balls:size:p}]
For any $R$, writing $\DIAM_\qq$ for the $L_\qq$ diameter (in $\mbr^k$), we have
\[
	B_\qq^k\rbb{ R - \DIAM_\qq [-\tfrac12, \tfrac12)^k }
\subseteq
	B_\qq^k(R)
\subseteq
	B_\qq^k\rbb{ R + \DIAM_\qq [-\tfrac12, \tfrac12)^k }.
\]
Note that
\(
	\DIAM_\qq [-\tfrac12, \tfrac12)^k = k^{1/\qq}.
\)
Hence, for $R$ with $R \ge k^{1+1/\qq}$, we have
\[
	\absb{ B_\qq^k(R) }
=
	\rbb{ 1 + \Oh{k^{1/\qq}/R} }^k
=
	1 + \Ohb{ k^{1+1/\qq}/R }.
\]
Cf \cite[Lemma~2.5]{H:cutoff-cayley-<}, where the case $\qq = 2$ is considered; there, convolutions are employed.
\end{Proof}

\begin{Proof}[Proof of \cref{res-p0:balls:size:inf}]
In the $L_\infty$ norm, the coordinates are independent.
The claim follows.
\end{Proof}

We use this lemma to find an $R_{k,p}$ from \cref{def-p0:balls:R}, which is the minimal integer satisfying $\abs{B_{k,\qq}(R_{k,p})} \ge n e^\omega$.
Recall that $\mfm_{k,\qq} = k^{1/\qq} n^{1/k} / C_\qq$, and that $C_\qq = 2 \, \Gamma(1/\qq + 1) (\qq e)^{1/\qq}$.
The next lemma shows that the difference between $M$ and $\mfm$ is only by subleading order terms.
Also, let $K$ be a constant, assumed to be as large as required, and let $\xi \cq 1 - e^{-K\omega/k}$ when $k \ll \log n$. (As such, we can always replace $1 \pm \xi$ by $e^{\pm\xi}$.)

\begin{subequations}
	\label{eq-p0:balls:R}
\begin{subtheorem}{thm}
	\label{res-p0:balls:R}

\begin{lem}
\label{res-p0:balls:R:1}
	For $k \ll \log n$ and $\qq = 1$,
	we have
	\[
		R_{k,1} \le \ceilb{ \mfm_{k,1} (1 + \xi) }
	\Qand
		\absb{ B_{k,1}\rbb{ \mfm_{k,1} (1 - \xi) } } \ll n.
	\label{eq-p0:balls:R:1}
	\nt
	\]
\end{lem}

\begin{lem}
\label{res-p0:balls:R:p}
	For $k \le \log n/\log\log n$ and all $\qq \in [1,\infty)$,
	we have
	\[
		R_{k,\qq} \le \floorb{ \mfm_{k,\qq} (1 + \xi) }
	\Qand
		\absb{ B_{k,\qq}\rbb{ \mfm_{k,\qq} (1 - \xi) } } \ll n.
	\label{eq-p0:balls:R:p}
	\nt
	\]
\end{lem}

\begin{lem}
\label{res-p0:balls:R:inf}
	For $\qq = \infty$,
	we have
	\[
		R_{k,\infty}
	=
		\ceilb{ \tfrac12 n^{1/k} e^{\omega/k} - \tfrac12 }
	\Qand
		\absb{ B_{k,\infty}\rbb{ \mfm_{k,\infty} (1 - \xi) } } \ll n.
	\label{eq-p0:balls:R:inf}
	\nt
	\]
	Moreover, if $k \ll \log n$ then $R_{k,\infty} \eqsim \mfm_{k,\infty}$.
\end{lem}

\begin{lem}
\label{res-p0:balls:R:logn}
	For all $\lambda \in (0,\infty)$,
	for $k \eqsim \lambda \log n$,
	there exists a function $\omega \gg 1$ and a constant $\alpha$ so that, for all $\eps \in (0,1)$, the minimal integer $M_1$ satisfying $\abs{B_{k,1}(M_1)} \ge n e^\omega$ satisfies
	\[
		R_{k,1} \eqsim \alpha k \eqsim \alpha \lambda \log n
	\Qand
		\absb{ B_{k,1}\rbb{ \alpha k(1 - \eps) } } \ll n.
	\label{eq-p0:balls:R:logn}
	\nt
	\]
	In fact, the result holds for any $1 \ll \omega \ll k$.
\end{lem}

\end{subtheorem}
\end{subequations}

%
%


\begin{Proof}[Proof of \cref{res-p0:balls:R:1}]
\emph{Upper Bound.}
Write $M \cq \ceil{ e^\xi k n^{1/k}/(2e) }$.
Note that $k \ll \log n$, and so $n^{1/k} \gg 1$, and so $M \gg k$.
Then, by \cref{eq-p0:balls:size:1} and Stirling's formula,
we have
\[
	\absb{ B_{k,1}(M) }
\ge
	2^k \binomt{M}{k}
&
\ge
	2^k (M-k)^k/k!
\gtrsim
	k^{-1/2} (1-k/M)^k (2eM/k)^k
\\&
\ge
	k^{-1/2} \expb{-k(2k/M + \xi)} \cdot n.
\]
Take $\xi \cq 2 \omega / k$:
	then $k/M \asymp n^{-1/k} \ll n^{-1/(2k)} \le \xi$
	and $e^{-\xi k} \gg k^{1/2}$.
Hence $\abs{ B_{k,1}(M) } \ge n e^\omega$.

\smallskip

\emph{Lower Bound.}
Set
\(
	M \cq k n^{1/k} e^{-K\omega/k} / (2e).
\)
Using $\binom Nk \le (eN/k)^k$ and \cref{eq-p0:balls:size:1}, we have
\[
	\absb{ B_{k,1}(M) }
\le
	\rbb{ 2e(M/k + 1) }^k
\le
	n e^{-K\omega} \expb{ 6 k / n^{1/k} }
\ll
	n,
\]
using $1 + x \le e^x$ with $x = k/M$, $\binom Nr \le (eN/r)^r$ and $\omega \ge k / n^{1/(2k)} \gg k / n^{1/k}$ as $k \ll \log n$.
\end{Proof}

\begin{Proof}[Proof of \cref{res-p0:balls:R:p}]
\emph{Upper Bound.}
From the formula \cref{eq-p0:balls:vol}, we see that
\[
	R_{k,\qq}
\cq
	n^{1/k} e^{2\omega/k} / V_{k,\qq}^{1/k}
=
	\tfrac12 n^{1/k} e^{2\omega/k} \Gamma(k/\qq+1)^{1/k} / \Gamma(1/\qq+1)
\]
satisfies $V_\qq^k(R_{k,\qq}) = n e^{2\omega}$.
Using Stirling's formula, and the fact that $k \gg 1$, we then deduce that
\[
	R_{k,\qq} \le n^{1/k} k^{1/\qq} e^{\xi}/ C_\qq.
\]
Observe that $k^{1+1/\qq}/R_{k,\qq} \asymp k/n^{1/k} \ll k/n^{1/(2k)} \le \omega$.
Applying \cref{res-p0:balls:size:p} with $R \cq R_{k,\qq}$, which is valid since $k \le \log n/\log\log n$, implying $n^{1/k} \gg k$ and hence $R_{k,\qq} \gg k^{1+1/\qq}$, gives
\[
	\absb{ B_\qq^k(R_{k,\qq}) } / V_\qq^k(R_{k,\qq})
&
=
	1 + \Ohb{ k^{1+1/\qq} / R_{k,\qq} }
=
	\expb{\oh\omega}.
\]
Noting that $V_\qq^k(R_{k,\qq}) = ne^{2\omega}$, we hence deduce that $\abs{B_\qq^k(R_{k,\qq})} \ge n e^\omega$.

\smallskip

\emph{Lower Bound.}
Set
\(
	M \cq k^{1/\qq} n^{1/k} e^{-K\omega/k} / C_\qq.
\)
Then, by \cref{eq-p0:balls:vol} and Stirling, we have
\[
	V_\qq^k(M)
=
	C_k M^k
=
	n e^{-K\omega} k^{k/\qq} \big/ \rbb{ \Gamma(k/\qq+1) (\qq e)^{k/\qq} }
\ll
	n.
\]
Note that $M \gg k^{1+1/\qq}$ since $k \le \log n/\log\log n$, and hence $\abs{ B_\qq^k(M) } \ll n$ by \cref{res-p0:balls:size:p}.
\end{Proof}

\begin{Proof}[Proof of \cref{res-p0:balls:R:inf}]
\emph{Upper Bound.}
This is immediate from \cref{eq-p0:balls:size:inf} and the fact that $n^{1/k} \gg 1$.

\smallskip

\emph{Lower Bound.}
Recall \cref{res-p0:balls:size:inf}.
Observe that
\[
	(2M+1)^k \le n e^{-\nu}
\Quad{if and only if}
	k \log(2M) + k \log\rbb{1 + 1/(2M)} \le \log n - \nu.
\]
Let us set $M \cq \tfrac12 n^{1/k} e^{-K\omega/k}$, for a constant $K$.
Then
\[
	(2M+1)^k \le n e^{-\nu}
\Quad{if and only if}
	\log n - K\omega + k\log\rbb{1 + 1/(2M)} \le \log n - \nu.
\]
Recall that $\omega \ge k/n^{1/(2k)} \gg k/n^{1/k} \asymp k/M$.
Hence, for any constant $K$, we have
\[
	\absb{ B_\infty^k(M) }
\le
	(2M+1)^k
\ll
	n,
\]
by choosing $\nu \gg 1$ but with $\nu = \oh{\omega}$.
Also, $k \ll \log n$, so $\floor{M} \gg 1$.
\end{Proof}

\begin{Proof}[Proof of \cref{res-p0:balls:R:logn}]
We first prove that there exists a strictly increasing function $c : (0,\infty) \to (0,\infty)$ so that, for all $a > 0$, omitting here and below all ceiling signs, we have
\[
	\absb{ B_{k,1}(ak) }
=
	\expb{ k \rbb{ c(a) + \oh1 } }.
\]
By considering the number $i$ of coordinates which equal 0, we have
\(
	\abs{ B_{k,1}(ak) }
=
	\sumt[k]{i=0} A_i,
\)
where
\[
	A_i
\cq
	A_i(k,a)
\cq
	\binomt ki 2^{k-i} \binomt{k-i+ak}{ak}.
\]
Choose $i_* \cq i_*(k,a)$ that maximises $A_i$.
Then
\(
	A_{i_*} \le \abs{ B_{k,1}(ak) } \le (k+1) A_{i^*}.
\)
Observe that
\[
	\frac{A_{i+1}}{A_i}
=
	\frac{(k-i)^2}{2(i+1)(k(1+a)-i)},
\]
and hence one can determine $i_*$ as a function of $k$ and $a$, conclude that $i_*(a,k)/k$ converges as $k \to \infty $ and thus determine $c(a)$ (in terms of the last limit). We omit the details.
Knowing this limit allows us to plug this into the definition of $A_i$ and use Stirling's approximation to get
\[
	A_{i_*} = \expb{ k \rbb{ c(a) + \oh1 } },
\]
for some strictly increasing function $c : (0,\infty) \to (0,\infty)$.
Since $k+1 = e^{\oh{k}}$, the claim follows.

\smallskip

\emph{Upper Bound.}
Since $k \eqsim \lambda \log n$, we have $M_1/k \to c^{-1}(1/\lambda)$ as $\ninf$; set $\alpha \cq c^{-1}(1/\lambda)$.

\smallskip

\emph{Lower Bound.}
It follows from the exponential increase in the size of the $L_1$ ball that\linebreak $\abs{ B_{k,1}( (1 - \eps) \alpha k ) } = \oh{n}$ for all $\eps > 0$, where $M_1 \eqsim \alpha k$ and $\alpha = c^{-1}(1/\lambda)$.
\end{Proof}

\section{Some Further Deferred Proofs}
\label{sec-p0:deferred}

\subsection{Uniformity of Linear Combination of Uniform Random Variables}
\label{sec-p0:deferred:unif-gcd}

\begin{lem}
\label{res-p0:deferred:unif-gcd}
	Let $G$ be Abelian, $k \in \mbn$ and $v \in (\mbz \setminus \bra{0})^k$.
	Draw $Z_1, ..., Z_k \sim^\iid \Unif(G)$.
	Then
	\[
		v \bcdot Z
	=
		\sumt[k]{1} v_i Z_i
	\sim
		\Unif(\mfgcd G)
	\Quad{where}
		\mfgcd \cq \gcd(v_1, ..., v_k, \abs G).
	\]
\end{lem}

\begin{Proof}
Decompose $G$ as $\oplus_1^d \: \mbz_{m_j}$.
Write $\mfgcd_j \cq \gcd(v_1, ..., v_k, m_j)$ for each $j \in [d]$.
Then, for each $i \in [k]$, we may write $Z_i = (\zeta_{i,1}, ..., \zeta_{i,d})$ with $\zeta_{i,j} \sim \Unif(\mbz_{m_j})$ with all the $\zeta_{i,j}$ independent.
Then
\[
	(v \bcdot Z)_j = \sumt[k]{i=1} v_i \zeta_{i,j},
\]
where $(v \bcdot Z)_j$ is the $j$-th component of $v \bcdot Z \in \mbz^d$, and in particular $( (v \bcdot Z)_j )_{j=1}^d$ are independent.
Assuming the $d = 1$ case, the above then shows that $(v \bcdot Z)_j \sim \Unif(\mfgcd_j \mbz_{m_j/\mfgcd_j})$ for each $j$.
Hence it suffices to prove the $d = 1$ case.

\smallskip

We now prove the $d = 1$ case.
Since any $i \in [k]$ with $v_i \equiv 0$ mod $m$ does not contribute to the sum, by passing to a subsequence, we may assume that $v_i \not\equiv 0$ mod $n$ for all $i \in [k]$.

We use induction on $\abs{\mci}$.
Let $U \sim \Unif\bra{1,...,n}$ and set $R \cq mU$ where $m \in \bra{1,...,n}$. Define
\[
	\mfgcd \cq \gcd(m,n) \text{ and } r \cq m/\mfgcd
\Quad{so that}
	R = m U = \mfgcd \cdot (r U).
\]
We then have $\gcd(r,n) = 1$, and so $r U \sim \Unif\bra{1,...,n}$: indeed, for any $x \in \bra{1,...,n}$, we have
\[
	\pr{ r U = x } = \pr{ U = x r^{-1} } = \tfrac1n
\Quad{where}
	\text{$r^{-1}$ is the inverse of $r$ mod $n$}.
\]
Thus we have $R = \mfgcd \cdot (r U) \sim \Unif\bra{\mfgcd,2\mfgcd,...,n}$, since $\mfgcd \wr n$.
This proves the base case $\abs{\mci} = 1$.

Now consider independent $X,Y \sim \Unif\bra{1,...,n}$ and set $R \cq aX + bY$. By pulling out a constant as above, we may assume that $a,b \wr n$. Write $c \cq \gcd(a,b,n)$. Then there exist $r,s \in \bra{1,...,n}$ with
\[
	ar + bs \equiv c \mod n,
\Quad{and hence}
	a(mr) + b(ms) \equiv cm \mod n
\text{ for any }
	m \in \bra{1,...,n}.
\]
Thus $\bra{c,2c,...,n} \subseteq \supp(R)$. By writing $R \cq c(ac^{-1} X + bc^{-1} Y)$, with $c^{-1}$ the inverse $\MOD n$, we see that in fact $\supp(R) = \bra{c,2c,...,n}$. It remains to show that $R$ is uniform on its support.

Pulling out the factor $c$, it is enough to consider $\gcd(a,b,n) = 1$.
For $m \in \bra{0,1,...,n-1}$, set
\[
	\Omega_m \cq \brb{ (x,y) \in [n]^2 \mid ax + by \equiv m \mod n }.
\]
We show that $\abs{\Omega_m}$ is the same for all $m$, and hence deduce that $R$ is uniform on $\bra{c,2c,...,n}$.
Indeed, for every $m$ there exists a pair $(x_m,y_m) \in [n]^2$ so that $ax_m + by_m = m$. If also $(x,y) \in \Omega_m$, then letting $x' = x - x_m$ and $y' = y - y_m$, we see that $(x',y') \in \Omega_0$.
This proves the case $\abs{\mci} = 2$.

Now suppose that $X_1,...,X_L \sim^\iid \Unif\bra{1,...,n}$ and $a_1,...,a_L \in \bra{1,...,n-1}$. By the hypothesis,
\[
	\sumt[L-1]{\ell=1} a_\ell X_\ell \sim c_0 U
\Quad{where}
	U \sim \Unif\bra{1,...,n}
\Qand
	c_0 \cq \gcd(a_1,...,a_{L-1},n).
\]
Now, $X_L$ is independent of this sum, and so the previous case applies to say that
\[
	\sumt[L]{\ell=1} a_\ell X_\ell \sim c U
\Quad{where}
	U \sim \Unif\bra{1,...,n}
\Qand
	c \cq \gcd(c_0,a_L,n) = \gcd(a_1,...,a_L,n).
\]
This completes the induction, and hence proves the claim.
\end{Proof}

\subsection{Decomposition for Product of Upper Triangular Matrices}
\label{sec-p0:heisenberg:matrix-product}

Write $H_{p,d}$ for the set of $d \times d$ unit upper triangular matrices with entries in $\mbz_p$.

\begin{lem}
\label{res-p0:heisenberg:matrix-product}
Let $Z_1,...,Z_k \in H_{p,d}$.
Let
	$\gamma \in [k]^L$ and $\sigma \in \bra{\pm1}^L$.
For $i,j \in [k]$,
set
\[
	C_{i,j}(\gamma, \sigma)
\cq
	\sumt[L]{\ell=0}
	\sumt[\ell-1]{m=0}
	\sigma_m \sigma_\ell \one{\gamma_m = i, \, \gamma_\ell = j }
+	\one{i = j}
	\sumt[L]{\ell=0}
	\one{\gamma_\ell = i, \, \sigma_\ell = -1}.
\]
Set
\(
	M \cq Z_{\gamma_1}^{\sigma_1} \cdots Z_{\gamma_L}^{\sigma_L}.
\)
Then,
for all $a \in [d]$,
we have
\[
	M(a,a) = 1
\Qand
	M(a,a+1) = \sumt[L]{\ell=1} \sigma_{\gamma_\ell} Z_{\gamma_\ell}(a,a+1),
\]
and,
for all $a,b \in [d]$ with $b \ge a+2$,
we have
\[
	M(a,b)
=
	\sumt{\ell \in [L]}
	Z_{\gamma_\ell}(a,b)
+	\sumt{i,j \in [k]}
	C_{i,j}(\gamma, \sigma) Z_i(a,a+1) Z_j(a+1,b)
+	g_{a,b}(\gamma, \sigma; Z_1,...,Z_k),
\]
where $g_{a,b}(\gamma, \sigma; Z_1, ..., Z_k)$ is a polynomial in
\(
	\rbr{ Z_i(x,y) : i \in [k], \, x \in [d-1], \, y > x }.
\)
Further, in this polynomial, each monomial contains the term $Z_i(a,a+1)$ either 0 times or exactly once and no monomial contains a term of the form $Z_i(a,a+1)Z_j(a+1,b)$ for $i,j \in [k]$.
\end{lem}

\begin{Proof}
Given $M_\ell \in H_{p,d}$, we can write $M_\ell = I + N_\ell$ with $N_\ell$ strictly upper triangular.
Consider now $M_1, ..., M_L \in H_{p,d}$; write $M \cq M_1 \cdots M_L$.
From the above expression for $M_\ell^{\sigma_\ell}$ and the fact that $N_\ell$ is strictly upper triangular, the claimed expression for $M(a,a+1)$ is immediate---specifically,
\(
	\rbr{
		\prodt[\ell]{r=1}
		N_{m_r}
	}\rbr{a,a+1}
=
	0
\)
for all $m_1, ..., m_\ell \in [L]$ and $a \in [L-1]$ when $\ell \ge 2$.

Herein we consider the terms above the super-diagonal, ie $(a,b)$ with $b \ge a + 2$.
Observe that
\[
	\prodt[L]{\ell=1}
	\rbr{ I + N_\ell }
=
	\sumt[L]{\ell=0}
	\sumt{m_1 < \cdots < m_\ell}
	\prodt[\ell]{r=1}
	N_{m_r}	
\]
where the indices $m_1, ..., m_\ell$ run over all of $[L]$.
Then, for $(a,b)$ with $b \ge a + 2$, we have
\[
	\rbb{
		\prodt[\ell]{r=1}
		N_{m_r}
	}\rbr{a,b}
&
=
	\sumt{c_0, ..., c_\ell \in [d]}
	\one{c_0 = a, c_\ell = b}
	\prodt[\ell]{r=1}
	N_{m_r}(c_{r-1}, c_r)
\\&
=
	\sumt{a = c_0 < c_1 < \cdots < c_{\ell-1} < c_\ell = b}
	\prodt[\ell]{r=1}
	M_{m_r}(c_{r-1}, c_r),
\]
using the strict upper triangular property of the $N_\ell$.
Similarly, for $(a,b)$ with $b \ge a + 2$, we have
\[
	\rbb{
		N_{m_1} N_{m_2}
	}\rbr{a,b}
=
	M_{m_1}(a,a+1) M_2(a+1,b)
+	\sumt[b-1]{c=a+2}
	M_{m_1}(a,c) M_2(c,b).
\]
Next observe that $N_\ell^d = 0$ as $N_\ell$ is strictly upper triangular.
Hence, for any $\sigma \in \bra{\pm1}$, we have
\[
	M_\ell^{\sigma_\ell}
=
	I
+	\sigma_\ell N_\ell
+	N_\ell^2 \one{\sigma_\ell = -1}
+	\sumt[d]{t=3}
	(-1)^t N_\ell^t \one{\sigma_\ell = -1}.
\]
Recall that $M = M_1 \cdots M_L$.
Then,
for $(a,b)$ with $b \ge a + 2$, we may write
\[
	M(a,b)
&
=
	\sumt{m}
	M(a,b)
+	\sumt{m_1 < m_2}
	M_{m_1}(a,a+1) M_{m_2}(a+1,b)
+	R(a,b)
\\&\qquad
+	\sumt{m}
	M_m(a,a+1) M_m(a+1,b) \one{\sigma_m = -1}
\]
where $R(a,b)$ is a `remainder' polynomial, containing the matrix products of degree 2 and higher \emph{except for} those of the form $M_{m_1}(a,a+1) M_{m_2}(a+1,b)$.
Indeed,
since the sequence $(c_0, ..., c_\ell)$ is strictly increasing,
	each monomial in $R(a,b)$ contains the term $M_m(a,a+1)$ for $m \in [L]$ either 0 times or exactly once
and,
since $\ell \ge 3$ and $c_\ell = b$,
	no monomial in $R(a,b)$ contains a term of the form $M_{m_1}(a,a+1) M_{m_2}(a+1,b)$ for $m_1, m_2 \in [L]$.

Suppose now that $M_\ell = Z_{\gamma_\ell}$ for some
$\gamma \cq (\gamma_\ell)_1^L \in [k]^L$.
By the above analysis,
	the `first order' term $\sumt[L]{\ell=1} Z_{\gamma_\ell}(a,b)$ has the desired form
and
	the `remainder' term has the desired property.
Thus it only remains to check that the `second order' term has the desired form.
Writing $\sumt{m_1 < m_2}$ as $\sumt[L]{m_2=1} \sumt[m_2-1]{m_1=1}$,
the $i \ne j$ case follows from some simple algebra; cf \cite[(\ref{eq-p1:cutoff:3:product}, \ref{eq-p1:cutoff:3:f-Cij}, \ref{eq-p1:cutoff:3:Cij-def})]{HOt:rcg:matrix} for the $d = 3$ case.
The analysis of $C_{i,i}$ is similar (and depends on whether or not inverses are allowed).
\end{Proof}

\subsection{Maximal Local Times for Random Bridges}
\label{sec-p0:deferred:local-times}

A \textit{random bridge} $R = (R_\ell)_{\ell=0}^L$ from $0$ to $r$ of length $L$ is a simple random walk started from $0$ and run for $L$ steps conditioned to have $R_L = 0$.

\begin{lem}[Local Time for Random Bridge]
\label{res-p0:cutoff:undir:local-time:bridge}
	Let $r \in \mbz$ and let $R = (R_\ell)_{\ell=0}^L$ be a random bridge from $0$ to $r$ of length $L$.
	Denote by $\mcl(r,L)$ the maximal local time of $R$.
	Let $C \ge 1$.
	Then
	\[
		\maxt{r}
		\ex{ \mcl(r,L)^C }
	\lesssim
		C^C L^{3C/5},
	\]
	where the implicit constant is absolute.
\end{lem}

We prove this via a reduction from a random bridge to a random walk.
We then use the following lemma on the tails of the local time at $0$ for SRW on $\mbz$.

\begin{lem}[Local Time for Random Walk]
\label{res-p0:cutoff:undir:local-time:walk}
	Let $X = (X_\ell)_{\ell=0}^\infty$ be a SRW on $\mbz$ run for $L$ steps.
	Denote by $\mcl'(0,L)$ the local time of $X$ at $0$ in the first $L$ steps.
	Then, for all $k \ge 0$, we have
	\[
		\pr{ \mcl'(0,L) > k }
	\le
		\expb{ - \tfrac12 k^2 / L }.
	\]
\end{lem}

We now prove \cref{res-p0:cutoff:undir:local-time:bridge} given \cref{res-p0:cutoff:undir:local-time:walk}.
Finally we prove \cref{res-p0:cutoff:undir:local-time:walk}.

\begin{Proof}[Proof of \cref{res-p0:cutoff:undir:local-time:bridge}]
We are estimating the local times for the random bridge $R = (R_\ell)_{\ell=0}^L$, where $R_L$ is given.
We reduce first to bridges with endpoint $0$ and then from bridges to SRWs.

Write $\mcl(x,r,L)$ for the local time at $x$ of a random bridge from $0$ to $r$ of length $L$---this is the local time of $R$ at $x$ conditional on $R_L = r$.
Write $\tau_r$ for the first hitting time of $r$.
Thus
\[
	\mcl(x,r,L) =^d \mcl(0,0,L-\ell)
\Quad{conditional on}
	\tau_r = \ell
\]
and further the local time $\mcl(0,0,L-\ell)$ is independent of $\tau_r$.
Hence
\[
	\pr{ \mcl(x,r,L) > k }
=
	\sumt{\ell}
	\pr{ \mcl(x,r,L) > k \mid \tau_r = \ell } \pr{ \tau_r = \ell }
=
	\sumt{\ell}
	\pr{ \mcl(0,0,L-\ell) > k } \pr{ \tau_r = \ell }.
\]
We have thus reduced to considering the local time for a bridge from $0$ to $0$.

Let $X = (X_\ell)_{\ell=0}^L$ be the usual discrete-time SRW on $\mbz$ with $X_0 = 0$ and run for $L$ steps.
It is not immediate that the map $L' \mapsto \mcl(0,0,L')$ is stochastically increasing, due to the conditioning on the endpoint of the bridge.
Instead, we first convert from a statement about the random bridge $R$ to one about the usual SRW $X$, where the corresponding map is increasing:
	let $\mcl'(0,L)$ denote the local time at $0$ of $X$;
	then $L' \mapsto \mcl(0,L')$ is stochastically increasing.
The random bridge $R$ is simply the random walk $X$ conditioned on $X_L$.
It is standard that $\pr{ X_L = 0 } \asymp 1/L^{1/2}$.
Hence
\[
	\pr{ \mcl(0,0,L-\ell) > k }
&
=
	\pr{ \mcl'(0,L-\ell) > k \mid X_L = 0 }
\\&
\lesssim
	L^{1/2} \, \pr{ \mcl'(0,L-\ell) > k }
\le
	L^{1/2} \, \pr{ \mcl'(0,L) > k }.
\]
Plugging this into the above formula, we find that
\[
	\maxt{x,r}
	\pr{ \mcl(x,r,L) > k }
\lesssim
	L^{1/2} \, \pr{ \mcl'(0,L) > k }.
\]
Denote by $\mcl(r,L)$ the maximal local time for a random bridge from $0$ to $r$.
A union bound then~gives
\[
	\maxt{r}
	\pr{ \mcl(r,L) > k }
\lesssim
	L^{3/2} \, \pr{ \mcl'(0,L) > k }.
\]

We control this latter probability via \cref{res-p0:cutoff:undir:local-time:walk}, which says that
\[
	\pr{ \mcl'(0,L) > k }
\le
	\expb{ - \tfrac12 k^2 / L }.
\]
Take
\(
	k
\cq
	\sqrt{ 2 (C+2) L^{1/2} \log L }.
\)
This gives
\[
	\pr{ \mcl'(0,L) > \textstyle \sqrt{ 2 (C+2) L^{1/2} \log L } }
\le
	L^{-C-2}.
\]
We thus deduce that
\(
	\max_{r}
	\pr{ \mcl(r,L) > \sqrt{ 2 (C+2) L^{1/2} \log L } }
\lesssim
	L^{-C}.
\)
Thus
\[
	\maxt{r}
	\ex{ \mcl(r,L)^C }
\lesssim
	L^C \cdot L^{-C}
+	\rbb{ 2 (C+2) L^{1/2} \log L }^{C/2}
\lesssim
	C^C L^{3C/5}.
\qedhere
\]
\end{Proof}

It remains to prove \cref{res-p0:cutoff:undir:local-time:walk}.

\begin{Proof}[Proof of \cref{res-p0:cutoff:undir:local-time:walk}]
Let $X = (X_\ell)_{\ell=0}^\infty$ be a SRW on $\mbz$.
Denote by $0 = \sigma_0, \sigma_1, \sigma_2, ...$ the return times to $0$ and set $\tau_i \cq \sigma_i - \sigma_{i-1}$, which is the length of the $i$-th excursion.
Then $\tau_1, \tau_2, ...$ are iid.
The SRW on $\mbz$ returns to $0$ infinitely often (almost surely), so $\sigma_i < \infty$ for all $i$ (almost surely) and thus $\tau_1, \tau_2, ...$ are all well-defined (almost surely).
Note that $\ex{\tau_1} = \infty$.
Then, by definition,
\[
	\mcl'(0,L)
=
	\sup\brb{ k \in \mbn_0 \mid \sigma_k \le L } + 1,
\Quad{or equivalently}
	\brb{ \mcl'(0,L) > k }
=
	\brb{ \sumt[k]{i=1} \tau_i \le L }.
\]

We use a standard Chernoff-type bound to control this sum:
\[
	\pr{ \mcl'(0,L) > k }
=
	\pr{ \sumt[k]{i=1} \tau_i \le L }
=
	\pr{ \expb{ - \lambda \sumt[k]{i=1} \tau_i } \ge e^{- \lambda L} }
\le
	\ex{ e^{-\lambda \tau_1} }^k / e^{-\lambda L}.
\]
It is easy to show via standard techniques that
\[
	\ex{ e^{-\lambda \tau_1} }
=
	e^{-\lambda} \rbb{ 1 - \sqrt{ 1 - e^{-2\lambda} } }
\le
	e^{-\sqrt{2\lambda}},
\]
where the equality holds for all $\lambda > 0$ and the inequality for all $\lambda \in (0,1)$.
We then optimise over $\lambda \in (0,1)$: set $\lambda \cq \tfrac12 (k/L)^2$.
Some simple algebra gives the desired tail bound.
	%
\end{Proof}

\subsection{Uniform Random Variables in Nilpotent Groups}
\label{sec-p0:deferred:nil}

\begin{lem}
\label{res-p0:nil:unif-rep}
	For each $\ell \in [L]$,
	let $Y_\ell \sim^\iid \Unif(R_\ell)$.
	Then $Y \cq Y_1 \cdots Y_L \sim \Unif(G)$.
\end{lem}

\begin{Proof}
Let $r_0 \in G$ and consider the event $\bra{Y = r_0}$.

If $r_0 = Y_1 \cdots Y_L$, then $r_1 \cq Y_1^{-1} r_0 = Y_2 \cdots Y_L$.
Clearly the right-hand side is in $G_1$, and so the left-hand side must be too.
Hence $r_0 \equiv Y_1$ mod $G_1$, ie $\pi_1(r_0) = Y_1$.
Since $Y_1 \sim \Unif(R_1)$, the probability of this is $1/\abs{R_1} = 1/\abs{G_0/G_1}$.
Similarly, $r_2 \cq Y_2^{-1} r_1 = Y_3 \cdots Y_L$, and we deduce that $r_2 \equiv Y_2$ mod $G_2$, the probability of which is $1/\abs{R_2} = 1/\abs{G_1/G_2}$.

Iterating this argument,
recalling that the $Y_\ell$ are independent,
we deduce that
\[
	\pr{ Y = r_0 }
=
	\prodt[L]{1} 1/\abs{G_{\ell-1}/G_\ell}
=
	\prodt[L]{1} \abs{G_\ell} / \abs{G_{\ell-1}}
=
	\abs{G_L} / \abs{G_0}
=
	1/\abs G.
\]
Since $r_0 \in G$ was arbitrary, we deduce that $Y \sim \Unif(G)$.
\end{Proof}

This gives the following corollary.

\begin{cor}
\label{res-p0:nil:sample-Z}
	For each $(i, \ell) \in [k] \times [L]$,
	sample $Z_{i,\ell} \sim \Unif(R_\ell)$ independently and set $Z_i \cq Z_{i,1} \cdots Z_{i,L}$.
	Then $Z_1, ..., Z_L \sim^\iid \Unif(G)$.
	Further, $Z_{i,\ell} G_\ell \sim \Unif(Q_\ell)$ independently for each~$(i, \ell)$.
\end{cor}

For the remainder of the section, assume that $Z$ is drawn in this way.
\begin{Proof}
All the independence claims are immediate.
The first claim is immediate from \cref{res-p0:nil:unif-rep}.

For the second claim, we have $Z_{i,\ell} \sim \Unif(R_\ell)$ and $\abs{R_\ell} = \abs{Q_\ell}$.
Now, $x G_\ell = y G_\ell$ if and only if $y^{-1} x G_\ell = G_\ell$.
If $X \sim \Unif(R_\ell)$ and $H \in Q_\ell$, say $H = y G_\ell$ with $y \in R_\ell$, then $y^{-1} X \sim \Unif(R_\ell)$ independently of $y$. So $\pr{ X G_\ell = y G_\ell } = 1/\abs{R_\ell}$.
Hence $X G_\ell \sim \Unif(Q_\ell)$.
\end{Proof}

\subsection{A Bound on the Number of Divisors of an Integer}
\label{sec-p0:deferred:divisor-bound}

In this section, we prove the following number-theoretic result.

\begin{lem}
\label{res-p0:deferred:divisor-bound:i(i|n)}
	For all $\eps > 0$,
	there exists a density-$(1-\eps)$ set $\mba_\eps \subseteq \mbn$ such that,
	for
		all $n \in \mba_\eps$,
		all $m \ge 2$
	and
		all $\lambda > 0$,
	we have
	\[
		\sumt{i \in [m]}
		i \, \one{i \wr n}
	\le
		40 (\lambda \eps)^{-1} m (\log m)^2.
	\]
\end{lem}

\begin{Proof}
Choose $N \in \mbn$ (large) and sample $n \sim \Unif(\bra{1,...,N})$; let $\eps, \lambda \in (0,1)$.
We prove that
\[
	\pr{ \cap_{m \in [1,n]} \brb{
			\sumt{i \in [m]} i \one{i \wr n} \le 20 (\eps \lambda)^{-1} m \logm[1+\lambda]
		} }
\ge
	1 - \eps.
\]
This implies the lemma.
We have $\pr{i \wr n} \le 1/i$ for each $i \in [N]$.
For $i \in [\floor{\log_2 N}]$,
defining
\[
	N_i
\cq
	\sumt{j \in [2^{i-1}, 2^i-1]} \one{ j \wr n },
\Quad{we have}
	\ext{ N_i }
\le
	\sumt[2^i-1]{j=2^{i-1}} 1/j
\le
	1.
\]
By Markov's inequality, for any $\lambda, C > 0$, we then have
\[
	\pr{ N_i \ge C i^{1 + \lambda} } \le 1/(C i^{1+\lambda}).
\]
Using the union bound, this gives
\[
	\prt{ \mce }
\ge
	1 - 2 \lambda^{-1} /C
\Qwhere
	\mce
\cq
	\cap_{i\ge1} \bra{ N_i \le C i^{1+\lambda} }.
\]
Set $r \cq \ceil{\log_2 m}$; then $m \le 2^r$.
On the event $\mce$, we then have
\[
	\sumt{i \in [m]} i \one{i \wr n}
\le
	\sumt[r]{i=1} 2^i N_i
\le
	C \sumt[r]{i=1} 2^i i^{1+\lambda}
\le
	C r^{1+\lambda} 2^{r+1}
\le
	20 C m \logm[1+\lambda].
\]
Now let $\eps \in (0,1)$ and set $C \cq 2 / (\eps \lambda)$, so then $\pr{\mce} \ge 1 - \eps$.
The result follows.
\end{Proof}

Exactly the same argument can be used to show the following result.

\begin{lem}
\label{res-p0:deferred:divisor-bound:(i|n)}
	For all $\eps > 0$,
	there exists a density-$(1 - \eps)$ set $\mba_\eps \subseteq \mbn$ such that,
	for
		all $n \in \mba_\eps$,
		all $m \ge 2$
	and
		all $\lambda > 0$,
	we have
	\[
		\sumt{i \in [m]} \one{i \wr n} \le 10 (\lambda \eps)^{-1} \logm[2+\lambda].
	\]
\end{lem}

\renewcommand{\bibfont}{\sffamily}
\renewcommand{\bibfont}{\sffamily\small}
\printbibliography[heading=bibintoc]

\end{document}